\documentclass[12pt,psfig]{article}
\usepackage{graphics,epsfig,amssymb}

\usepackage{amsmath}
\usepackage{amsthm}
\usepackage{amsfonts}
\usepackage{amssymb,latexsym}
\usepackage[all]{xy}
\usepackage{graphicx}
\usepackage[mathscr]{eucal}
\usepackage{verbatim}
\usepackage{hyperref}
\usepackage{color}

\theoremstyle{plain}
    \newtheorem{thm}{Theorem}[section]
    \newtheorem{prop}{Proposition}[section]
    \newtheorem{lemma}{Lemma}[section]

    \newtheorem{cor}{Corollary}[section]

    \newtheorem{defn}{Definition}[section]

    \newtheorem{rem}{Remark}[section]

\numberwithin{equation}{section}

\usepackage{graphicx}
\usepackage{pict2e}

\begin{document}
\title{Propagation Dynamics for Monotone Evolution Systems without Spatial Translation Invariance}
\author{Taishan Yi\\
School of Mathematics (Zhuhai)\\
Sun Yat-Sen University\\
Zhuhai, Guangdong 519082, China
\and
Xiao-Qiang Zhao\thanks{Corresponding author.
	E-mail: zhao@mun.ca}\\
Department of Mathematics and Statistics\\
Memorial University of Newfoundland\\
St.John's, NL A1C 5S7, Canada}

\date {}
\maketitle

\begin{abstract}
In this paper,  under an abstract setting we establish the existence of spatially inhomogeneous steady states and the asymptotic propagation properties for a large class of monotone
evolution systems without spatial translation invariance.   Then we apply the
developed theory to study traveling  waves  and spatio-temporal propagation patterns for time-delayed nonlocal  equations,  reaction-diffusion
equations in a cylinder,  and asymptotically homogeneous KPP-type equations.
We also obtain the existence of steady state solutions and asymptotic spreading properties of solutions for a time-delayed reaction-diffusion equation subject to the Dirichlet boundary condition.
\end{abstract}


\noindent {\bf AMS(2010) Subject Classification.} 35B40, 37C65, 37L15, 92D25.

\baselineskip=18pt

\section {Introduction}

This paper is devoted to the study of the propagation dynamics for nonlinear evolution equations admitting the comparison principle.
Since the seminal works of  Fisher~\cite{f1937} and KPP~\cite{kpp1937}, there have been extensive investigations on travelling wave solutions and propagation phenomena for various evolution equations
(see, e.g., \cite{fm1977,Shen2010, vvv1994} and references therein). A fundamental feature of propagation problem is the asymptotic spreading speed introduced by Aronson and Weinberger~\cite{aw1978}. Under an abstract setting, Weinberger~\cite{w1982} established the theory of traveling waves and spreading speeds for monotone discrete-time systems with spatial translation invariance. This theory has been greatly developed in~\cite{FYZ2017,fz2014,gpt2014,LYZ,lz2007,lz2010,l1989,w2002,ycw2013,yz2015} for more general {\it monotone} and some {\it non-monotone} semiflows so that it can be applied to a variety of discrete and continuous-time  evolution systems in homogeneous or periodic media.
By using the Harnack inequality up to boundary and the strict positivity of solutions, Berestycki \textit {et al.}~\cite{bhn2005,bhn2010} studied the asymptotic spreading speed for KPP equations in  periodic or  non-periodic spatial domains.  Differently from these two approaches, the spreading speed and asymptotic propagation were obtained in \cite{yc2017}
for the Dirichlet problem of monostable  reaction-diffusion equations on the half line by employing the iterative properties of travelling wave maps. Note that the solution maps of such a Dirichlet problem have no spatial translation invariance and  the Harnack inequality  cannot be extended  to the boundary.

With an increasing interest in  impacts of climate changes (see, e.g., \cite{bbltc2012,k2008,wpcmpbfhb2002}), there have been quite a few works on
traveling waves and asymptotic behavior for evolution  equation models with  a shifting environment, see \cite{bdnz2009, bf2018, br2008,BR2009, BG2019, DWZ2018, flw2016,FangPengZhao2018,hl2015,hz2017,lbsf2014, lwz2018, pl2004, wz2018,WuWangZou2019, ZhangZhao2019} and references therein.
Another class of evolution equations consists of those in locally spatially inhomogeneous media (see \cite{ks2011}).
 We should point out that these  evolution equations admit the comparison principle, but their solution maps no longer possess the spatial translation invariance.  This motivated us to develop the theory of spreading speeds and traveling waves for the monotone semiflows without  spatial translation invariance.  As a starting point, we assume that the given  monotone system has two  limiting
 systems in certain translation sense, and then establish  the existence of steady state solutions and asymptotic propagation properties for monotone semiflows without translation invariance.

 In order to overcome the difficulty induced by the lack of translation invariance,
 we first introduce two limiting systems admitting  the translation invariance under an abstract setting. Then for a special class of  initial functions having compact
 supports, we obtain certain estimates of their orbits under translations for
 the limiting system with the upward convergence property,  and  further carry them to the given system without the translation invariance by comparison arguments (see Section 2). Combining
 these estimates with the asymptotic annihilation property of the other limiting
 system, we are able to characterize the propagation dynamics for the given system.

 The rest of the paper is organized as follows.  In Section 2,  we  present notations and preliminary results.    In order to avoid using traveling wave mappings, we directly establish the links between the system without translation invariance and its  limiting systems.
In Section 3, we prove the existence of fixed points and asymptotic propagation properties for discrete-time semiflows. In Sections 4 and 5,
we extend these results to continuous-time semiflows and  a class of
nonautonomous  evolution systems without translation invariance, respectively.  In Section 6, we apply the developed  theory to two types of  time-delayed nonlocal equations with a shifting habitat, a reaction-diffusion equation in a cylinder,  the Dirichlet problem for a
time-delayed equation on the half line, and a KPP-type equation  in spatially inhomogeneous media.
We expect that our developed theory and methods in this paper may be applied to other monotone evolution systems including cooperative and competitive models
with spatio-temporal heterogeneity.

\section {Preliminaries}
Let $\mathbb{Z}$,
$\mathbb{N}$, $\mathbb{R}$, $\mathbb{R}_+$, $\mathbb{R}^N$, and $\mathbb{R}_+^N$ be the sets of all
integers, nonnegative integers, reals,   nonnegative reals,  N-dimensional real vectors,  and  N-dimensional  nonnegative real vectors,
respectively. We equip  $\mathbb{R}^N$ with the norm $||\xi||_{\mathbb{R}^N}\triangleq
\sqrt{\sum \limits_{n=1}^{N}\xi_n^2}$.  Let $X=BC(\mathbb{R},\mathbb{R}^N)$ be the normed
vector space of all bounded and continuous functions from
$\mathbb{R}$ to $\mathbb{R}^N$ with the norm $||\phi||_{X}\triangleq
\sum \limits_{n=1}^{\infty}2^{-n}\sup\limits_{|x|\leq n}\{||\phi(x)||_{\mathbb{R}^N}\}$. Let $X_+=\{\phi\in
X:\phi(x)\in \mathbb{R}_+^N, \forall x\in \mathbb{R}\}$ and
$X_+^\circ=\{\phi\in X:\phi(x)\in Int(\mathbb{R}_+^N), \forall x\in \mathbb{R}\}$.

For a given compact topological space $M$,  let $C=C(M,X)$ be the
normed vector space of all continuous functions from $M$
into $X$ with the norm $||\varphi||_{C}\triangleq
\sup\limits_{\theta\in M}\{||\varphi(\theta)||_{X}\}$,  $C_+=C(M,X_+)$ and
$C_{+}^{o}=C(M,X_{+}^\circ)$. It follows that $C_+$ is a closed
cone in the normed vector space $C$. Note that $C_+^\circ\neq Int(C_+)$ due to
the non-compactness of the spatial domain $\mathbb{R}$. Also, let
$Y=C(M,\mathbb{R}^N)$ be the normed vector space of all
continuous functions from $M$ into $\mathbb{R}^N$ with the
norm $||\beta||_{Y}\triangleq  \sup\limits_{\theta\in
M} \{||\beta(\theta)||_{\mathbb{R}^N}\}$ and
$Y_{+}=C(M,\mathbb{R}_+^N)$.

For the sake of convenience, we identify an element $\varphi\in C$
with a bounded and continuous function from $M\times
\mathbb{R}$ into $\mathbb{R}^N$. For $a \in \mathbb{R}^N$, $\hat{a}\in
X$ is defined as ${\hat{a}}(x)=a$ for all $x\in \mathbb{R}$.
Similarly, $\hat{\hat{a}}\in C$ is defined as
$\hat{\hat{a}}(\theta)=\hat{a}$ for all $\theta\in M$.
Moreover, for any  $\phi\in X$ and $\beta\in Y$, we define
$\tilde{\phi}\in C$ and $\tilde{\beta}\in C$, respectively,  by
$\tilde{\phi}(\theta,x)=\phi(x)$  and
$\tilde{\beta}(\theta,x)=\beta(\theta)$ for all $(\theta,x)\in
M\times \mathbb{R}$. In the following, we identify
$\hat{a}$ or $\hat{\hat{a}}$ with $a$ for $a\in \mathbb{R}$.
Furthermore, we identify $\phi\in X$ and $\beta\in Y$ with
$\tilde{\phi}\in C$ and $\tilde{\beta}\in C$, respectively. Accordingly, we can regard $X$ and $Y$ as subspaces of $C$.

For any $\xi$, $\eta \in X$,  we write $\xi\geq_X \eta$ if $\xi-\eta
\in X_+$, $\xi >_X \eta$ if $\xi\geq_X \eta$ and $\xi\neq \eta$,
$\xi \gg_X \eta$ if $\xi- \eta \in X_+^\circ$. Similarly, for any $\xi$, $\eta \in \mathbb{R}^N$,  we write $\xi\geq_{\mathbb{R}^N} \eta$ if $\xi-\eta
\in \mathbb{R}^N_+$, $\xi >_{\mathbb{R}^N} \eta$ if $\xi\geq_{\mathbb{R}^N}\eta$ and $\xi\neq \eta$,
$\xi \gg_{\mathbb{R}^N} \eta$ if $\xi- \eta \in Int(\mathbb{R}^N_+)$; for any
$\varphi$, $\psi\in C$, we write $\varphi\geq_C \psi$ if $\varphi-
\psi \in C_+$, $\varphi
>_C \psi$ if $\varphi\geq_C\psi$ and $\varphi \neq \psi$, $\varphi\gg_C
\psi$ if $\varphi- \psi \in C_+^\circ$; for any $\varrho$, $\rho\in Y$,
we write $\varrho\geq_{Y} \rho$ if $\varrho-\rho \in Y_+$, $\varrho
>_{Y} \rho$ if $\varrho\geq_{Y}\rho$
and $\varrho \neq \rho$, $\varrho\gg_{Y} \rho$ if $\varrho-\rho \in
Int(Y_+)$, where $Int(Y_+)=\{\beta\in Y: \beta (\theta)\in Int(\mathbb{R}_+^N)  \mbox{ for
all } \theta\in M\}$. For simplicity, we  write
$\geq$, $>$, $\gg$, and $||\cdot||$, respectively, for $\geq_*$,
$>_*$, $\gg_*$, and $||\cdot||_*$, where $*$ stands  for one of ${\mathbb{R}^N}$, $X$,
$C$, and $Y$.

For any two vectors $s, r\in Int(\mathbb{R}_+^N)$ with $s\geq r$, define $C_r=\{\varphi\in
C:0\leq\varphi\leq r\}$ and $C_{r,s}=\{\varphi\in C:r\leq\varphi\leq
s \}$. For $\varphi\in C_+$, define $C_{\varphi}=\{\psi\in C:
0\le \psi\le \varphi\}$.
We also  define $[\varphi,\psi]_*=\{\xi\in *:\varphi\leq_*\xi\leq_*\psi\}$ and  $[[\varphi,\psi]]_*=\{\xi\in *:\varphi\ll_*\xi\ll_*\psi\}$ for $\varphi,\psi\in *$ with $\varphi\leq_*\psi$, where $*$ stands  for one of ${\mathbb{R}^N}$, $X$,
$C$, and $Y$.

For any given $y\in\mathbb{R}$, define the translation operator
$T_{y}$ by $T_{y}[\varphi](\theta,x)=\varphi(\theta,x-y)$ for all $\varphi\in C$, $\theta\in M, x\in\mathbb{R}$. Let  $Q:  C_+\to C_+$ be a given map.
Throughout this paper, we assume that
\begin{enumerate}
\item [(A1)] $T_{-y}\circ Q [\varphi]\geq Q\circ T_{-y}[\varphi]$ for all $\varphi\in C_{+}$ and
$y\in  \mathbb{R}_+$.

\item [(A2)] For any vector $r \in Int(\mathbb{R}_+^N)$, $Q|_{C_r}:C_r\to C_{+}$ is continuous, and monotone in the sense that
		$Q[\phi]\leq Q[\psi]$ whenever $\phi,\psi\in  C_{r}$ with
		$\phi\leq \psi$.
\end{enumerate}

By virtue of (A1), it is easy to see that   for any $(\theta,x,\varphi)\in M\times \mathbb{R} \times C_{+}$,  $T_{-y}\circ Q \circ T_y[\varphi](\theta,x)$  is nondecreasing in  $y\in\mathbb{R}$.   We introduce two
maps $Q_\pm:C_+\to  L^\infty(M\times \mathbb{R},\mathbb{R}^N)$ by
\begin{equation}\label{Qminus}
Q_-[\varphi](\theta,x):=\lim\limits_{y\to-\infty}T_{-y}\circ Q \circ T_y[\varphi](\theta,x), \forall   \varphi\in C_+, (\theta,x)\in
M\times  \mathbb{R}
\end{equation}
and
\begin{equation}\label{Qplus}
Q_+[\varphi](\theta,x):=\lim\limits_{y\to\infty}T_{-y}\circ Q \circ T_y[\varphi](\theta,x), \forall  \varphi\in C_+, (\theta,x)\in
M\times  \mathbb{R}.
\end{equation}
In view of the definitions of $Q_\pm$ and the assumptions (A1) and (A2),  we have the following observation.

\begin{lemma}\label{lemm2.1} The maps  $Q_\pm:C_+\to L^\infty(M\times \mathbb{R},\mathbb{R}^N)$ admit the following properties:
\begin{enumerate}
\item [(i)] $T_y[Q_\pm[\phi]]=Q_\pm[T_y[\phi]]$ for all $(y,\phi)\in
\mathbb{R}\times C_{+}$;

\item   [(ii)]  $Q_\pm$ is monotone in the sense that
$Q_\pm[\phi](\theta,x)\leq Q_\pm[\psi](\theta,x)$ for all $(\theta,x)\in
M\times  \mathbb{R}$ and $\phi,\psi\in  C_{+} $ with
$\phi\leq \psi$.
\end{enumerate}
\end{lemma}

In the following, we further assume that
\begin{enumerate}
\item [(A3)] $T_{-y}\circ Q^n \circ T_y[\varphi]\to Q_-^n[\varphi]$ in $C$ as $y\to-\infty$
and  $T_{-y}\circ Q^n \circ T_y[\varphi]\to Q_+^n[\varphi]$  in $C$ as $y\to\infty$  for all $\varphi\in C_+$ and $n\in \mathbb{N}$.
\item [(A4)] For any vector  $r\in Int(\mathbb{R}_+^N)$, $Q_\pm|_{C_r}:C_r\to C_{+}$ is continuous with  $Q_+[0]=0$ and $Q_+[r^*]=r^*$ for some $r^*\in Int(Y_+)$.
\end{enumerate}

It easily follows from (A3) and (A4) that $Q_\pm[C_+]\subseteq C_+$, $Q_-\leq Q \leq Q_+$, and $Q_\pm$ satisfies (A2).

Next we present several hypotheses about strong positivity, upward convergence, asymptotic annihilation, and uniform asymptotic annihilation of   maps $Q$ and  $Q_\pm$.

\begin{enumerate}
\item [(SP)]  There exists $N^*>0$ such that $Q^n[\varphi](\theta,x)\in Int(\mathbb{R}_+^N)$ for all  $n\geq N^*$, $(\theta,x)\in M\times (0,\infty)$ and $\varphi\in
C_+\setminus\{0\}$.

\item [(UC)]  There exist $c_-^*,c_+^*\in \mathbb{R}$  such that $c_+^*+c_-^*>0$ and $\lim\limits_{n\rightarrow \infty}
\max\limits_{x\in \mathcal{A}_{\varepsilon,n}^+}
||Q_+^n[\varphi](\cdot,x)-r^*(\cdot)||=0$ for all $\varepsilon\in (0,\frac{c_+^*+c_-^*}{2})$ and  $\varphi\in
C_+\setminus\{0\}$, where $\mathcal{A}_{\varepsilon,n}^+= n{[-c_-^*+\varepsilon,c_+^*-\varepsilon]}$.

\item [(AA)] There exist $\bar{c}_-,\bar{c}_+\in \mathbb{R}$  such that
 $\bar{c}_+ +\bar{c}_->0$ and
$\lim\limits_{n\rightarrow \infty}
\sup\limits_{x\notin \mathcal{A}_{\varepsilon,n}^-}
||Q_+^n[\varphi](\cdot,x)||=0$  for all $\varepsilon>0$ and  all $\varphi\in
C_{r^*}$ having compact supports, where  $\mathcal{A}_{\varepsilon,n}^-= n{[-\bar{c}_--\varepsilon,\bar{c}_++\varepsilon]}$.

\item [(UAA)]  $\lim\limits_{n\rightarrow \infty}
\sup\limits_{(\theta,x)\in M\times\mathbb{R}}
Q_-^n[\varphi](\theta,x)=0$ for all   $\varphi\in
C_+$.
\end{enumerate}

If $Q_+$ satisfies  (UC) and (AA) with  $c_-^*=\bar{c}_-$ and $c_+^*=\bar{c}_+$, then $c_-^*$ and $c_+^*$
are called the leftward and rightward spreading speeds, respectively, for the discrete-time monotone system $\{Q_+^n\}_{n\geq 0}$. For the general results on the existence
of spreading speeds for monotone semiflows, we refer to \cite{w1982,lz2007, LYZ,lz2010,fz2014}.

Let us define a function $h:M\times \mathbb{R}
\rightarrow \mathbb{R}$ by
\[
h(\theta,x)=\left\{
\begin{array}{ll}
  1, &  (\theta,x)\in M \times [-1,1],  \\
x+2, &  (\theta,x)\in M \times [-2, -1),
  \\
2-x, \qquad& (\theta,x)\in M \times  (1,2],
  \\
 0, &  (\theta,x)\in M \times  (\mathbb{R}_+\setminus [-2, 2]).
\end{array}
\right.
\]

\begin{prop} \label{prop2.1}
Assume  that $Q_+$ satisfies (UC).
Then the following statements are valid:
\begin{itemize}
\item [{\rm (i)}] $\lim\limits_{n\rightarrow \infty}
\min\limits_{x\in\mathcal{A}_{\frac{\varepsilon}{3},n}^+}
Q_+^n[\frac{r^*}{16}h](\theta,x)\geq \frac{2r^*(\theta)}{3}$ for all $\varepsilon\in (0,\frac{3c_+^*+3c_-^*}{2})$ and $\theta\in M$.

\item [{\rm (ii)}] For any $\varepsilon\in (0,\frac{3c_+^*+3c_-^*}{2})$, there exists $n_0:=n_0({\varepsilon})\geq \max\{1,\frac{6}{\varepsilon}\}$ such that  $T_{-n c}\circ Q_+^{n}[\frac{r^*}{16}h](\cdot,0)\geq \frac{r^*(\cdot)}{2}$ for all
$c\in [-c_-^*+\frac{\varepsilon}{3},c_+^*-\frac{\varepsilon}{3}]$ and $n\geq n_0$.

\item [{\rm (iii)}] For any $\varepsilon\in (0,\frac{3c_+^*+3c_-^*}{4})$, there exists $y_0=y_0({\varepsilon})>2$ such that $T_{-n c}\circ T_{-y_0} \circ Q^{n}\circ T_{y_0}[\frac{r^*}{16}h]\geq \frac{r^*}{4}h$ for all
$c\in [-c_-^*+\frac{2\varepsilon}{3},c_+^*-\frac{2\varepsilon}{3}]$ and $n\in [n_0,2n_0]\cap\mathbb{N}$,  where $n_0=n_0({\varepsilon})$ is
defined as in statement (ii).
\end{itemize}

\end{prop}

\noindent
{\bf Proof.}  (i) follows from  (UC) and the fact that $\frac{r^*}{16}h\in C_+\setminus \{0\}$.

(ii) By (i), there exists $N_0:=N({\varepsilon})\geq 1$ such that $Q_+^n[\frac{r^*}{16}h](\cdot,x)\geq  \frac{r^*}{2}$ for all $n\geq N_0$ and $-nc_-^*+\frac{n\varepsilon}{3}\leq x\leq nc_+^*-\frac{n\varepsilon}{3}$. Then $n_0:=\max\{N_0,  \frac{6}{\varepsilon}\}$, as required.

(iii) In view of (A3), $r^*\in Int(Y_+)$, and the definition of $Q_+$, we know that  $y_0:=y_0({\varepsilon})>2$ such that
$$ T_{-y_0} \circ Q^{n}\circ T_{y_0}[\frac{r^*}{16}h](\cdot,x)\geq Q_+^n[\frac{r^*}{16}h](\cdot,x)-\frac{r^*}{4}$$  for all $n\in[1,2n_0]$ and $-nc_-^*+\frac{n\varepsilon}{3}\leq x\leq nc_+^*-\frac{n\varepsilon}{3}$. This, together with (ii) and the choices of $h, n_0$, implies that
$$
T_{-n c}\circ T_{-y_0} \circ Q^{n}\circ T_{y_0}[\frac{r^*}{16}h]\geq \frac{r^*}{4}h
$$
for all $c\in [-c_-^*+\frac{2\varepsilon}{3},c_+^*-\frac{2\varepsilon}{3}]$ and $n\in [n_0,2 n_0]\cap\mathbb{N}$. \qed

\

Following \cite{Zhaobook}, we say $Q$  is a subhomogeneous map on $[0,r^*]_C$ if $Q[\kappa \phi]\geq \kappa Q[\phi]$ for all $(\kappa,\phi)\in [0,1]\times [0,r^*]_C$.

\begin{prop} \label{prop2.2}
Assume  that $Q_+$ satisfies (UC), $c_+^*>0$, and   $\varepsilon\in (0,\min\{c_+^*,\frac{c_+^*+c_-^*}{2}\})$, and let  $n_0:=n_0({\varepsilon})$ and  $y_0:=y_0({\varepsilon})$ be  defined as in Proposition \ref{prop2.1}. Then the following statements are valid:
\begin{itemize}

\item [{\rm (i)}]  $T_{-n c}\circ T_{-y_0} \circ Q^{n}\circ T_{y_0}[\frac{r^*}{16}h]\geq \frac{r^*}{4}h$ for all
$c\in [\max\{0,-c_-^*+\frac{2\varepsilon}{3}\},c_+^*-\frac{2\varepsilon}{3}]$ and $n\geq n_0$. If, in addition,   $Q$ is subhomogeneous on $[0,r^*]_C$, then
for any $\delta\in [0,1]$,  there holds $T_{-n c}\circ T_{-y_0} \circ Q^{n}\circ T_{y_0}[\frac{\delta r^*}{16}h]\geq \frac{ \delta r^*}{4}h$ for all
$c\in [\max\{0,-c_-^*+\frac{2\varepsilon}{3}\},c_+^*-\frac{2\varepsilon}{3}]$ and all $n\geq n_0$.

\item [{\rm (ii)}]  If $Q$ satisfies (SP) and is  subhomogeneous on $[0,r^*]_C$, then there exists $\kappa_0=\kappa_0(\varepsilon)\in (0,1]$ such that
for any $\delta\in (0,1)$,  there holds $T_{-n c} \circ Q^{n}[\frac{\delta  r^*}{16}h]\geq \frac{\delta \kappa_0 r^*}{4}h$ for all
$c\in [\max\{\varepsilon,-c_-^*+\varepsilon\},c_+^*-\varepsilon]$ and $n\geq N_0(\varepsilon):=N^*+\max\{n_0,\frac{6y_0}{\varepsilon},\frac{6N^* c_+^*}{\varepsilon}\}$.

\end{itemize}

\end{prop}

\noindent
{\bf Proof.}
(i) Define
\begin{eqnarray*}
&&n^*:=\sup\left\{k\geq n_0: \quad T_{-n c-y_0} \circ Q^{n}\circ T_{y_0}[\frac{r^*}{16}h]\geq \frac{r^*}{4}h, \right.\\
&&\qquad \qquad \quad  \left.\mbox{ for all }
c\in [\max\{0,-c_-^*+\frac{2\varepsilon}{3}\},c_+^*-\frac{2\varepsilon}{3}], n\in [n_0,k]\cap\mathbb{N}\right\}.
\end{eqnarray*}
Clearly, $n^*\geq 2n_0$ due to Proposition~\ref{prop2.1}-(iii) . It suffices to prove $n^*=\infty$. Otherwise, we have $n^*<\infty$.  It follows from (A1), (A2), Proposition~\ref{prop2.1}-(iii) and the choices  of $y_0,n_0,n^*$ that for any $
c\in [\max\{0,-c_-^*+\frac{2\varepsilon}{3}\},c_+^*-\frac{2\varepsilon}{3}]$, we have
\begin{eqnarray*}
&&T_{-(n^*+1) c-y_0} \circ Q^{n^*+1 }\circ T_{y_0}[\frac{r^*}{16}h]
\\
&&=T_{-n_0c-y_0}\circ [T_{-(n^*+1-n_0)c}\circ Q^{n_0}\circ T_{(n^*+1-n_0)c}]\circ T_{-(n^*+1-n_0)c}\circ Q^{n^*+1-n_0}\circ T_{y_0}[\frac{r^*}{16}h]
\\
&&\geq  T_{-n_0c-y_0}\circ Q^{n_0}\circ T_{y_0} [T_{-(n^*+1-n_0)c-y_0}\circ Q^{n^*+1-n_0}\circ T_{y_0}[\frac{r^*}{16}h] ]
\\
&&\geq T_{-n_0c-y_0}\circ Q^{n_0}\circ T_{y_0} [\frac{r^*}{4}h ]\geq \frac{r^*}{4}h,
\end{eqnarray*}
which  contradicts the choice of  $n^*$.
For any $\delta\in [0,1]$, we see from the subhomogeneity  of $Q$ that
$$
T_{-n c}\circ T_{-y_0} \circ Q^{n}\circ T_{y_0}[\frac{\delta r^*}{16}h]
\geq \delta T_{-n c}\circ T_{-y_0} \circ Q^{n}\circ T_{y_0}[\frac{ r^*}{16}h]
\geq \delta \frac{  r^*}{4}h,
$$
where
$c\in [\max\{0,-c_-^*+\frac{2\varepsilon}{3}\},c_+^*-\frac{2\varepsilon}{3}]$ and all $n\geq n_0$.

(ii)  It suffices to consider the case of $\delta=1$ due to the subhomogeneity of $Q$. By (SP),  there exists $\kappa_0:=\kappa_0(\varepsilon)\in (0,1]$ such that $T_{-y_0} \circ Q^{N^* } [\frac{r^*}{16}h] \geq \frac{\kappa_0 r^*}{16}h $. It follows  from (i) that for any $c\in [\max\{\varepsilon,-c_-^*+\varepsilon\},c_+^*-\varepsilon]$ and $n\geq N^*+\max\{n_0,\frac{6y_0}{\varepsilon},\frac{6N^* c_+^*}{\varepsilon}\}$, we have
\begin{eqnarray*}
T_{-nc} \circ Q^{n }[\frac{r^*}{16}h]
&=&T_{-(n-N^*)\times \frac{nc-y_0}{n-N^*}} \circ T_{-y_0} \circ Q^{n-N^* }\circ T_{y_0} [ T_{-y_0} \circ Q^{N^* } [\frac{r^*}{16}h] ]
\\
&\geq &T_{-(n-N^*)\times \frac{nc-y_0}{n-N^*}} \circ T_{-y_0} \circ Q^{n-N^* }\circ T_{y_0}  [\frac{\kappa_0r^*}{16}h]
\\
&\geq &\frac{\kappa_0 r^*}{4}h,
\end{eqnarray*}
where $\frac{nc-y_0}{n-N^*}\in  [\max\{0,-c_-^*+\frac{2\varepsilon}{3}\},c_+^*-\frac{2\varepsilon}{3}]$ is used.
 \qed

\

In the following, we introduce the asymptotic subhomogeneity, (SP), and (UC)   hypothesis to replace the subhomogeneity, (SP), and (UC) of $Q$ in Proposition~\ref{prop2.2}.

\begin{enumerate}
\item [(ASH-UC-SP)] There exist sequences $\{c_{\pm l}^*\}_{l=1}^\infty$ in $\mathbb{R}$, $\{r_l^*\}_{l=1}^\infty$ in $Int(Y_+)$, and  $\{Q_{l}:C_+ \to C_+\}_{l=1}^\infty$ such that for all positive integers $k,l$, there hold
\begin{enumerate}

\item [(i)]   $Q_{l}$ satisfies (A1) and (A2),   $Q_l[0]=0$, and
$Q\geq Q_l$ in $C_{r^*}$  for all  $l\in \mathbb{N}$.

\item [(ii)]    $T_{-y}\circ Q_l \circ T_y[\varphi]\to Q_{l}^+[\varphi]$  in $C$ as $y\to\infty$  for all $\varphi\in C_+$, where
$Q_{l}^+[r^*_l]=r^*_l\leq r^*$.

\item [(iii)]  $Q_l$ is subhomogeneous on $[0,r_l^*]_ C$.

\item [(iv)]   $Q_{l}$  satisfies (SP) and    $Q_{l}^+$  satisfies (UC) with  $c_{-l}^*,c_{+l}^*$ and $r^*_l$.
\end{enumerate}
\end{enumerate}

\begin{prop} \label{prop2.2-cor}
Assume  that $Q$ satisfies (ASH-UC-SP).
Let $c_+^*=\lim\limits_{l\to \infty}c_{+l}^*$ and  $c_-^*=\lim\limits_{l\to \infty}c_{-l}^*$.
If $c_+^*>0$ and $\varepsilon\in (0,\min\{c_+^*,\frac{c_+^*+c_-^*}{2}\})$, then there exist  $y_1=y_1({\varepsilon})>2$, $N_1=N_1({\varepsilon})>1$, and $\kappa_1=\kappa_1(\varepsilon),\kappa_2=\kappa_2(\varepsilon)\in (0,1]$ such that  the following statements are valid:
\begin{itemize}

\item [{\rm (i)}]
For any $\delta\in (0,1)$,  there holds $T_{-n c}\circ T_{-y_1} \circ Q^{n}\circ T_{y_1}[\frac{\delta r^*}{16}h]\geq \frac{ \delta \kappa_2 r^*}{4}h$ for all $c\in [\max\{0,-c_-^*+\frac{2\varepsilon}{3}\},c_+^*-\frac{2\varepsilon}{3}],
n\geq N_1$.

\item [{\rm (ii)}]  For any $\delta\in (0,1)$,  there holds $T_{-n c} \circ Q^{n}[\frac{\delta r^*}{16}h]\geq \frac{\delta \kappa_1  \kappa_2 r^*}{4}h$ for all $ c\in [\max\{\varepsilon,-c_-^*+\varepsilon\},c_+^*-\varepsilon],
n\geq N_1$.
\end{itemize}
\end{prop}

\noindent
{\bf Proof.}
By the choices of $c_\pm^*$, there exists a positive integer $l_0:=l_0(\varepsilon)$ such that  $c_{+l_0}^*>0$, $\frac{\varepsilon}{2}\in (0,\min\{c_{+l_0}^*,\frac{c_{+l_0}^*+c_{-l_0}^*}{2}\})$, $|c_{+l_0}^*-c_+^*|<\frac{\varepsilon}{3}$, and $|c_{-l_0}^*-c_-^*|<\frac{\varepsilon}{3}$.
Applying  Proposition~\ref{prop2.2} to $Q_{l_0}$,  we know that  there exist  $y_1:=y_0(\frac{\varepsilon}{2})>2$, $N_1:=N_0(\frac{\varepsilon}{2})>1$, and $\kappa_1:=\kappa_0(\frac{\varepsilon}{2})\in (0,1)$ such that  for any $\delta\in (0,1)$ and $n\in [N_1,\infty)\cap\mathbb{N}$, there hold
 $$
 T_{-n c}\circ T_{-y_0} \circ (Q_{l_0})^{n}\circ T_{y_0}[\frac{\delta  r_{l_0}^*}{16}h]\geq \frac{ \delta  r_{l_0}^*}{4}h,\, \,  \forall  c\in [\max\{0,-c_{-_{l_0}}^*+\frac{\varepsilon}{3}\},c_{+_{l_0}}^*-\frac{\varepsilon}{3}]
 $$
 and
 $$
 T_{-n c} \circ  (Q_{l_0})^{n}[\frac{\delta r_{l_0}^*}{16}h]\geq \frac{\delta \kappa_1 r_{l_0}^*}{4}h, \, \,  \forall c\in [\max\{\frac{\varepsilon}{2},-c_{-_{l_0}}^*+\frac{\varepsilon}{2}\},c_{+_{l_0}}^*-\frac{\varepsilon}{2}].
 $$
 Choose   $\kappa_2\in (0,1]$ such that ${\kappa_2} r^*(\cdot)\leq r_{l_0}^*(\cdot)\leq r^*(\cdot)$. Then  the desired results  follow from the fact
 that  $Q\geq Q_{l_0}$ and the choices of $l_0$, $y_1$, $N_1$, $\kappa_1$, and $\kappa_2$. \qed

\begin{rem} \label{rem2.1}
		The subhomogeneity assumption
		is used only in the proof of  Propositions~\ref{prop2.2} and~\ref{prop2.2-cor}. Consequently,  if we find some other sufficient conditions for the conclusions of Propositions~\ref{prop2.2} and~\ref{prop2.2-cor}, then all the results in the rest of this paper are still valid.
\end{rem}

It is easy to verify the  following properties for $Q_-$.
\begin{prop} \label{prop2.3}
Assume that $Q_-$ satisfies (UAA).
Then the following statements are valid:
\begin{itemize}
\item [{\rm (i)}]  $Q_-^n[r^*](\theta,x)=Q_-^n[r^*](\theta,0)$ for all $(n,\theta,x)\in \mathbb{N}\times M\times \mathbb{R}$.

\item [{\rm (ii)}]  $Q_-^n[r^*](\theta,x)$ is nonincreasing in $n\in  \mathbb{N}$.

\item [{\rm (iii)}]  $Q_-^n[r^*]\to 0$ in $L^\infty(M\times \mathbb{R},\mathbb{R}^N)$ as $n\to \infty$.

\end{itemize}
\end{prop}

\begin{prop} \label{prop2.4}
Assume  that $Q_-$ satisfies (UAA).  Then  the following statements are valid:
\begin{itemize}
\item [{\rm (i)}]  $Q^n[r^*](\cdot,x)$   is nonincreasing in $n\in  \mathbb{N}$ and nondecreasing in $x\in  \mathbb{R}$.

\item [{\rm (ii)}] If  $Q[r^{**}]\leq r^{**}$ for some $r^{**}\in Int(Y_+)$, then
$\lim\limits_{(n,x)\to(\infty, -\infty)}Q^n[r^{**}](\cdot,x)= 0$ in  $Y$. In particular, $\lim\limits_{(n,x)\to(\infty, -\infty)}Q^n[r^{*}](\cdot,x)= 0$ in  $Y$.
\end{itemize}

\end{prop}

\noindent
{\bf Proof.}  (i) is obvious.

(ii)   It suffices to prove the statement for $r^*$ since the proof for $r^{**}$ is similar.  Let $\varepsilon>0$ be given.  By Proposition~\ref{prop2.3}-(iii),  there exists an integer $n^*>0$ such that  $||Q_-^n[r^*]||_{L^\infty(M\times \mathbb{R},\mathbb{R}^N)}\leq \frac{\varepsilon}{3} $ for all $n\geq n^*$. In particular, $||Q_-^{n^*}[r^*](\cdot,0)||\leq \frac{\varepsilon}{3} $.
It follows from  the definition of $Q_-$ that  there is $z^*>0$ such that $$
||T_{-z}\circ Q^{n^*} \circ T_{z}[r^*](\cdot,0)-Q_-^{n^*}[r^*](\cdot,0)||\leq \frac{\varepsilon}{3} \mbox{ for all } z\in (-\infty,-z^*].
$$
Thus, $||T_{-z}\circ Q^{n^*} \circ T_{z}[r^*](\cdot,0)|| < \varepsilon \mbox{ for all }   z\in (-\infty,-z^*]$, that is, $||Q^{n^*} [r^*](\cdot,z) ||< \varepsilon  \mbox{ for all }  z\in (-\infty,-z^*]$. This, together with statement  (i),  implies that $||Q^{n} [r^*](\cdot,z) ||< \varepsilon  \mbox{ for all }  (n,z)\in [n^*,\infty)\times (-\infty,-z^*]$ with $n\in \mathbb{N}$.
Now statement  (ii) follows from the
arbitrariness of $\varepsilon$.   \qed

\

As a consequence of Proposition~\ref{prop2.4}-(ii), we have the following
observation.
\begin{prop} \label{prop2.5}
Assume that $Q_-$ satisfies (UAA) and there exists $\phi\in C_{r^*}$ such that $Q[\phi]=\phi$.
Then $\phi(\cdot,x)\to 0$ in $Y$ as $x\to -\infty$.
\end{prop}

It is easy to verify the  following properties for $Q_+$.
\begin{prop} \label{prop2.30000}
Assume that $Q_+$ satisfies (UC).
Then the following statements are valid:
\begin{itemize}
\item [{\rm (i)}]  $Q_+^n[\alpha r^*](\theta,x)=Q_+^n[\alpha r^*](\theta,0)$ for all $(n,\alpha,\theta,x)\in \mathbb{N} \times [0,1]\times M \times \mathbb{R}$.

\item [{\rm (ii)}]   For any given $\alpha\in (0,1]$, $Q_+^n[\alpha r^*]\to r^*$ in $L^\infty(M\times \mathbb{R},\mathbb{R}^N)$ as $n\to \infty$.

\end{itemize}
\end{prop}

\begin{prop} \label{prop2.6}
Assume that $Q_+$ satisfies (UC).  If   $\min\{c_-^*,c_+^*\}>0$,  then
$\lim\limits_{x\rightarrow \infty}\lim\limits_{n\rightarrow \infty}
||Q^n[r^*](\cdot, x)-r^*||= 0$, and hence, for any $\varepsilon \in (0,1)$, there exists $x_0:=x_0(\varepsilon)>0$ such that  $Q^n[r^*](\cdot, x)\gg (1-\varepsilon )r^* \mbox{ for all }  (x,n)\in  [x_0,\infty)\times \mathbb{N}$.

\end{prop}

 \noindent
{\bf Proof.}  Define
$$
q(x)=\lim\limits_{n\rightarrow \infty}\Big[\max\{ \alpha\in [0,1]:
Q^n[r^*](\cdot,x)\geq \alpha r^* \}\Big]
\mbox{ for all }x\in \mathbb{R}.
$$
According to  the definition of $q(\cdot) $,  we easily see that $0\leq q(x) \leq 1$  and $q(x)$ is nondecreasing in $x\in \mathbb{R}$.
 By  Proposition~\ref{prop2.2}-(i), we have
 $0< q(x) \leq 1$ for all large  $x$.

To finish the proof, it suffices to prove $\lim\limits_{x\to \infty}q(x)=1$ due to Proposition~\ref{prop2.4}-(i). Otherwise, there is $\alpha^{*}\in (0,1)$ such that  $q(x)\leq \alpha^{*}:=\lim\limits_{x\to \infty}q(x)<1$ for  all $x\in\mathbb{R}$.
By Proposition~\ref{prop2.30000}-(ii),   there exist $\gamma_0\in (0, \min\{\alpha^{*}, \frac{1- \alpha^{*}}{3}\})$ and  $N_0>0$ such that $$
Q_+^{N_0}[(\alpha^{*}-\gamma_0)r^{*}](\cdot,0)\geq (\alpha^{*}+3\gamma_0)r^{*}.
$$ In view of  (A3),   there exists  $z_0>0$ such that $$
T_{-z_0}\circ Q^{N_0} \circ T_{z_0} [(\alpha^{*}-\gamma_0)r^{*}](\cdot,0)\geq Q_+^{N_0}[(\alpha^{*}-\gamma_0)r^{*}](\cdot,0)-\gamma_0r^{*}\geq  (\alpha^{*}+2\gamma_0)r^{*}.
$$
 It follows from Proposition~\ref{prop2.4}-(i) and the choice of $\alpha^{*}$ that  there is $x^*>0$ such that $q(x)\geq \alpha^{*}-\gamma_0$,
 and hence, $Q^n[r^*](\cdot,x)\geq (\alpha^{*}-\gamma_0)r^{*}$ for all  $x\geq x^*$ and $n\in \mathbb{N}$.

 Let us define a function $\xi_d:M\times \mathbb{R}
\rightarrow \mathbb{R}$ by
\[
\xi_d(\theta,x)=\left\{
\begin{array}{ll}
  1, &  (\theta,x)\in M \times [-d,d],  \\
x+1+d, &  (\theta,x)\in M \times [-d-1, -d),
  \\
d+1-x, \qquad& (\theta,x)\in M \times  (d,d+1],
  \\
 0, &  (\theta,x)\in M \times  (\mathbb{R}_+\setminus [-d-1, d+1]).
\end{array}
\right.
\]
It then follows that
$$
\lim\limits_{d\to \infty}T_{-z_0}\circ Q^{N_0} \circ T_{z_0} [(\alpha^{*}-\gamma_0)r^{*}\xi_d](\cdot,0)\geq (\alpha^{*}+2\gamma_0)r^{*}\gg  (\alpha^{*}+\gamma_0)r^{*}
$$
and
$$
T_{-(x^*+d+1)}\circ Q^n[r^*]\geq (\alpha^{*}-\gamma_0)r^{*}\xi_d \mbox{ for all } n\in \mathbb{N}.
$$
Thus, there is $d_0>0$ such that
$$T_{-z}\circ Q^{N_0} \circ T_{z} [T_{-(x^*+d+1)}\circ Q^n[r^*]](\cdot,0)\geq (\alpha^{*}+\gamma_0)r^{*}
$$
for all $z\geq z_0$, $d\geq d_0$, and $n\in \mathbb{N}$.  In particular, by taking $z_1=\max\{z_0,x^*+d_0+1\}$, we have
$$T_{-z_1}\circ Q^{N_0} \circ T_{z_1-(x^*+d_0+1)}\circ Q^n[r^*](\cdot,0)\geq (\alpha^{*}+\gamma_0)r^{*}.$$
This, together with (A1), implies that
\begin{eqnarray*}
T_{-z_1}\circ Q^{N_0+n} [r^*](\cdot,0)&=&T_{-z_1}\circ Q^{N_0} \circ T_{z_1-(x^*+d_0+1)}\circ T_{-z_1+(x^*+d_0+1)}\circ Q^n\circ T_{z_1-(x^*+d_0+1)}[r^*](\cdot,0)\\
&\geq& T_{-z_1}\circ Q^{N_0} \circ T_{z_1-(x^*+d_0+1)}\circ Q^n[r^*](\cdot,0)\\
&\geq& (\alpha^{*}+\gamma_0)r^{*} \mbox{ for all } n\in \mathbb{N}.
\end{eqnarray*}
It follows that $Q^{N_0+n} [r^*](\cdot,x)\geq (\alpha^{*}+\gamma_0)r^{*}$, and hence,  $q(x)\geq \alpha^{*}+\gamma_0 > \alpha^{*}$ for all $n\in \mathbb{N}$ and $x\geq z_1$, which contradicts the choice of $ \alpha^{*}$.\qed

\section {Discrete-time  semiflows}
In this section, we study the upward convergence, asymptotic annihilation,
and the existence of fixed points for discrete-time  maps.

\begin{thm} \label{thm3.1}
Assume that $Q_+$ satisfies (UC) and either $Q$ is a subhomogeneous map with (SP), or  $Q$ satisfies (ASH-UC-SP).
If $c_+^*>0$, then for any $\varepsilon\in (0,\min\{c_+^*,\frac{c_+^*+c_-^*}{2}\})$ and $\varphi\in C_+\setminus\{0\}$, there holds $\lim\limits_{n\rightarrow \infty}
\max\{||Q^n[\varphi](\cdot, x)-r^*||:n\max\{\varepsilon,-c_-^*+\varepsilon\}\leq x \leq n(c_+^*-\varepsilon)\}= 0$.

\end{thm}

\noindent
{\bf Proof.}
For any $\varepsilon\in (0,\min\{c_+^*,\frac{c_+^*+c_-^*}{2}\})$, we define
$$
I^\varepsilon=[\max\{\varepsilon,-c_-^*+\varepsilon\},c_+^*-\varepsilon],
$$
$$
U_-(\varepsilon)=\liminf\limits_{n\rightarrow \infty}\Big[\sup\{\alpha\in \mathbb{R}_+:
Q^n[\varphi](\cdot,x)\geq \alpha r^*  \mbox { for all }x\in n I^\varepsilon\}\Big],$$
and
$$
U_+(\varepsilon)=\limsup\limits_{n\rightarrow \infty}\Big[\inf\{\alpha\in \mathbb{R}_+:
Q^n[\varphi](\cdot,x)\leq \alpha r^*  \mbox { for all }x\in n I^\varepsilon\}\Big].
$$
According to  the definitions of $U_-(\varepsilon) , U_+(\varepsilon) $,  we easily see that $0\leq U_-(\varepsilon) \leq U_+(\varepsilon)$,  $U_-(\varepsilon)$ is non-decreasing in $\varepsilon\in (0,\min\{c_+^*,\frac{c_+^*+c_-^*}{2}\})$, and
$U_+(\varepsilon)$ is non-increasing in $\varepsilon\in (0,\min\{c_+^*,\frac{c_+^*+c_-^*}{2}\})$.
 In view of  (UC), $Q\leq Q_+$  and the definitions of $U_+(\varepsilon) $,  we have  $U_+(\varepsilon) \leq 1$ for all $(0,\min\{c_+^*,\frac{c_+^*+c_-^*}{2}\})$.  By  $\varphi\in C_+\setminus\{0\}$ and (SP), there exists $\delta_0:=\delta_0(\varphi)\in (0,1)$ such that $T_{-3}\circ Q^{N^*}[\varphi]\geq \frac{\delta_0 r^*}{16}h$. It follows from (A1) that
 $$
 Q^n[\varphi]\geq  T_{3}\circ  Q^{n-N^*}\circ T_{-3}\circ Q^{N^*}[\varphi]\geq T_{3}\circ Q^{n-N^*}[ \frac{\delta_0 r^*}{16}h] \mbox{ for all } n\geq N^*.
 $$
 This, together with Proposition~\ref{prop2.2}-(ii) or Proposition~\ref{prop2.2-cor}-(ii), implies that
 $0< U_-(\varepsilon) \leq U_+(\varepsilon)\leq 1$ for all $\varepsilon\in (0,\min\{c_+^*,\frac{c_+^*+c_-^*}{2}\})$.

 To finish the proof, we only need to prove $U_-(\varepsilon)=1$ for all $\varepsilon\in (0,\min\{c_+^*,\frac{c_+^*+c_-^*}{2}\})$. Otherwise, $U_-(\varepsilon_0)<1$ for some $\varepsilon_0\in (0,\min\{c_+^*,\frac{c_+^*+c_-^*}{2}\})$. Thus, $U_-(\varepsilon)<1$ for all $\varepsilon\in (0,\min\{c_+^*,\frac{c_+^*+c_-^*}{2}\})$.
 Due to the
monotonicity of $U_\pm$, we easily see that $U_\pm(\varepsilon)$ are
continuous in $\varepsilon \in [0,\varepsilon_0]$, except possibly
for $\varepsilon$ from a countable set of $[0,\varepsilon_0]$. We
may assume,  without loss of generality, that for some
$\varepsilon_1\in (0,\varepsilon_0)$, $U_-$ is continuous at
$\varepsilon_1$.

By Proposition~\ref{prop2.30000}-(ii),   there exist $\gamma_0\in (0, \min\{U_-(\varepsilon_1), \frac{1- U_-(\varepsilon_1)}{3}\})$ and  $N_0>0$ such that $$
Q_+^{N_0}[U_-(\varepsilon_1) r^{*}](\cdot,0)\geq (U_-(\varepsilon_1)+3\gamma_0)r^{*}.
$$ In view of  (A3),   there exists  $z_0>0$ such that $$
T_{-z_0}\circ Q^{N_0} \circ T_{z_0} [U_-(\varepsilon_1)r^{*}](\cdot,0)\geq Q_+^{N_0}[U_-(\varepsilon_1)r^{*}](\cdot,0)-\gamma_0r^{*}\geq  (U_-(\varepsilon_1)+2\gamma_0)r^{*}.
$$
According to the definition of $U_-(\tau)$, for any $\tau\in
(\varepsilon_1,\varepsilon_0)$, there exist $\theta_0\in M$ and sequences $n_k\rightarrow \infty$ and
$(\theta_k,x_k)\in M\times n_k I^\tau$ such that $\lim\limits_{k\rightarrow \infty}\theta_k=\theta_0$ and
$\lim\limits_{k\rightarrow \infty}Q^{n_k}[\varphi](\theta_k,x_k)-U_-(\tau)r^*(\theta_0)\in \partial(\mathbb{R}^N_+)$. Since
 for any bounded subset $\mathcal{B}$ of $\mathbb{R}$,
$x_k+\mathcal{B}\subseteq(n_k-N_0) I^{\varepsilon_1}$ for all
large $k$, we obtain
$$
\liminf\limits_{k\rightarrow \infty}\Big[
\sup\{\alpha\in \mathbb{R}_+:T_{-x_k}\circ Q^{n_k-N_0}[\varphi](\cdot,x)\geq \alpha r^* \mbox{ for all } x\in \mathcal{B}\}\Big]\in [U_-(\varepsilon_1),U_+(\varepsilon_1)]
$$
 and
$$
\limsup\limits_{k\rightarrow \infty}
\Big[\inf\{\alpha\in \mathbb{R}_+:T_{-x_k}\circ Q^{n_k-N_0}[\varphi](\cdot,x)\leq \alpha r^* \mbox{ for all }  x\in \mathcal{B}\}\Big] \in
[U_-(\varepsilon_1),U_+(\varepsilon_1)].
$$
In other words,
$$
\limsup\limits_{k\rightarrow \infty}\Big[{\rm dist}(T_{-x_k}\circ Q^{n_k-N_0}[\varphi], [U_-(\varepsilon_1)r^*,U_+(\varepsilon_1)r^*]_C)\Big]=0.
$$
 It then follows from (A1) and (A2) that

\begin{eqnarray*}
&&\liminf\limits_{k\to \infty }\min\limits_{\theta\in M}[T_{-x_k}\circ Q^{N_0} \circ T_{x_k} [T_{-x_k}\circ Q^{n_k-N_0}[\varphi]](\theta,0)-(\gamma_0+U_-(\varepsilon_1))r^*(\theta)]
\\
&&\geq \liminf\limits_{k\to \infty }\min\limits_{\theta\in M}[T_{-z_0}\circ Q^{N_0} \circ T_{z_0} [T_{-x_k}\circ Q^{n_k-N_0}[\varphi]](\theta,0)-(\gamma_0+U_-(\varepsilon_1))r^*(\theta)]
\\
\\
&&\geq  \min\limits_{\theta\in M}[ T_{-z_0}\circ Q^{N_0} \circ T_{z_0} [U_-(\varepsilon_1)r^*](\theta,0)-(\gamma_0+U_-(\varepsilon_1))r^*(\theta)]
\\
&&\geq \min\limits_{\theta\in M} [\gamma_0r^*(\theta)]
>0.
\end{eqnarray*}
This implies that
$$
\lim\limits_{k\rightarrow \infty}Q^{n_k}[\varphi](\theta_k,x_k)-(\gamma_0+U_-(\varepsilon_1))r^*(\theta_0)=\lim\limits_{k\rightarrow \infty}\Big[Q^{n_k}[\varphi](\theta_k,x_k)-(\gamma_0+U_-(\varepsilon_1))r^*(\theta_k)\Big]\in \mathbb{R}^N_+.
$$
Thus, by the choices of $\theta_0,n_k$, and $r^*(\cdot)$, there holds
$$
U_-(\tau)\geq U_-(\varepsilon_1)+\gamma_0> U_-(\varepsilon_1).
$$
By
the continuity of $U_-$ at $\varepsilon_1$, letting
$\tau\rightarrow \varepsilon_1$, we then have $U_-(\varepsilon_1)\geq
U_-(\varepsilon_1)+\gamma_0>U_-(\varepsilon_1),$ a
contradiction.
This shows that $U_-(\varepsilon)=U_+(\varepsilon)=r^*$ for all $\varepsilon\in
(0,c_+^*)$, and hence, the conclusion holds true.  \qed

\

By applying the arguments in the proof of Theorem~\ref{thm3.1}, combined with Proposition~\ref{prop2.2}-(i) or Proposition~\ref{prop2.2-cor}-(i), we have the following result.
\begin{cor} \label{cor3.1}
Assume that $Q_+$ satisfies (UC) and either  $Q$ is a subhomogeneous map with (SP), or  $Q$ satisfies (ASH-UC-SP).
If $\min\{c_-^*,c_+^*\}>0$, then for any $\varepsilon\in (0,c_+^*)$ and $\varphi\in C_+\setminus\{0\}$, there holds $\lim\limits_{\alpha \rightarrow \infty}
\Big[\sup\{||Q^n[\varphi](\cdot, x)-r^*||:n\geq \alpha \mbox{ and } \alpha \leq x \leq n(c_+^*-\varepsilon)\}\Big]= 0$.

\end{cor}



\begin{thm} \label{thm3.2}
Assume that $Q_+$ satisfies (AA) and $Q_-$ satisfies (UAA). If $\varphi\in
C_{r^*}$, then the following statements are valid:
\begin{itemize}
\item [{\rm (i)}] If $\varphi
$ has a compact support, then $\lim\limits_{n\rightarrow \infty}
\Big[\sup\{
Q^n[\varphi](\theta,x):(\theta,x)\in M\times n(\mathbb{R}\setminus [-\bar{c}_--\varepsilon,\bar{c}_++\varepsilon])\}\Big]=0$ for all $\varepsilon>0$.

\item [{\rm (ii)}] If $Q[r^{**}]\leq r^{**}$ for some $r^{**}\in Int(Y_+)$, then
$\lim\limits_{n\rightarrow \infty}
\Big[\sup\{
||Q^n[\varphi](\cdot,x)||:x\in  (-\infty,-n\varepsilon]\}\Big]=0$ for all $\varepsilon>0$ and $\varphi\in C_{r^{**}}$. In particular, $\lim\limits_{n\rightarrow \infty}
\Big[\sup\{
||Q^n[\varphi](\cdot,x)||:x\in  (-\infty,-n\varepsilon]\}\Big]=0$ for all $\varepsilon>0$ and $\varphi\in C_{r^*}$.
\end{itemize}
\end{thm}

\noindent
{\bf Proof.}  (i) follows from (AA) and $Q\leq Q_+$.

(ii) Fix $\varepsilon>0$ and $\varphi\in C_{r^{**}}$. By Proposition~\ref{prop2.4}-(ii), it follows that for any  $\gamma>0$, there exists $\rho_0:=\rho_0(\gamma)>0$ such that $Q^n[r^{**}](\cdot,x)\ll \gamma r^{**}$ for all $n>\rho_0$, and $x<-\rho_0$.
Thus, $Q^n[r^{**}](\cdot,x)\ll \gamma r^{**}$ for all $n>\max\{\rho_0,\frac{\rho_0}{\varepsilon}\}$ and $x<-n\varepsilon $.
Since $\varphi\leq r^{**}$ and $\gamma$ is arbitrary, we have
$\lim\limits_{n\rightarrow \infty}\Big
[\sup\{
Q^n[\varphi](\cdot,x):x\in  (-\infty,-n\varepsilon]\}\Big]=0$.
\qed

\

By the arguments  supporting Theorem~\ref{thm3.2}-(ii), we also have the following result.
\begin{cor} \label{cor3.2}
Assume that  $Q_-$ satisfies (UAA).
If $Q[r^{**}]\leq r^{**}$ for some $r^{**}\in Int(Y_+)$, then
$\lim\limits_{\alpha\rightarrow \infty}
\Big[\sup\{
||Q^n[\varphi](\cdot,x)||:x\in  (-\infty,-\alpha] \mbox{ and } n\geq	 \alpha\}\Big]=0$ for all $\varphi\in
C_{r^{**}}$. In particular, $\lim\limits_{\alpha\rightarrow \infty}
\Big[\sup\{
||Q^n[\varphi](\cdot,x)||:x\in  (-\infty,-\alpha] \mbox{ and } n\geq	 \alpha\}\Big]=0$ for all $\varphi\in
C_{r^*}$.
\end{cor}

\medskip
Recall that $\phi$ is a  nontrivial  fixed point  of the map $Q$
 if $\phi\in C_+\setminus \{0\}$  and
$Q[\phi]=\phi$. We say that $\phi(\theta,x)$ connects  $0$ to $r^*$ if
$\phi(\cdot,-\infty):=\lim\limits_{s\to-\infty}\phi(\cdot,s)=0$ and $\phi(\cdot,\infty):=\lim\limits_{s\to\infty}\phi(\cdot,s)=r^*$.
It easily follows  that for any $(\theta,x)\in M\times \mathbb{R}$, $W(\theta,x):=\lim\limits_{n\to \infty}Q^n[r^*](\theta,x)$ is well defined.
Further, we have the following observation.

\begin{lemma}\label{lemm3.1}   If  $W$ is  a  nontrivial  fixed point  of the map $Q$ in $C_{r^*}\setminus \{0\}$ and  $W(\cdot,s)$ is nondecreasing in $s\in \mathbb{R}$, then the following statements are valid:

\begin {itemize}
\item [{\rm (i)}] If   (UAA) holds, then $W(\cdot,-\infty)=0$.

\item [{\rm (ii)}] If  $(-c_-^*,c_+^*)\cap (0,\infty)\neq \emptyset$  and either $Q_+$ satisfies (UC) and $Q$ is a subhomogeneous map with (SP), or  $Q$ satisfies (ASH-UC-SP), then $W(\cdot,\infty)=r^*$.

\item [{\rm (iii)}] If $\{(T_{-z}\circ Q)^n [r^*]:n\in \mathbb{N}\}$ is
 precompact in $C$ for some  $z\in (0,\infty)$, then  $T_{-z}\circ Q$ has a  nontrivial  fixed point $W_z$ in $C_+$ such that $W_z(\cdot,s)$ is nondecreasing in $s\in \mathbb{R}$. If, in addition, $W(\cdot,\infty)=r^*$, then $W_z(\cdot,\infty)=r^*$.


\end {itemize}
\end{lemma}

\noindent
{\bf Proof.} (i) and (ii) follow from   Proposition~\ref{prop2.5} and  Theorem~\ref{thm3.1}, respectively.

(iii) Since $Q[W]=W$, we have $r^*>T_{-z}\circ Q[W]=T_{-z}[W]\geq W>0$.  We can define $W_z:=\lim\limits_{n\to \infty} (T_{-z}\circ Q)^n [r^*]$ due to the compactness and monotonicity. Then $T_{-z}\circ Q[W_z]=W_z$, $W\leq W_z<r^*$, and $W_z(\cdot,x)$ is nondecreasing in $x$.    Moreover, if  $W(\cdot,\infty)=r^*$, then by $W\leq W_z\leq r^*$, we obtain $W_z(\cdot,\infty)=r^*$. \qed

\begin{thm} \label{thm3.3}
Assume that  (UAA) and (UC) hold,  and either  $Q$ is a subhomogeneous map with (SP) or  $Q$ satisfies (ASH-UC-SP). If  $c_-^*>0$  and  $\{Q^n [r^*]:n\in \mathbb{N}\}$ is precompact in $C$,
   then the following statements are valid:

\begin {itemize}
\item [{\rm (i)}] If $c_+^*>0$, then
 $Q$ has a  nontrivial  fixed point $W$ in $C_+$ such that $W(\cdot,-\infty)=0$,
$W(\cdot,\infty)=r^*$, and $W(\cdot,s)$ is nondecreasing in $s\in \mathbb{R}$.

\item [{\rm (ii)}] If $c_+^*\leq 0$ and $\{(T_{x_0}\circ Q)^n [r^*]:n\in \mathbb{N}\}$ is precompact in $C$ for some $x_0\in (-c_+^*,c_-^*)$,  then $Q$ has a  nontrivial  fixed point $W$ in $C_+$ such that $W(\cdot,-\infty)=0$,
$W(\cdot,\infty)=r^*$, and $W(\cdot,s)$ is nondecreasing in $s\in \mathbb{R}$.

\end {itemize}
\end{thm}

\noindent
{\bf Proof.} 
(i) Define $W_n(\cdot,x)=Q^n[r^*](\cdot,x)  \mbox{ for all } (x,n)\in \mathbb{R}\times \mathbb{N}$.  Then for $n\in \mathbb{N}$, $0<W_{n+1}\leq W_n<r^*$ and  $W_n(\cdot,x)$ is nondecreasing in $x$.  By taking $\varepsilon=\frac{\min\{c_-^*,c_+^*\}}{3}$ and applying Proposition~\ref{prop2.2}-(i) or Proposition~\ref{prop2.2-cor}-(i) with $c=0$, we see that there exist $n_0,y_0>0$ such that
$W_n(\cdot,y_0)=T_{-y_0} \circ Q^{n}[r^*](\cdot,0)\geq \frac{r^*}{4}>0$ for all  $n\geq n_0$. Thus,  $W(\theta,x):=\lim\limits_{n\to\infty }W_n(\theta,x)\not\equiv0$. By the compactness of $\{Q^n [r^*]:n\in \mathbb{N}\}$ in $C$,   it follows that $W_n$ tends to $W$ in $C$ and hence, $Q[W]=W$ and $W$ is  nondecreasing.   In view of  Lemma~\ref{lemm3.1}-(i,ii), we have $W(\cdot,\infty)=r^*$ and $W(\cdot,-\infty)=0$.

(ii) By the choice of $x_0$,  it follows that
 $x_0>0$ and $\min\{c_+^*+x_0,c_-^*-x_0\}>0$.
Using (A1) repeatedly, we easily verify  that $(T_{x_0}\circ Q)^n[r^*]\geq T_{n x_0}[ Q^n[r^*]]$ for all $n\in \mathbb{N}$. This, together with the argument similar to that in (i) (after slight modifications), implies that $T_{x_0}\circ Q$ has a  nontrivial  fixed point $W_{-x_0}$ in $C_+$ such that $W_{-x_0}(\cdot,-\infty)=0$,
$W_{-x_0}(\cdot,\infty)=r^*$, and $W_{-x_0}(\cdot,s)$ is nondecreasing in $s\in \mathbb{R}$. Thus, the desired conclusion follows from Lemma~\ref{lemm3.1}-(iii). \qed

\section{Continuous-time semiflows}

In this section, we extend our results on spreading speeds and asymptotic
behavior to a continuous-time semiflow on $C_+$. A
map $Q:\mathbb{R}_{+}\times{C}_+\to{C}_+$ is said to be a
continuous-time semiflow on $C_+$ if  for any vector $r\in Int(\mathbb{R}_+^N)$,
$Q|_{\mathbb{R}_{+}\times C_r}: \mathbb{R}_{+}\times{C}_r\to{C}_+$
is continuous,  $Q_0=Id|_{{C}_+}$,  and $Q_{t}\circ
Q_s=Q_{t+s}$ for all $t,s\in\mathbb{R}_+$, where $Q_t\triangleq
Q(t,\cdot)$ for all $t\in\mathbb{R}_+$.

In the  section, we need the following assumption for some results.

\begin{enumerate}
\item [(SC)]
 For any  $\phi_k\in C_+$ with
$\lim\limits_{k\to \infty}\phi_k=0$ and $\sup\limits_{k\in \mathbb{N}}||\phi_k||_{L^\infty(M\times \mathbb{R}, \mathbb{R})}<\infty$,
there holds
$\lim\limits_{k\to \infty}T_{-y}\circ Q_t \circ T_{y}[\phi_k]=0$ in $C$ uniformly for $ (t,y)\in [0,1] \times \mathbb{R}_+$.
\end{enumerate}

\begin{rem} \label{rem4.1}
It is easy to see that if $\lim\limits_{y\to \infty}T_{-y}\circ Q_t \circ T_{y}[\phi]= Q_t^+[\phi]$ in $C$ for any $(t,\phi)\in \mathbb{R}_+\times C_+$ and $\{Q_t^+\}_{t\in \mathbb{R}_+}$ is  a continuous-time semiflow on $C_+$,  then the property (SC) holds for  $Q_t$.

\end{rem}

\begin{thm}
\label{thm4.1}  Let $\{Q_t\}_{t\geq 0}$ be a continuous-time semiflow on $C_+$ such that each $Q_t$ satisfies  (A1) and (A2).  Assume that there exists $t_0>0$ such that $Q_+$ satisfies (AA) and (UC),   $Q_-$ satisfies (UAA),  and either  $Q_{t_0}$ is a subhomogeneous map with (SP), or  $Q_{t_0}$ satisfies (ASH-UC-SP), where $Q_\pm$ are defined as in \eqref{Qminus} and \eqref{Qplus}  with $Q$ replaced by $Q_{t_0}$.
 Then the following statements are valid:
\begin{itemize}
\item [{\rm (i)}] If $c_+^*>0$, then for any $\varepsilon\in (0,\frac{1}{t_0}\min\{c_+^*,\frac{c_+^*+c_-^*}{2}\})$ and $\varphi\in C_+\setminus\{0\}$, we have $$\lim\limits_{t\rightarrow \infty}
\max\{||Q_t[\varphi](\cdot, x)-r^*||:t\max\{\varepsilon,-\frac{c_-^*}{t_0}+\varepsilon\}\leq x  \leq t(\frac{c_+^*}{t_0}-\varepsilon)\}= 0.$$
\item [{\rm (ii)}] If $\varphi
$ has a compact support and (SC) holds,  then $\lim\limits_{t\rightarrow \infty}
\Big[\sup\{
||Q_t[\varphi](\cdot,x)||:x\in t(\mathbb{R}\setminus [-\frac{\bar{c}_-}{t_0}-\varepsilon,\frac{\bar{c}_+}{t_0}+\varepsilon])\}\Big]=0$ for all $\varepsilon>0$.

\item [{\rm (iii)}]
$\lim\limits_{t\rightarrow \infty}
\Big[\sup\{
||Q_t[\varphi](\cdot,x)||:x\in  (-\infty,-t\varepsilon]\}\Big]=0$ for all $\varepsilon>0$ and $\varphi\in C_{r^*}$.

\item [{\rm (iv)}] If there exists a sequence  of points $\{\phi_k\}_{k \in \mathbb{N}}$ in $Int(Y_+)$ such that  $C_+\subseteq \bigcup\limits_{k\in\mathbb{N}}(\phi_k-C_+)$ and $Q_{t_0}[\phi_k]\leq \phi_k$ for
 all $k\in\mathbb{N}$, then
$\lim\limits_{t\rightarrow \infty}
\Big[\sup\{
||Q_t[\varphi](\cdot,x)||:x\in  (-\infty,-t\varepsilon]\}\Big]=0$ for all $\varepsilon>0$ and $\varphi\in C_+$.

\item [{\rm (v)}]  If $\min\{c_-^*,c_+^*\}>0$, then  for any  $\varepsilon\in (0,\frac{1}{t_0}c_+^*)$ and $\varphi\in C_+\setminus\{0\}$, we have $$\lim\limits_{\alpha \rightarrow \infty}
\Big[\sup\{||Q_t[\varphi](\cdot, x)-r^*||:  t\geq t_0\alpha \mbox{ and } \alpha \leq x \leq t(\frac{c_+^*}{t_0}-\varepsilon)\}\Big]= 0.$$

\end{itemize}
\end{thm}

\textbf {Proof.}
Since  $\{Q_t\}_{t\geq 0}$  is an autonomous semiflow, we assume that  $t_0=1$ in our proof. Otherwise, we consider the autonomous semiflow $\{\hat{Q}_t\}_{t\geq 0}:=\{Q_{t_0t}\}_{t\geq 0}$ instead of
$\{Q_t\}_{t\geq 0}$.

(i)  Without loss of generality, we may assume that $\varphi\leq r^*$.

Given any $\varepsilon\in (0,\min\{c_+^*,\frac{c_+^*+c_-^*}{2}\})$ and $\gamma>0$. We prove statement  (i) by distinguishing two cases.

{\it Case 1. } $\min\{c_+^*,c_-^*\}>0$.

By applying Theorem~\ref{thm3.1} to $Q_1$, we have
\begin{equation}\label{4.1}
\lim\limits_{n\rightarrow \infty}
\max \left \{||Q_n[\varphi](\cdot, x)-r^*||:x\in  \left [n\max\{\frac{\varepsilon}{3},-c_-^*+\frac{\varepsilon}{3}\}, n(c_+^*-\frac{\varepsilon}{3})\right ]\right \}= 0.
\end{equation}
Since $Q_t[r^*](\cdot, x)$ is nonincreasing in $t\in \mathbb{R}_+$ and nondecreasing $x\in \mathbb{R}$, it follows from  Proposition~\ref{prop2.6} that there exists
$y_0:=y_0(\gamma)>0$ such that $||Q_t[r^*](\cdot,y)-r^*||<\frac{\gamma}{3}$ for all $t\in \mathbb{R}_+$ and $y\geq y_0$.

By the uniform continuity of $Q$ at $r^*$ for $t\in[0,1]$, there
exist $\delta=\delta(\gamma)>0$ and $d=d(\gamma)>0$ such
that if $\psi\in C_{r^{*}}$ with $||\psi(\cdot,x)-r^*||<\delta$ for all
$x\in [-d,d]$, then
$||Q_{t}\circ T_{y_0}[\psi](\cdot,y_0)-Q_t[r^*](\cdot,y_0)||<\frac{\gamma}{3}$ for all $t\in[0,1]$.

It follows from~\eqref{4.1} that there is an integer $n_1>0$
such that $||Q_{n}[\varphi] (\cdot,x)-r^*||<\delta$ for all $x\in  \left [n\max\{\frac{\varepsilon}{3},-c_-^*+\frac{\varepsilon}{3}\}, n(c_+^*-\frac{\varepsilon}{3})\right ]$
and  $n>n_1$. Let $n_2=\max\{n_1,\frac{3d}{\varepsilon}\}$. Then for any $n>n_2$ and  $y\in  \left [n\max\{\frac{2\varepsilon}{3},-c_-^*+\frac{2\varepsilon}{3}\}, n(c_+^*-\frac{2\varepsilon}{3})\right ]$, we have
$||T_{-y}\circ Q_{n}[\varphi](\cdot,x)-r^*||<\delta$ for all $x\in
[-d,d]  $. Since $Q_t$ satisfies
$(A1)$ for any $t\geq 0$, it follows that for any $n>\max\{n_2,\frac{3 y_0}{2\varepsilon}\}$, $y\in  \left [n\max\{\frac{2\varepsilon}{3},-c_-^*+\frac{2\varepsilon}{3}\}, n(c_+^*-\frac{2\varepsilon}{3})\right ]$, and $t\in [0,1]$, there holds
\begin{eqnarray*}
&&||Q_{t+n}[\varphi](\cdot,y)-r^*||\\
&&\leq ||T_{-y_0} \circ Q_{t}\circ T_{y_0} \circ T_{-y} \circ Q_{n}[\varphi](
\cdot,0)-r^*||\\
&&\leq ||Q_{t}\circ T_{y_0} [ T_{-y} \circ Q_{n}[\varphi]](
\cdot,y_0)- Q_{t}[r^*](\cdot,y_0)||+||Q_{t}[r^*](\cdot,y_0)-r^*||\\
&&<\frac{\gamma}{3}+\frac{\gamma}{3}<\gamma.
\end{eqnarray*}
Thus,
$||Q_{t}[\varphi](\cdot,y)-r^*||<\gamma$ for all $t>1+\max\{n_2,\frac{3y_0}{2\varepsilon},\frac{3c_+^*}{\varepsilon}\}
$ and $y\in  \left [t\max\{\varepsilon,-c_-^*+\varepsilon\}, t(c_+^*-\varepsilon)\right ]$.
Hence, (i) follows from the arbitrariness of $\gamma$.




{\it Case 2. }  $c_-^*\leq 0$ and  $\min\{c_+^*+c_-^*,c_+^*\}>0$.

Take $c\in (-c_-^*,\min\{c_+^*-\varepsilon,-c_-^*+\frac{\varepsilon}{6}\})$.  Then $0\leq -c_-^*<c<c_+^*$. Note that
\begin{eqnarray*}
T_{-ct}\circ Q_t\circ T_{-cs}\circ Q_s &\leq& T_{-ct}\circ [T_{-cs} \circ  Q_t\circ T_{cs}]\circ T_{-cs}\circ Q_s\\
&=& T_{-c(t+s)} \circ  Q_{t+s} \mbox{ for all } t,s\in \mathbb{R_+}.
\end{eqnarray*}
It then follows that both $T_{-c(t+n)} \circ  Q_{t+n}[r^*]\geq T_{-ct} \circ  Q_{t}\circ T_{-cn} \circ  Q_{n}[r^*]$ and $T_{-cn} \circ  Q_{n}[r^*]\geq (T_{-c} \circ  Q_{1})^n[r^*]$ hold true for all $t\in[0,1]$ and $n\in \mathbb{N}$.

In the following, we apply  the arguments for  Case 1  to $T_{-ct}\circ Q_t$ with some necessary modifications  since $T_{-ct}\circ Q_t$ does not admit the semigroup property.

By Proposition~\ref{prop2.6}, as applied to $T_{-c}\circ Q_1$, we know that there exists
$y_0:=y_0(\gamma)>0$ such that $||(T_{-c}\circ Q_1)^n[r^*](\cdot,y)-r^*||<\frac{\gamma}{3}$ for all $n\in \mathbb{N}$ and $y\geq y_0$.  It then follows that
$$
||T_{-c}\circ Q_1[r^*](\cdot,y)-r^*||<\frac{\gamma}{3} \mbox{ for all } y\geq y_0.
$$
This,  together with the fact that
$$
T_{-c} \circ  Q_{1}[r^*]= T_{-c(1-t)} \circ T_{-ct} \circ  Q_{t} \circ  Q_{1-t}[r^*]\leq T_{-c(1-t)} \circ T_{-ct} \circ  Q_{t}[r^*]   \mbox{ for all }  t\in [0,1],
$$
implies that  $||T_{-ct} \circ  Q_{t}[r^*](\cdot,y)-r^*||<\frac{\gamma}{3}  \mbox{ for all } y\in [y_0+(1-t)c,\infty)$ and $t\in [0,1]$.  In particular,  $||T_{-ct} \circ  Q_{t}[r^*](\cdot,y_1)-r^*||<\frac{\gamma}{3}   \mbox{ for all }  t\in [0,1]$, where $y_1=y_0+c$.

In view of the uniform continuity of $Q$ at $r^*$  for $t\in[0,1]$, there
exist $\delta=\delta(\gamma)>0$ and $d=d(\gamma)>0$  \mbox{ for all }  such
that if $\psi\in C_{r^{*}}$ with $||\psi(\cdot,x)-r^*||<\delta$ for all
$x\in [-d,d] $, then
$$||T_{-ct}\circ Q_{t}\circ T_{y_1}[\psi](\cdot,y_1)-T_{-ct}\circ Q_t[r^*](\cdot,y_1)||<\frac{\gamma}{3}   \mbox{ for all }  t\in[0,1].
$$
By applying Theorem~\ref{thm3.1} to $T_{-c}\circ Q_1$, we have
\begin{equation}\label{4.2}
\lim\limits_{n\rightarrow \infty}
\max \left \{||(T_{-c}\circ Q_1)^n[\varphi](\cdot, x)-r^*||:x\in  \left [n\max\{\frac{\varepsilon}{3},-c_-^*-c+\frac{\varepsilon}{3}\}, n(c_+^*-c-\frac{\varepsilon}{3})\right ]\right \}= 0.
\end{equation}
It follows from~\eqref{4.2} that there is an integer $n_1>0$
such that $||(T_{-c}\circ Q_{1})^n[\varphi] (\cdot,x)-r^*||<\delta,$ and hence,
$$
||T_{-nc}\circ Q_{n}[\varphi] (\cdot,x)-r^*||<\delta$$ for all  $x\in  \left [n\max\{\frac{\varepsilon}{3},-c_-^*-c+\frac{\varepsilon}{3}\},\, n(c_+^*-c-\frac{\varepsilon}{3})\right ], \, n>n_1.
$

Let $n_2=\max\{n_1,\frac{3d}{\varepsilon}\}$. Then for any $n>n_2$ and  $y\in  \left [n\max\{\frac{2\varepsilon}{3},-c_-^*-c+\frac{2\varepsilon}{3}\}, n(c_+^*-c-\frac{2\varepsilon}{3})\right ]$, we have
$||T_{-y}\circ T_{-nc}\circ Q_{n}[\varphi](\cdot,x)-r^*||<\delta$ for all $x\in
[-d,d]  $. Since $Q_t$ satisfies $(A1)$ for any $t\geq 0$, it follows that for any $n>\max\{n_2,\frac{2y_1}{\varepsilon}\}$, $y\in  \left [n\max\{\frac{2\varepsilon}{3},-c_-^*-c+\frac{2\varepsilon}{3}\}, n(c_+^*-c-\frac{2\varepsilon}{3})\right ]$, and $t\in [0,1]$,
\begin{eqnarray*}
&&||T_{-c(t+n)}Q_{t+n}[\varphi](\cdot,y)-r^*||\\
&&\leq ||T_{-y} \circ T_{-ct}\circ Q_{t}\circ T_{y} \circ T_{-y}\circ T_{-cn} \circ Q_{n}[\varphi](
\cdot,0)-r^*||\\
&&\leq ||T_{-y_1}\circ T_{-ct} \circ Q_{t}\circ T_{y_1} \circ T_{-y}\circ T_{-cn} \circ Q_{n}[\varphi](
\cdot,0)-r^*||\\
&&\leq ||T_{-ct} \circ Q_{t}\circ T_{y_1}[ T_{-y} \circ T_{-cn} \circ Q_{n}[\varphi]](
\cdot,y_1)- T_{-ct} \circ Q_{t}[r^*](\cdot,y_1)||+||T_{-ct} \circ Q_{t}[r^*](\cdot,y_1)-r^*||\\
&&<\frac{\gamma}{3}+\frac{\gamma}{3}<\gamma.
\end{eqnarray*}
Thus, we have
$$||T_{-ct} \circ Q_{t}[\varphi](\cdot,y)-r^*||<\gamma$$ for all $t>1+\max\{n_2,\frac{2y_1}{\varepsilon},\frac{6(c_+^*-c_-^*)}{\varepsilon}\}
$ and $y\in  \left [t\max\{\frac{5\varepsilon}{6},-c_-^*-c+\frac{5\varepsilon}{6}\}, t(c_+^*-c-\frac{5\varepsilon}{6})\right ]$.
In other words,
$$|| Q_{t}[\varphi](\cdot,y)-r^*||<\gamma$$ for all $t>1+\max\{n_2,\frac{2y_1}{\varepsilon},\frac{6(c_+^*-c_-^*)}{\varepsilon}\}
$ and $y\in  \left [t\max\{\frac{5\varepsilon}{6}+c,-c_-^*+\frac{5\varepsilon}{6}\}, t(c_+^*-\frac{5\varepsilon}{6})\right ]$.
This, together with the choices of $c,\gamma$, implies that
$$\lim\limits_{n\rightarrow \infty}
\max\{||Q_t[\varphi](\cdot, x)-r^*||:t\max\{\varepsilon,-c_-^*+\varepsilon\}\leq x \leq t(c_+^*-\varepsilon)\}= 0.$$

(ii) Given any $\varepsilon>0$, $\gamma>0$ and $r^{**}>0$.   It follows from (A1) and (SC) that
there
exist $\delta=\delta(\gamma)>0$ and $d=d(\gamma)>0$ such
that if $\psi\in C_{r^{**}}$ with $||\psi(\cdot,x)||<\delta$ for all
$x\in [-d,d] $,  then
$||T_{-y}\circ Q_{t}\circ T_{y}[\psi](\cdot,0)||<\gamma$ for all $t\in[0,1]$ and $y\in \mathbb{R}$.

In view of  Theorem~\ref{thm3.2}-(i),  we obtain
\begin{equation}\label{4.3}
\lim\limits_{n\rightarrow \infty}
\Big[\sup\{
Q_n[\varphi](\theta,x):(\theta,x)\in M\times n(\mathbb{R}\setminus [-\bar{c}_--\frac{\varepsilon}{3},\bar{c}_++\frac{\varepsilon}{3}])\}\Big]=0.
\end{equation}
It follows from~\eqref{4.3} that there is an integer $n_1>0$
such that $||Q_{n}[\varphi] (\cdot,x)||<\delta$ for all $x\in  n(\mathbb{R}\setminus [-\bar{c}_--\frac{\varepsilon}{3},\bar{c}_++\frac{\varepsilon}{3}])$
and  $n>n_1$. Let $n_2=\max\{n_1,\frac{3d}{\varepsilon}\}$. Then for any $n>n_2$ and  $y\in  n(\mathbb{R}\setminus [-\bar{c}_--\frac{2\varepsilon}{3},\bar{c}_++\frac{2\varepsilon}{3}])$, we have
$||T_{-y}\circ Q_{n}[\varphi](\cdot,x)||<\delta$ for all $x\in
[-d,d]  $. According to the previous discussions,  we know that for any $n>n_2$, $y\in  n(\mathbb{R}\setminus [-\bar{c}_--\frac{2\varepsilon}{3},\bar{c}_++\frac{2\varepsilon}{3}])$, and $t\in [0,1]$,
\begin{eqnarray*}
||Q_{t+n}[\varphi](\cdot,y)||&=&||T_{-y} \circ Q_{t}\circ T_{y} \circ T_{-y} \circ Q_{n}[\varphi](\cdot,0)||\\
&=&||T_{-y} \circ Q_{t}\circ T_{y} [ T_{-y} \circ Q_{n}[\varphi]](
\cdot,0)||<\gamma.
\end{eqnarray*}
 In particular,
$||Q_{t}[\varphi](\cdot,y)||<\gamma$ for all $t>1+\max\{n_2,\frac{3 (|\bar{c}_+|+|\bar{c}_-|)}{\varepsilon}\}
$ and $y\in t(\mathbb{R}\setminus [-\bar{c}_--\varepsilon,\bar{c}_++\varepsilon])$.
Thus, (ii) follows from the arbitrariness of $\gamma$.

(iii) Let $n_t$ be the integer part of $t$ and $r_t$ be the  nonnegative fraction part of $t$.  Then $Q_t[\varphi]=Q_{n_t}[Q_{r_t}[\varphi]]\leq Q_{n_t}[Q_{r_t}[r^*]]\leq Q_{n_t}[r^*]$. This, together with Theorem~\ref{thm3.2}-(ii), implies that
$\lim\limits_{t\rightarrow \infty}
\Big[\sup\{
||Q_t[\varphi](\cdot,x)||:x\in  (-\infty,-n_t\varepsilon]\}\Big]=0$. It follows that $\lim\limits_{t\rightarrow \infty}
\Big[\sup\{
||Q_t[\varphi](\cdot,x)||:x\in  (-\infty,-t\varepsilon]\}\Big]=0$.

(iv) follows from the proof similar to that of (iii) with $r^*$  replaced by $\phi_k$ .

(v)
By using Corollary~\ref{cor3.1} and  the arguments similar to those in  the proof  of  statement (i), we can  obtain  (v).   \qed

\

In the rest of this section, we consider the nontrivial equilibrium points for the
continuous-time semiflow $\{Q_t\}^{\infty}_{t=0}$ on $C_+$.
We say that $W$ is an equilibrium point of $\{Q_t\}^{\infty}_{t=0}$ if
$W: M\times\mathbb{R}\to \mathbb{R}_+$ is a bounded and
continuous function and
$Q_t[W](\theta,x)=W(\theta,x)$ for all $(\theta,x)\in
M\times \mathbb{R}$ and $t\in \mathbb{R}_+$,  and  that $W$  connects
$0$ to $r^*$ if $W(\cdot,-\infty)=0$ and
$W(\cdot,\infty)=r^*$.

\begin{thm} \label{thm4.2} Suppose that $Q:=Q_{t_0}$ and its associated $Q_\pm$ satisfy all conditions in Theorem~\ref{thm3.3} for some $t_0>0$. If $ c_-^{*}>0 $ and $Q_t[r^*]\leq r^*$ for all $t\in \mathbb{R}_+$, then $\{Q_t\}_{t\in\mathbb{R}_+}$ has an  equilibrium point $W$ connecting
$0$ to $r^*$.
\end{thm}

\noindent
\textbf {Proof.}  Since $ c_-^{*}> 0$ and  $ c_-^{*}+c_+^*>0$,  it is easy to see  that $\min\{c_-^{*},c_+^*\}>0$ or $c_+^*\leq 0$.  By Theorem~\ref{thm3.3},  there exists a nondecreasing function $W\in C_+^\circ$ such that $W(\cdot,-\infty)=0$, $W(\cdot,\infty)=r^*$, $W=\lim\limits_{n\to \infty}Q_{nt_0}[r^*]$, and $Q_{nt_0}[W]=W$.  Since  $Q_t[r^*]$ is nonincreasing in $t\in \mathbb{R}_+$ due to the fact that $Q_t[r^*]\leq r^*$ for all $t\in \mathbb{R}_+$, we have $W=\lim\limits_{t\to \infty}Q_t[r^*]$.  Thus, $Q_t[W]=W$ for all $t\in \mathbb{R}_+$.  \qed

\begin{cor} \label{cor4.1} Suppose that
$Q:=Q_{t_0}$ and its associated $Q_\pm$  satisfy all conditions in
Theorem~\ref{thm3.3} for some $t_0>0$.  Let  $Q_t[r^*]\leq r^*$ for all $t\in \mathbb{R}_+$ and $\min\{ c_-^{*}, c_+^{*}\}>0 $. If  $\{Q_t[\psi]:
t\geq t_0\}$ is precompact  in $C$ and $Q_t$ has at least one equilibrium point in $\mathcal{K}$ for any closed, convex, and positively invariant set $\mathcal{K}$ of $Q_t$,  then $\{Q_t\}_{t\in\mathbb{R}_+}$ has a unique equilibrium point $W$ connecting
$0$ to $r^*$ if and only if  for any $\varepsilon\in (0,\frac{1}{t_0}c_+^*)$ and $\varphi\in C_+\setminus\{0\}$, there holds $\lim\limits_{t\rightarrow \infty}
\Big[\sup\{||Q_t[\varphi](\cdot, x)-W(\cdot,x)||: x \leq t(\frac{c_+^*}{t_0}-\varepsilon)\}\Big]= 0$.
\end{cor}

\noindent
\textbf {Proof.}  It suffices to prove the necessity since the sufficiency is obvious.  Let  $W=\lim\limits_{t\to \infty}Q_t[r^*]$. According to Theorem~\ref{thm4.1}-(v) and Theorem~\ref{thm4.2}, we have $$\lim\limits_{\alpha\rightarrow \infty}
\Big[\sup\{||Q_t[\varphi](\cdot, x)-W(\cdot,x)||: t\geq t_0\alpha \mbox{ and }\alpha\leq x \leq t(\frac{c_+^*}{t_0}-\varepsilon)\}\Big]= 0.$$ Thus, we only need to prove $\omega(\varphi)=\{W\}$. Otherwise,  there exist $\xi,\eta\in C_{r^*}\setminus \{0\}$ such that $0<\eta\leq\omega(\xi) \leq W$
and $\omega(\xi) \neq W$. Then $\omega(\eta)\leq \omega(\xi)\leq W$,  $\mathcal{K}:=\{\phi\in C_{r^*}\setminus \{0\}:\omega(\eta)\leq \phi \leq \omega(\xi)\}<W$, and $\mathcal{K}$ is a closed, convex, and positively invariant set of $Q_t$. Thus, $\{Q_t\}_{t\in\mathbb{R}_+}$ has another equilibrium point $W_-$ in $\mathcal{K}$. By $W_-<W$, we have $W_-(\cdot,-\infty)=0$. Then Theorem~\ref{thm4.1}-(i) gives  rise to $W_-(\cdot,\infty)=r^*$, a contradiction. \qed

\begin{rem} \label{rem4.1}
Let $t_0>0$ and $\mathcal{X}_+=\Big\{\varphi\in L^\infty(M\times \mathbb{R},\mathbb{R}^N):\varphi(\theta,\cdot)\in X_+,  \mbox{ for all } \theta \in M \mbox{ and } \sup\limits_{(\theta,x)\in M\times \mathbb{R}}||\varphi(\theta,x)||<\infty\Big\}$ with the norm $||\varphi||_{\mathcal{X}}\triangleq
	\sup\limits_{\theta\in M}||\varphi(\theta,\cdot)||_{X}$. For a given map $Q: \mathbb{R}_+\times\mathcal{X}_+\to \mathcal{X}_+$, let
$Q_t:=Q[t,\cdot]$ and define  $Q_\pm$  as in \eqref{Qminus} and \eqref{Qplus}  with $Q$ replaced by $Q_{t_0}|_{C_+}$.
Instead of  (A1), (A2) and (SC), we assume that
\begin{itemize}
\item [{\rm (H1)}]
 $Q_s\circ Q_t=Q_{s+t}$, $Q_t[\xi]\geq Q_t[\eta]$, and $T_{-y}\circ Q_t [\eta]\geq Q\circ T_{-y}[\eta]$
for all $(t,s,y,\eta)\in \mathbb{R}_+^3\times \mathcal{X}_+$ and $\xi\in \eta+\mathcal{X}_+$.

\item [{\rm (H2)}]  $Q[[t_0,\infty)\times \mathcal{X}_+]\subseteq C_+$.

\item [{\rm (H3)}]  For any $r\in Int(\mathbb{R}_+^N)$, $Q[t,\cdot]$ is continuous in $C_r$ uniformly for  $t$ in any compact interval on $\mathbb{R}_+$.

\item [{\rm (H4)}]  There is a continuous-time semiflow $\{Q_t^+\}_{t\in \mathbb{R}_+}$ on $C_+$ such that
for any $(t,\phi)\in \mathbb{R}_+\times C_+$, $\lim\limits_{y\to \infty}T_{-y}\circ Q_t \circ T_{y}[\phi]= Q_t^+[\phi]\in C_+$ with respect to the compact open topology.
 \end{itemize}
It then follows that under assumptions (H1)--(H4) and the other conditions in
Theorem~\ref{thm4.1}, Theorem~\ref{thm4.2}, and Corollary~\ref{cor4.1}, respectively,  all the results in  this section  still hold true.
\end{rem}

\section{Nonautonomous systems}

In this section, we extend our results on spreading speeds and asymptotic
behavior to  a class of nonautonomous
evolution systems.  Assume that
 $P:\mathbb{R}_{+}\times  {C}_+\to{C}_+$ is a map such that for any 
 vector $r\in Int(\mathbb{R}_+^N)$,
$P|_{\mathbb{R}_{+}\times C_r}: \mathbb{R}_{+}\times {C}_r\to{C}_+$
is continuous.  For any given $c\in \mathbb{R}$, we define a family of
mappings $Q_t:=T_{-ct}\circ P[t,\cdot]$  with parameter
$t\in \mathbb{R}_+$.

By Theorem~\ref{thm4.1} and the definition of $Q_t$, we have the following result.

\begin{thm}
\label{thm5.1}
 Assume  that there exist $t_0>0$ and  $c\in \mathbb{R}$ such that $Q_t:=T_{-ct}\circ P[t,\cdot]$ is a continuous-time semiflow on $C_+$,
 and  $Q_t$  satisfies all the conditions in Theorem~\ref{thm4.1}.
 Then the following statements are valid:
\begin{itemize}
\item [{\rm (i)}] If $c_+^*>0$, then for any $\varepsilon\in (0,\frac{1}{t_0}\min\{c_+^*,\frac{c_+^*+c_-^*}{2}\})$ and $\varphi\in C_+\setminus\{0\}$, we have $$\lim\limits_{t\rightarrow \infty}
\max\{|P[t,\varphi](\theta, x)-r^*|:\theta\in M\mbox{ and }t\max\{c+\varepsilon,c-\frac{c_-^*}{t_0}+\varepsilon\}\leq x \leq t(c+\frac{c_+^*}{t_0}-\varepsilon)\}= 0.$$
\item [{\rm (ii)}] If $\varphi
$ has a compact support  and (SC) holds, then $\lim\limits_{t\rightarrow \infty}
\Big[\sup\{
P[t,\varphi](\theta,x):(\theta,x)\in M\times t(\mathbb{R}\setminus [c-\frac{c_-^*}{t_0}-\varepsilon,c+\frac{c_+^*}{t_0}+\varepsilon])\}\Big]=0$ for all $\varepsilon>0$.

\item [{\rm (iii)}]
$\lim\limits_{t\rightarrow \infty}
\Big[\sup\{
||P[t,\varphi](\cdot,x)||:x\in  (-\infty,t(c-\varepsilon)]\}\Big]=0$ for all $\varepsilon>0$ and $\varphi\in C_{r^*}$.

\item [{\rm (iv)}]  If there exists a sequence of points  $\{\phi_k\}_{k \in \mathbb{N}}$ in $Int(Y_+)$  such that  $C_+\subseteq \bigcup\limits_{k\in\mathbb{N}}(\phi_k-C_+)$ and $Q_{t_0}[\phi_k]\leq \phi_k$ for
 all $k\in\mathbb{N}$, then
$\lim\limits_{t\rightarrow \infty}
\Big[\sup\{
||P[t,\varphi](\cdot,x)||:x\in  (-\infty,t(c-\varepsilon)]\}\Big]=0$ for all $\varepsilon>0$ and $\varphi\in C_+$.
\end{itemize}
\end{thm}

\noindent \textbf {Proof.} We only prove (i), since the other cases can be dealt with in a  similar way.
According to Theorem \ref{thm4.1}-(i), we have
$$\lim_{t\rightarrow\infty}\max\big\{ \| Q_{t}[\varphi](\cdot,x)-r^{*} \|:~t\max\{ \varepsilon,-\frac{c_{+}^{*}}{t_0}+\varepsilon\}\leq x\leq t(\frac{c_{+}^{*}}{t_0}-\varepsilon)  \big\}=0.$$
This, together with the definition of $Q_{t}$, implies that
$$\lim_{t\rightarrow\infty}\max\big\{ \| T_{-ct}[P(t,\varphi)](\cdot,x)-r^{*} \|:~t\max\{ \varepsilon,-\frac{c_{+}^{*}}{t_0}+\varepsilon\}\leq x\leq t(\frac{c_{+}^{*}}{t_0}-\varepsilon)  \big\}=0.$$
In other words,
$$\lim_{t\rightarrow\infty}\max\big\{ \| P(t,\varphi)(\cdot,x)-r^{*} \|:~t\max\{ c+\varepsilon,c-\frac{c_{+}^{*}}{t_0}+\varepsilon\}\leq x\leq t(c+\frac{c_{+}^{*}}{t_0}-\varepsilon)  \big\}=0.$$
This completes the proof. \qed

\

We say that $W$ is a travelling wave of $P$ if
$W:M\times\mathbb{R}\to \mathbb{R}_+$ is a bounded and
nonconstant continuous function such that
$P[t,W](\theta,x)=W(\theta,x-tc)$ for all $(\theta,x)\in
M\times \mathbb{R}$ and $t\in \mathbb{R}_+$, and that
$W$ connects $0$ to $r^*$ if $W(\cdot,-\infty)=0$ and
$W(\cdot,\infty)=r^*$. We should point out that our method  for the existence
of travelling waves is quite different from those  in~\cite{bf2018,flw2016,hz2017,wz2018}.

As a consequence of  Theorem~\ref{thm4.2} and the definition of $Q_t$, we have the following result.

\begin{thm} \label{thm5.2} Assume that
 $Q_t:=T_{-ct}\circ P[t,\cdot]$ is a continuous-time semiflow on $C_+$,
 and $Q_t$  satisfies all the conditions in Theorem~\ref{thm4.2}. If $c_-^*>0$, then
$\{P[t,\cdot]\}_{t\in\mathbb{R}_+}$ has  a  travelling wave $W(x-ct,\cdot)$ connecting
$0$ to $r^*$.

\end{thm}

\noindent \textbf {Proof.}  By Theorem \ref{thm4.2}, it follows that $Q_{t}[W]=W$ for some $W\in C_{+}^{\circ}$ with $W$ being nondecreasing,
$W(\cdot,-\infty)=0$ and $W(\cdot,+\infty)=r^{*}$. Thus, $T_{-ct}[P(t,W)]=W$,
that is, $P(t,W)(\theta,x)=W(\theta,x-ct)$ for all $(\theta,x,t)\in M\times\mathbb{R}\times\mathbb{R}_{+}$. \qed

\

According to \cite[Section 3.1]{Zhaobook}, we say a map $Q:\mathbb{R}_{+}\times{C}_+\to{C}_+$ is a
 continuous-time  $\omega$-periodic semiflow on $C_+$ if  for any vector $r\in Int(\mathbb{R}_+^N)$,
 $Q|_{\mathbb{R}_{+}\times C_r}: \mathbb{R}_{+}\times{C}_r\to{C}_+$
 is continuous,  $Q_0=Id|_{{C}_+}$,  and $Q_{t}\circ
 Q_{\omega}=Q_{t+\omega}$ for  some number $\omega>0$ and all $t\in\mathbb{R}_+$,  where $Q_t\triangleq
 Q(t,\cdot)$ for all $t\in\mathbb{R}_+$.

\begin{rem}
In the case where $Q_t:=T_{-ct}\circ P[t,\cdot]$ is a continuous-time
$\omega$-periodic semiflow on $C_+$, we can apply Theorems \ref{thm3.1},
\ref{thm3.2} and \ref{thm3.3} to the Poincar\'e map $Q_{\omega}$  to
establish the spreading properties and the forced time-periodic traveling waves with speed $c$ for the nonautonomous evolution system
$P[t,\cdot]$. We refer to\cite{LYZ} for the Poincar\'e map approach to
monotone periodic semiflows. As an application, one can use the obtained
abstract results to investigate the propagation dynamics of nonlinear
evolution equations in a time-periodic shifting habitat.
\end{rem}

\section{Applications }

In this section, we apply the results obtained in Sections 4 and 5 to
four classes of  monotone evolution equations.
We start with  the definition of KPP property (see, e.g.,~\cite{bhn2005,bhn2010}).
\begin{defn}\label{defn2.1}
Let  $u^{*}\in(0,\infty)$ and $F: \mathbb{R}_{+}\rightarrow\mathbb{R}_{+}$
be  a continuously differentiable function.
We say that $F$ satisfies the KPP property
with respect to $u^{*}$,  or $(F,u^{*})$ has the KPP property
if
\begin{itemize}
\item[{\rm (i)}] $F(0)=0$, $F(u^{*})=0$, and $F'(0)>0$,

\item[{\rm (ii)}] $F(u)(u-u^{*})<0$ for all $u\in(0,\infty)\backslash\{u^{*}\}$,

\item[{\rm (iii)}] $F(u)<F'(0)u$ for all $u\in (0,\infty)$.

\end{itemize}

\end{defn}

For simplicity, we always assume in this section  that $f$ satisfies the following conditions:

\begin{enumerate}
\item[(B1)] $f\in C(\mathbb{R}\times\mathbb{R}_{+},\mathbb{R}) \mbox{ and  } f(s,\cdot)\in C^{1}(\mathbb{R}_{+},\mathbb{R})  \mbox{  for all }
s\in \mathbb{R}$;
\item[(B2)] $f(s,\cdot)\to f_{\pm}^{\infty}(\cdot) \mbox{ in  } C_{loc}^{1}(\mathbb{R}_{+},\mathbb{R})  \mbox{  as } s\to \pm\infty,  \mbox{  where } f_{\pm}^{\infty}(u)\triangleq\lim\limits_{s\to\pm\infty}f(s, u)$;
\item[(B3)] $f_{-}^{\infty}(u_1)\leq f(s_1,u_1)\leq f(s_2,u_2)\leq f_{+}^{\infty}(u_2) \mbox{ for  all } (s_1,u_1),(s_2,u_2)\in \mathbb{R}\times\mathbb{R}_{+}
 \mbox{  with } s_1\leq s_2  \mbox{   and  } u_1\leq u_2$;
\item[(B4)] $\mbox{There exists  }u^{*} >0 \mbox{  such that } f_{+}^{\infty}(u)-u \mbox{  satisfies  the KPP property with respect to } u^{*}$;
\item[(B5)]  $ f_{-}^{\infty}(0)=0 \mbox{   and  }  f_{-}^{\infty}(u)< u  \mbox{  for all  } u\in (0,\infty)$.
\end{enumerate}

By elementary analysis, we have  the following result.
\begin{lemma}\label{lem6.1}
For any $u^{**}\geq u^{*}$ and $\gamma\in(0,\frac{{\rm d}f_+^{\infty}(0)}{{\rm d}u}-1)$, there exist $r=r_{\gamma,u^{**}}(\cdot)\in C(\mathbb{R},\mathbb{R})$ and $K=K_{\gamma,u^{**}}>0$ such that
\begin{itemize}
\item[{\rm (i)}]$r$ is a nondecreasing function with $-\frac{1}{K}<r(-\infty)\leq 0<r(\infty)=\frac{\frac{{\rm d}f_+^{\infty}(0)}{{\rm d}u}-1-\gamma}{K}$;

\item[{\rm (ii)}] $f(s,u)\geq f_{\gamma,u^{**}}(s,u)$ for all $(s,u)\in \mathbb{R}\times[0,u^{**}]$ and $f_{\gamma,u^{**}}$ satisfies (B1)--(B5), where
\[
f_{\gamma,u^{**}}(s,u)=\left\{
\begin{array}{ll}
\frac{(1+Kr(s))^2}{4K}, & (s,u)\in \mathbb{R}\times [\frac{1+Kr(s)}{2K},\infty),
 \\
u+K u(r( s)-u), \qquad & (s,u)\in \mathbb{R}\times [0,\frac{1+Kr(s)}{2K}).
\end {array}
\right.
\]
\end{itemize}
\end{lemma}

\subsection {A  time-delayed nonlocal equation with a shifting habitat \label {Sec6.1}}

Consider the following reaction-diffusion equation with time delay:
\begin {equation}
\left\{
\begin{array}{ll}
\frac{\partial u}{\partial t}(t,x)  =  d u_{xx}(t,x)-\mu u(t,x) +
\mu \int_{\mathbb{R}}f(y-ct,u(t-\tau,y))k(x-y)\mathrm{d}y, \, (t,x)\in (0,\infty)\times \mathbb{R}, \\
u(\theta,x) = \varphi(\theta,x),  \quad  (\theta,x) \in [-\tau,0]
\times \mathbb{R},
\end{array} \right.
\label{6.1}
\end {equation}
where  $c\in \mathbb{R}$,
$\mu>0$, $\tau\geq 0$, $f:\mathbb{R}\times \mathbb{R}_+\to \mathbb{R}_+$
satisfies (B1)--(B5), and the
initial data $\varphi$ belongs to $C([-\tau,0]\times \mathbb{R},\mathbb{R}_+)\cap L^\infty([-\tau,0]\times \mathbb{R},\mathbb{R})$.  Regarding  the kernel function, we always assume that
either $k(x)=\delta(x)$, or $k:\mathbb{R}\to [0,\infty)$ is
continuous with $\int_{\mathbb{R}}k(y)\mathrm{d}y=1$ and
$k(x)=k(-x)$ for all $x\in \mathbb{R}$.

A prototypical  kernel function is  $k(x)=\frac{1}{\sqrt{4\pi\alpha}}e^{-\frac{x^2}{4\alpha}}$, which was used in
\cite{swz2001} to describe the growth of the
matured population of a single species.

Let $M=[-\tau,0]$, $C=BC([-\tau,0]\times \mathbb{R},\mathbb{R})$, and  $C_+=BC([-\tau,0]\times \mathbb{R},\mathbb{R}_+)$.
It is well-known that for any given  $\phi\in  C_+$, equation~\eqref{6.1}  has a
unique solution on a maximal interval $[0, \eta_{\phi;f})$,
 denoted by $u^{\phi}(t,x;f)$ or $(u^{\phi;f})_t$, which is
also the classical solution of~\eqref{6.1} on  $(0,\eta_{\phi;f})$
with
$
[0,\eta_{\phi;f})\ni t\mapsto (u^{\phi;f})_t\in C_+$ being
continuous and $\limsup\limits_{t\to \eta_{\phi;f}^-
}||u^{\phi}(t,\cdot;f)||=\infty$ whenever $\eta_{\phi;f}<\infty$.

By the Phragm\'en-Lindel\"of type maximum principle~\cite{pw1967} and the standard comparison arguments, one can easily get the following result on the global existence, monotonicity, and boundedness of solutions to~\eqref{6.1}.

\begin {prop} \label{prop6.1} Assume that $u_+^*>0$ and $f_1,f_2\in C(\mathbb{R}\times  \mathbb{R}_+,\mathbb{R})$ satisfy   $f_2(\cdot,\cdot)\geq f_1(\cdot,\cdot)$  and  $ f_2(s,u)\leq u$ for all $(s,u)\in \mathbb{R}\times [u_+^*,\infty)$. Let $\psi$, $\phi \in C_+$ with $\phi\leq \psi$. Then  $0\leq  u^\phi(t,x;f_1)\leq u^\psi(t,x;f_2)\leq \max\{||\psi||,u_+^*\}
$
for all $(t,x)\in [0,\min\{\eta_{\phi;f_2},\eta_{\psi;f_1}\}) \times \mathbb{R}$, and hence,  $\eta_{\phi;f_2}=\eta_{\psi;f_1}=\infty$.
\end{prop}

Now we introduce  the following auxiliary equations:

\begin {equation}
\frac {\partial u}{\partial t}= cu_x+d u_{xx}(t,x)-\mu u(t,x) +
\mu \int_{\mathbb{R}}f(y,u(t-\tau,y))k(x+c\tau-y)\mathrm{d}y, \qquad t>0, x\in\mathbb {R},
\label {6.1.1}
\end {equation}

\begin {equation}
\frac {\partial u}{\partial t}= cu_x+d u_{xx}(t,x)-\mu u(t,x) +
\mu \int_{\mathbb{R}}f_\pm^\infty(u(t-\tau,y))k(x+c\tau-y)\mathrm{d}y, \qquad t>0, x\in\mathbb {R},
\label {6.1.2}
\end {equation}
and
\begin {equation}
\frac {\partial u}{\partial t}= d u_{xx}(t,x)-\mu u(t,x) +
\mu \int_{\mathbb{R}}f_\pm^\infty(u(t-\tau,y))k(x-y)\mathrm{d}y, \qquad t>0, x\in\mathbb {R}.
\label {6.1.3}
\end {equation}
Define  $P:\mathbb{R}_+\times C_+\to C_+$ by
$P[f;t,\phi](\theta,x)=u^{\phi}(t+\theta,x;f)$ for all $(t,\theta,\phi)\in \mathbb{R}_+\times  [-\tau,0] \times C_+$.  Let $Q[f;t,\phi]$, $Q_\pm[f_\pm^\infty;t,\phi]$ and $\Phi_\pm[f_\pm^\infty;t,\phi]$ be the mild solutions  of~\eqref{6.1.1}, \eqref{6.1.2}  and \eqref{6.1.3} with the initial value $u_0 = \phi\in C_+$, respectively.   For simplity, we denote
$P[f;t,\phi]$, $Q[f;t,\phi]$, $Q_\pm[f_\pm^\infty;t,\phi]$ and $\Phi_\pm[f_\pm^\infty;t,\phi]$ by $P[t,\phi]$, $Q[t,\phi]$, $Q_\pm[t,\phi]$ and $\Phi_\pm[t,\phi]$, respectively.

\begin{prop} \label{prop6.2} Let $t\in \mathbb{R}$ and $\phi\in C_+$. Then the following statements are valid:
\begin{itemize}
\item [{\rm (i)}]  $Q[t,\phi](\theta,x)=P[t,\phi](\theta,x+ct)$ for all $(\theta,x)\in [-\tau,0]\times \mathbb{R}$.

\item [{\rm (ii)}]  $ Q_\pm[t,\phi]=\lim\limits_{y\to \pm \infty} Q[t,T_{y}[\phi]](\cdot,\cdot+y),$ that is,  $Q[t,T_{y}[\phi]](\cdot,\cdot+y)\to  Q_\pm[t,\phi]$ in $C$ as $y\to \pm \infty$.

\item [{\rm (iii)}]   $\Phi_\pm[t,\phi]=Q_\pm[t,\phi](\cdot,\cdot-ct).$
\end{itemize}

\end{prop}

\noindent
{\bf Proof.}  By straightforward  computations, we can directly verify  (i) and (iii).

(ii)  It is easy to see that for any  $z\in \mathbb{R}$, $u(t+\theta,x):=Q[t,T_{z}[\phi]](\theta,x+z)$ satisfies
\begin {equation*}
\left\{
\begin{array}{ll}
\frac {\partial u}{\partial t}=cu_x+d u_{xx}(t,x)-\mu u(t,x) +
\mu \int_{\mathbb{R}}f(y+z,u(t-\tau,y))k(x+c\tau-y)\mathrm{d}y, \, t>0, x\in\mathbb {R},\\
 u(\theta,x)=\phi(x), \qquad (\theta,x)\in
    [-\tau,0]\times \mathbb{R}.
\end{array} \right.
\label{6.1}
\end {equation*}
This, together with the fact that $f(y+z,u) \to f_{\pm}(u)$ locally uniform for $y$ and $u$  as $z\to \infty$, yields (ii).  \qed

\begin{prop} \label{prop6.3} Let  $c^*>0$ be the spreading speed  of  system \eqref{6.1.3} with $f_+^\infty$ (see \cite{lz2007}).
 For a given number $t_0>\tau$,
let  $c^*_+=\bar{c}_+=t_0(c^*-c)$, $c_-^*=\bar{c}_-=t_0(c+c^*)$,  $r^*=u^*$, $Q_t:=Q[t,\cdot]$, and $Q_\pm:=Q_\pm[t_0,\cdot]$ for all $t\in \mathbb{R}_+$. Then the following statements hold true:
\begin{itemize}
\item [{\rm (i)}]  $\{Q_t\}^{\infty}_{t=0}$ is a continuous-time semiflow on $C_+$  and satisfies (A1-A4) and (SP).

\item [{\rm (ii)}] $\{(T_{-z}\circ Q_{t_0})^n[r^*]:n\in \mathbb{N}\}$ is precompact in $C$  for all $z\in \mathbb{R}$.

\item [{\rm (iii)}]    $Q_+$ satisfies  (UC) and (AA).

\item [{\rm (iv)}]    $Q_-$ satisfies (UAA).

\item [{\rm (v)}]   $Q_t$ satisfies (ASH-UC-SP) for all $t>0$.
\end{itemize}

\end{prop}

\noindent
{\bf Proof.}  Since (i), (ii), and (iv) are obvious, we only verify (iii) and (v).

(iii) Applying  \cite[Theorem 5.1]{lz2007} to system \eqref{6.1.3},   we know that
any $\varepsilon>0$ and $\phi\in C_+$,
\[
\lim\limits_{t\to \infty} \max\{\Phi[t,\phi](\theta, x): \theta\in [-\tau,0] \mbox{ and } |x|\geq t(c^*+\varepsilon)\}=0 \mbox{ if } \phi \mbox{ has a compact support, }
\]
\mbox{ and }
\[
\lim\limits_{t\to \infty} \max\{|\Phi[t,\phi](\theta, x)-r^*|:  \theta\in [-\tau,0] \mbox{ and } |x|\leq
t (c^*-\varepsilon)\}=0 \mbox{ if } \phi\neq 0.
\]
Thus,  Proposition~\ref{prop6.2}-(iii) implies (iii).

(v)  Take  $\gamma_l\in (0, \frac{{\rm d} f_+^\infty(0)}{{\rm d} u}-1)$ with $\lim\limits_{l\to \infty} \gamma_l=0$ and $\gamma_l>\gamma_k$ for all  positive integers $k>l$.
Let $r_l^*=\frac{\frac{{\rm d}f_+^{\infty}(0)}{{\rm d}u}-1-\gamma_l}{K_{\gamma_l,u^{*}}}$, $Q_l[t,\phi]=Q[f_{\gamma_l,u^*};t,\phi]$, and let $c_{\pm l}^*$ be the spreading speed  of  \eqref{6.1.3} with $f_{\gamma_l,u^{*}}$, where $K_{\gamma_l,u^{*}}$ and $f_{\gamma_l,u^{*}}$  are defined as in Lemma~\ref{lem6.1}.  Then (v) follows from  Lemma~\ref{lem6.1} and Propositions~\ref{prop6.1}, ~\ref{prop6.2}-(ii), and  ~\ref{prop6.3}-(i,iii). \qed

\

As a straightforward consequence of  Proposition~\ref{prop6.3} and Theorems~\ref{thm5.1} and~\ref{thm5.2}, we have the following result
for system  ~\eqref{6.1}.

\begin{thm}\label{thm6.1}
Assume that  $f$ satisfies  (B1)--(B5).   Let  $c^*>0$ be the spreading speed  of  system \eqref{6.1.3} with $f_+^\infty$. Then the following statements are valid:
\begin{itemize}
\item [{\rm (i)}] If $c<c^*$, then for any $\varepsilon\in (0,
\min\{c^*-c,c^*\})$ and $\varphi>0$, we have $$\lim\limits_{n\rightarrow \infty}
\max\{|P[t,\varphi](\theta, x)-r^*|:\theta\in  [-\tau,0]\mbox{ and }t\max\{c+\varepsilon,-c^*+\varepsilon\}\leq x \leq t(c^*-\varepsilon)\}= 0.$$

\item [{\rm (ii)}] If $\varphi\in C_+$ has a compact support, then $\lim\limits_{t\rightarrow \infty}
\Big[\sup\{
P[t,\varphi](\theta,x):(\theta,x)\in  [-\tau,0]\times t(\mathbb{R}\setminus [-c^*-\varepsilon,c^*+\varepsilon])\}\Big]=0$ for all $\varepsilon>0$.

\item [{\rm (iii)}]  For any $\varepsilon>0$ and $\varphi\in C_+$, we have
$\lim\limits_{t\rightarrow \infty}
\Big[\sup\{
||P[t,\varphi](\cdot,x)||:x\in  (-\infty,t(c-\varepsilon)]\}\Big]=0$.

 \item [{\rm (iv)}]  If $c>-c^{*} $, then $\{P[t,\cdot]\}_{t\in\mathbb{R}_+}$ has  a  travelling wave $W(x-ct)$ connecting
$0$ to $r^*$.

\end{itemize}

\end{thm}

In the case where  $\tau=0$ and $k(\cdot)=\delta$, we have the following result.
\begin{cor}\label{cor6.3}
Assume that  $f$ satisfies  (B1)--(B5). Then $c^*:=2\sqrt{\mu d(\lim\limits_{s\to \infty}\frac{{\rm d } f(s,u)}{{\rm d} u}|_{u=0}-1)}$ and the following statements are valid:
\begin{itemize}
\item [{\rm (i)}] If $c<c^*$, then for any $\varepsilon\in (0,
	\min\{c^*-c,c^*\})$ and $\varphi> 0$, we have $$\lim\limits_{n\rightarrow \infty}
\max\{|P[t,\varphi](0, x)-r^*|: t\max\{c+\varepsilon,-c^*+\varepsilon\}\leq x \leq t(c^*-\varepsilon)\}= 0.$$

\item [{\rm (ii)}] If $\varphi\in C_+$ has a compact support, then $\lim\limits_{t\rightarrow \infty}
\Big[\sup\{
P[t,\varphi](\theta,x): x\in t(\mathbb{R}\setminus [-c^*-\varepsilon,c^*+\varepsilon])\}\Big]=0$ for all $\varepsilon>0$.

\item [{\rm (iii)}]  For any $\varepsilon>0$ and $\varphi\in C_+$, we have
$\lim\limits_{t\rightarrow \infty}
\Big[\sup\{
|P[t,\varphi](0,x)|:x\in  (-\infty,t(c-\varepsilon)]\}\Big]=0$.

 \item [{\rm (iv)}]  If $c>-c^{*} $, then $\{P[t,\cdot]\}_{t\in\mathbb{R}_+}$ has  a  travelling wave $W(x-ct)$ connecting
$0$ to $r^*$.

\end{itemize}

\end{cor}

It is worthy  pointing  out that if $\mu f(s,u)-\mu u=u[{r}(s)-u]$ and $r\in C(\mathbb{R},\mathbb{R})$ is a nondecreasing function   with
$\lim\limits_{s\rightarrow\infty}{r}(s)>0$ and $\lim\limits_{s\rightarrow\infty}{r}(s)\leq 0$, then we can apply  Corollary~\ref{cor6.3}  with $f=f_\mu$ for all large $\mu$ to show that  all results of Corollary~\ref{cor6.3} still hold true for the non-monotone $f(s,\cdot)$, where \[
f_\mu(u)=\left\{
\begin{array}{ll}
u+\frac{u[r(s)-u]}{\mu}, & (s,u)\in \mathbb{R}\times [0,\frac{\mu+r(s)}{2}],
 \\
\frac{[\mu+r(s)]^2}{4\mu}, \qquad & (s,u)\in \mathbb{R}\times (\frac{\mu+r(s)}{2},\infty).
\end {array}
\right.
\]

We also remark that Corollary~\ref{cor6.3}-(i-iii) was  obtained in \cite{lbsf2014} in the case where $\lim\limits_{s\rightarrow-\infty}{r}(s)< 0$ and $c>0$;  while Corollary~\ref{cor6.3}-(i-iii) was established in \cite{hyz2019} in the case where $\lim\limits_{s\rightarrow-\infty}{r}(s)=0$ and $c>0$.  Further,  Corollary~\ref{cor6.3}-(iv) was proved in \cite{bf2018,flw2016,hz2017}  via the method of sup- and subsolutions.

In the rest of this subsection,  we consider the following time-delayed
nonlocal dispersal equation:
\begin {equation}
\left\{
\begin{array}{rcl}
\frac{\partial u}{\partial t}(t,x) & = &
d[\int_{\mathbb{R}}u(t,y)k(x-y)\mathrm{d}y-u(t,x)]-\mu u(t,x) + \mu
f(x-ct,u(t-\tau,x))
\\
& & \hspace*{2in} (t,x)\in (0,\infty)\times \mathbb{R},
\\
u(\theta,x)& =& \varphi(\theta,x),  \quad  (\theta,x) \in [-\tau,0]
\times \mathbb{R},
\end{array} \right.
\label{6.3}
\end {equation}
where $d,\mu>0$, $\tau\geq 0$, the initial data $\varphi$ belongs
to $C_+:=BC([-\tau,0]\times \mathbb{R},\mathbb{R}_+)$, $f:\mathbb{R}_+\to
\mathbb{R}_+$ satisfies (B1)--(B5), and the kernel  $k:\mathbb{R}\to (0,\infty)$ is continuous and symmetric with
 $\int_{\mathbb{R}}k(y)\mathrm{d}y=1$ and
$\int_{\mathbb{R}} e^{\rho y} k(y)\mathrm{d}y <\infty$ for $\rho \in
\mathbb{R}$.

It is well-known that for any given  $\phi\in  C_+$, equation~\eqref{6.3}  has a
unique solution on a maximal interval $[0, \eta_{\phi;f})$, denoted by $u^{\phi}(t,x;f)$ or $(u^{\phi;f})_t$, which is
also the classical solution of~\eqref{6.3} on  $(0,\eta_{\phi;f})$
with
$
[0,\eta_{\phi;f})\ni t\mapsto (u^{\phi;f})_t\in C_+$ being
continuous and $\limsup\limits_{t\to \eta_{\phi;f}^-
}||u^{\phi}(t,\cdot;f)||=\infty$ whenever $\eta_{\phi;f}<\infty$.

We introduce  the following auxiliary equations:
\begin {equation}
\frac {\partial u}{\partial t}= cu_x(t,x)+d[\int_{\mathbb{R}}u(t,y)k(x-y)\mathrm{d}y-u(t,x)]-\mu u(t,x) + \mu
f(x,u(t-\tau,x+c\tau)), \qquad t>0, x\in\mathbb {R},
\label {6.2.1}
\end {equation}

\begin {equation}
\frac {\partial u}{\partial t}= cu_x(t,x)+d[\int_{\mathbb{R}}u(t,y)k(x-y)\mathrm{d}y-u(t,x)]-\mu u(t,x) + \mu
f_\pm^\infty (u(t-\tau,x+c\tau)), \qquad t>0, x\in\mathbb {R},
\label {6.2.2}
\end {equation}
and
\begin {equation}
\frac {\partial u}{\partial t}= d[\int_{\mathbb{R}}u(t,y)k(x-y)\mathrm{d}y-u(t,x)]-\mu u(t,x) + \mu
f_\pm^\infty(u(t-\tau,x)), \qquad t>0, x\in\mathbb {R}.
\label {6.2.3}
\end {equation}
Define  $P:\mathbb{R}_+\times C_+\to C_+$ by
$P[t,\phi](\theta,x)=u^{\phi}(t+\theta,x;f)$ for all $(t,\theta,\phi)\in \mathbb{R}_+\times  [-\tau,0] \times C_+$.  Let  $Q[t,\phi]$, $Q_\pm[t,\phi]$, and $\Phi[t,\phi]$ be the mild solutions of~\eqref{6.2.1}, \eqref{6.2.2}  and \eqref{6.2.3} with the initial value $u_0 = \phi\in C_+$, respectively.

Since the map $P[t,\cdot]: C_+\to C_+$ is not compact  for  any $t>0$, we first reduce the 
existence of a traveling wave with speed $c$ of system \eqref{6.3} to that of  a  fixed point of an appropriate map with parameter $c$  on $\mathcal{C}_+:=BC(\mathbb{R},\mathbb{R}_+)$ equipped with the compact open topology.
Let  $c^*>0$ be the spreading speed  of system \eqref{6.2.3} with $f_+^\infty$
(see \cite{lz2007}). Then $c^*$ is also the minimum wave speed  of
monotone traveling waves for  system \eqref{6.2.3} with $f_+^\infty$
(see \cite{fz2014}).

Define
$$
\hat{k} (\rho)=\int_{\mathbb{R}}e^{\rho y}k(y)\, d y,
\quad l (c,\rho)=\frac{1}{c\rho +d+\mu}[d \hat{k}
(\rho)+\mu \frac{{\rm d } f_+^\infty}{{\rm d }u}(0)e^{-\rho c \tau}],\,\,
\forall  \rho,\,c  \in \mathbb{R},
$$
\[
l^\pm(c,\rho)=\left\{
\begin{array}{ll}
l(c,\pm\rho), & (c,\rho)\in J_\pm\triangleq \{(a,b)\in \mathbb{R}\times (0,\infty):
d+\mu\pm ab>0\},
 \\
\infty, \qquad & (c,\rho)\notin J_\pm.
\end {array}
\right.
\]
and 
$$c_\pm^*(c)=\inf\limits_{\rho>0}\frac{1}{\rho}\log l^\pm(c,\rho),
\,\, \forall  c\in \mathbb{R}. 
$$
In view of  \cite[Lemma 4.9]{yz2015}, we then have the following observation.

\begin{prop}\label{lemm6.8000}
 The following statements are valid:
\begin{itemize}
\item [{\rm (i)}] $c^*=\inf\{c\in \mathbb{R}:c_+^*(c)\leq 0\}=
-\sup\{c\in \mathbb{R}:c^*_-(c)\leq 0\}$.
\item [{\rm (ii)}] $c_+^*(c)+c_-^*(c)>0$ for all
$c\in \mathbb{R}$.
\item [{\rm (iii)}] $c_-^*(c)>0$ for all
$c>-c^*$,  and  $c_+^*(c)>0$ for all
$c<c^*$.
\end{itemize}
\end{prop}

For  any given $c\in \mathbb{R}$,
we define $K[\cdot;c], K_\pm[\cdot;c], L[\cdot;c],Q[\cdot;c],Q_\pm[\cdot;c]: \mathcal{C}_+\to
\mathcal{C}_+$ by
\begin {eqnarray*}
&& K[\phi;c](x)=\frac{1}{d+\mu}
\left[d\int_{\mathbb{R}}\phi(y)k(x-y)\, dy+\mu
f(x,\phi(x+c\tau))\right], \\
&& K_\pm[\phi;c](x)=\frac{1}{d+\mu}
\left[d\int_{\mathbb{R}}\phi(y)k(x-y)\, dy+\mu
f_\pm^\infty(\phi(x+c\tau))\right], \\
&& L[\phi;c](x)=\left\{
\begin{array}{ll}
\frac{d+\mu}{c}\int_x^{\infty} e^{\frac{d+\mu}{c}(x-y)}\phi(y)\, dy,\qquad & c>0,\\
\phi(x), & c=0,\\
-\frac{d+\mu}{c}\int^x_{-\infty} e^{\frac{d+\mu}{c}(x-y)}\phi(y)\, dy, & c<0,
\end{array}
\right.\\
&& Q[\phi;c](x)=L[K[\phi]](x),
\\
&& Q_\pm[\phi;c](x)=L[K_\pm[\phi]](x),\end{eqnarray*}
for all  $\phi\in \mathcal{C}_+$ and $x\in \mathbb{R}$.
By  using \cite[Lemma 4.7]{yz2015}, Lemma~\ref{lem6.1}, and  Proposition~\ref{lemm6.8000}, we can verify the following
properties for the  maps $Q$ and $Q_\pm$ defined above.

\begin{prop} \label{prop6.3000000000} Assume that  $f$ satisfies (B1)--(B5).   Let $c\in \mathbb{R}$, $c^*_+=c^*_+(c)$, $c_-^*=c^*_-(c)$,
and $r^*=u^*$. Then the following statements hold true:
\begin{itemize}
\item [{\rm (i)}]  $Q[\cdot;c]$ and $Q_\pm[\cdot;c]$ are  continuous maps on $\mathcal{C}_+$  and satisfy (A1-A4) and (SP).

\item [{\rm (ii)}]    $Q_+[\cdot;c]$ satisfies  (UC) and (AA).

\item [{\rm (iii)}]    $Q_-[\cdot;c]$ satisfies (UAA).

\item [{\rm (iv)}]   $Q[\cdot;c]$ satisfies (ASH-UC-SP).

\item [{\rm (v)}]   $Q[\cdot;c]$ is a compact map on $\mathcal{C}_+$ for any
$c\neq 0$.

\item [{\rm (vi)}] $\{Q^n[r^*;0]:n\in \mathbb{N}\}$ is precompact in $\mathcal{C}_+$
provided $\frac{{\rm{d}}f(x,u)}{{\rm{d}} u}<\frac{d+\mu}{\mu}, \,
\, \forall (x,u)\in \mathbb{R}\times (0,r^*).$

\item [{\rm (vii)}]  If $Q[\cdot;c]$ has a fixed point $W$ in
	$\mathcal{C}_+$, then  $W(x-ct)$ is a  travelling wave of 
	$\{P[t,\cdot]\}_{t\in\mathbb{R}_+}$. 
\end{itemize}
\end{prop}

Now we are ready to present  the result on the propagation dynamics of  system \eqref{6.3}.

\begin{thm}\label{thm6.2}
Assume that  $f$ satisfies (B1)--(B5).    Then the following statements are valid:
\begin{itemize}
\item [{\rm (i)}] If $c<c^*$, then for any
$\varepsilon\in (0,\min\{c^*-c,c^*\})$ and $\varphi> 0$, we have $$\lim\limits_{n\rightarrow \infty}
\max\{|P[t,\varphi](\theta, x)-r^*|:\theta\in [-\tau,0]\mbox{ and }t\max\{c+\varepsilon,-c^*+\varepsilon\}\leq x \leq t(c^*-\varepsilon)\}= 0.$$
\item [{\rm (ii)}] If $\varphi \in C_+
$ has a compact support, then $\lim\limits_{t\rightarrow \infty}
\Big[\sup\{
P[t,\varphi](\theta,x):(\theta,x)\in [-\tau,0]\times t(\mathbb{R}\setminus [-c^*-\varepsilon,c^*+\varepsilon])\}\Big]=0$ for all $\varepsilon>0$.

\item [{\rm (iii)}]  For any $\varepsilon>0$ and $\varphi \in C_+$, we have
$\lim\limits_{t\rightarrow \infty}
\Big[\sup\{
||P[t,\varphi](\cdot,x)||:x\in  (-\infty,t(c-\varepsilon)]\}\Big]=0$.

 \item [{\rm (iv)}]  If $c>-c^{*} $ and $\frac{{\rm{d}}f(x,u)}{{\rm{d}} u}<\frac{d+\mu}{\mu}$ for all $(x,u)\in \mathbb{R}\times (0,r^*)$ whence $c=0$, then $\{P[t,\cdot]\}_{t\in\mathbb{R}_+}$ has  a  travelling wave $W(x-ct)$ connecting
$0$ to $r^*$.
\end{itemize}
\end{thm}
\noindent
{\bf Proof.}
Statements (i), (ii)  and (iii) follow from  the essentially same arguments as those for system \eqref{6.1}, and  (iv) is a consequence
of  Theorem~\ref{thm3.3} and Propositions~\ref{lemm6.8000} and \ref{prop6.3000000000}. \qed

\

In the case where  $\tau=0$, we have the following result.

\begin{cor}\label{cor6.4}
Assume that  $f$ satisfies (B1)--(B5), and  let   $c^*>0$ be the spreading speed  of  system \eqref{6.2.3} with $\tau=0$ and $f_+^\infty$. Then  the following statements are valid:
\begin{itemize}
\item [{\rm (i)}] If $c<c^*$, then for any $\varepsilon\in (0,\min\{c^*-c,c^*\})$ and $\varphi>0$, we have $$\lim\limits_{n\rightarrow \infty}
\max\{|P[t,\varphi](0, x)-r^*|:t\max\{c+\varepsilon,-c^*+\varepsilon\}\leq x \leq t(c^*-\varepsilon)\}= 0.$$
\item [{\rm (ii)}] If $\varphi \in C_+
$ has a compact support, then $\lim\limits_{t\rightarrow \infty}
\Big[\sup\{
P[t,\varphi](0,x): x\in t(\mathbb{R}\setminus [-c^*-\varepsilon,c^*+\varepsilon])\}\Big]=0$ for all $\varepsilon>0$.

\item [{\rm (iii)}]  For any $\varepsilon>0$ and $\varphi \in C_+
$, we have
$\lim\limits_{t\rightarrow \infty}
\Big[\sup\{
||P[t,\varphi](0,x)||:x\in  (-\infty,t(c-\varepsilon)]\}\Big]=0$.

 \item [{\rm (iv)}]  If $c>-c^{*} $, then $\{P[t,\cdot]\}_{t\in\mathbb{R}_+}$ has  a  travelling wave $W(x-ct)$ connecting
$0$ to $r^*$.
\end{itemize}
\end{cor}

We note that in the case where $\mu f(s,u)-\mu u=u[{r}(s)-u]$, $c>0$, and $r\in C(\mathbb{R},\mathbb{R})$ is a nondecreasing function   with
$\lim\limits_{s\rightarrow\infty}{r}(s)>0$ and $\lim\limits_{s\rightarrow-\infty}{r}(s)< 0$,  Corollary~\ref{cor6.4}-(i-iii) was obtained in \cite{lwz2018},  and Corollary~\ref{cor6.4}-(iv) was proved in \cite{wz2018}  by using the method of sup- and subsolutions.

\subsection{A reaction-diffusion equation in a cylinder}
Consider the following reaction-diffusion equation in a cylinder and with a shifted habitat:
\begin {equation}\label{6.13}
\left\{
\begin{array}{ll}
&\frac{\partial u}{\partial t}  =  \frac{\partial^{2}u}{\partial x^{2}} + \Delta_{y}u +ug(x-ct,y,u),~~~~~x\in\mathbb{R},y\in\Omega\subseteq\mathbb{R}^{m},t>0\\
&\frac{\partial u}{\partial \nu} = 0,~~~~~\mbox{on}~\mathbb{R}\times\partial\Omega\times(0,\infty),\\
&u(t_{0},x,y)=\phi(x,y),~~~~~(x,y)\in\mathbb{R}\times\overline{\Omega},
\end{array} \right.
\end {equation}
where $\Omega$ is a bounded domain in $\mathbb{R}^{m}$ with smooth boundary $\partial\Omega$, $\Delta_{y}=\sum\limits_{i=1}^{m}\frac{\partial^{2}}{\partial y_{i}^{2}}$,
and $\nu$ is the outer unit normal vector to $\mathbb{R}\times\partial\Omega$.

Let $\lambda_0$ be the principal eigenvalue of the elliptic eigenvalue problem
\begin {equation*}
\left\{
\begin{array}{ll}
	&\lambda \varphi (y)= \Delta_{y}\varphi (y) +g^{\infty}_+(y,0)\varphi (y),~~~~~y\in\Omega,\\
	&\frac{\partial \varphi }{\partial \nu} = 0,~~~~~\mbox{on}~\partial\Omega.
\end{array} \right.
\end {equation*}
We assume that
\begin{enumerate}
\item[(C1)] $g\in C(\mathbb{R}\times \overline{\Omega}\times\mathbb{R}_{+},\mathbb{R}) \mbox{~and~}
      g(s,\cdot,\cdot)\in C^{1}(\overline{\Omega}\times\mathbb{R}_{+},\mathbb{R}) \mbox{~for~all~}s\in\mathbb{R}$;
 \item[(C2)] $g(s,\cdot,\cdot)\rightarrow g_{\pm}^{\infty}(\cdot,\cdot)$ in $C_{loc}^{1}(\overline{\Omega}\times\mathbb{R}_{+},\mathbb{R})$
 as $s\rightarrow\pm\infty$, where $g_{\pm}^{\infty}(y,u):=\lim\limits_{s\rightarrow\pm\infty}g(s,y,u)$
for all $(y,u)\in\overline{\Omega}\times\mathbb{R}_{+}$;
\item[(C3)]  $g_{-}^{\infty}(y,u)\leq g(s_{1},y,u)\leq g(s_{2},y,u)\leq g_{+}^{\infty}(y,u)\mbox{~for~all~}s_{1},s_{2}\in\mathbb{R} \mbox{~and~}
        (y,u)\in\overline{\Omega}\times\mathbb{R}_{+}
        \mbox{~with~}s_{1}\leq s_{2}$;
 \item[(C4)]  $g_{-}^{\infty}(y,u)\leq 0 \mbox{~and~}g_{-}^{\infty}(y,0)=0\mbox{~for~all~}(y,u)\in\overline{\Omega}\times\mathbb{R}_{+}$;
  \item[(C5)] $\frac{\partial g_{+}^{\infty}}{\partial u}(y,u)<0\mbox{~for~all~}(y,u)\in\overline{\Omega}\times\mathbb{R}_{+},
        \mbox{~and~there~is~}K>0\mbox{~such~that~}g_{+}^{\infty}(y,u)\leq0\mbox{~for~all~}
        (y,u)\in\overline{\Omega}\times[K,\infty)$;
 \item[(C6)]  $\lambda_0>0$.
 \end {enumerate}

Let $M=\overline{\Omega}$, $C=C(\mathbb{R}\times\overline{\Omega},\mathbb{R})\cap L^{\infty}(\mathbb{R}\times\overline{\Omega},\mathbb{R})$, and  $C_+=C(\mathbb{R}\times\overline{\Omega},\mathbb{R}_+)\cap L^{\infty}(\mathbb{R}\times\overline{\Omega},\mathbb{R})$.
It is well-known that for any given  $\phi\in  C_+$, equation~\eqref{6.13}  has a
unique solution on a maximal interval $[0, \eta_{\phi;g})$,
 denoted by $u^{\phi}(t,x,y;g)$ or $(u^{\phi;g})_t$, which is
also the classical solution of~\eqref{6.13} on  $(0,\eta_{\phi;g})$
with $[0,\eta_{\phi;g})\ni t\mapsto (u^{\phi;g})_t\in C_+$ being
continuous and $\limsup\limits_{t\to \eta_{\phi;g}^-}||u^{\phi}(t,\cdot,\cdot;g)||=\infty$ whenever $\eta_{\phi;g}<\infty$.

By the Phragm\'en-Lindel\"of type maximum principle~\cite{pw1967} and the standard comparison arguments, one can easily get the following result on the global existence, monotonicity, and boundedness of solutions to~\eqref{6.13}.

\begin {prop} \label{prop6.10}
Assume that $M^*>0$ and $g_1,g_2\in C(\mathbb{R}\times\overline{\Omega}\times  \mathbb{R}_+,\mathbb{R})$ satisfy   $g_2(\cdot,\cdot,\cdot)\geq g_1(\cdot,\cdot,\cdot)$
and  $ g_2(s,y,u)\leq 0$ for all $(s,y,u)\in \mathbb{R}\times \overline{\Omega} \times [M^*,\infty)$. Let $\psi$, $\phi \in C_+$ with $\phi\leq \psi$.
Then  $0\leq  u^\phi(t,x,y;g_1)\leq u^\psi(t,x,y;g_2)\leq \max\{||\psi||,M^*\}$
for all $(t,x,y)\in [0,\min\{\eta_{\phi;f_2},\eta_{\psi;g_1}\}) \times \mathbb{R}\times \overline{\Omega}$, and hence,  $\eta_{\phi;g_2}=\eta_{\psi;g_1}=\infty$.
\end{prop}

Now we introduce  the following auxiliary equations:
\begin {equation}\label{6.5.1}
\left\{
\begin{array}{ll}
\frac{\partial u}{\partial t}  &=  \frac{\partial^{2}u}{\partial x^{2}} + \Delta_{y}u
+c\frac{\partial u}{\partial x}+ug(x,y,u),~~~~~(x,y,t)\in\mathbb{R}\times\Omega\times(0,\infty),\\
\frac{\partial u}{\partial\overrightarrow{n}} &= 0,~~~~~\mbox{on}~\mathbb{R}\times\partial\Omega\times(0,\infty),
\end{array} \right.
\end {equation}

\begin {equation}\label{6.5.2}
\left\{
\begin{array}{ll}
\frac{\partial u}{\partial t}  &=  \frac{\partial^{2}u}{\partial x^{2}} + \Delta_{y}u
+c\frac{\partial u}{\partial x}+ug_{\pm}^{\infty}(y,u),~~~~~(x,y,t)\in\mathbb{R}\times\Omega\times(0,\infty),\\
\frac{\partial u}{\partial\overrightarrow{n}} &= 0,~~~~~\mbox{on}~\mathbb{R}\times\partial\Omega\times(0,\infty),
\end{array} \right.
\end {equation}
and
\begin {equation}\label{6.5.3}
\left\{
\begin{array}{ll}
\frac{\partial u}{\partial t}  &=  \frac{\partial^{2}u}{\partial x^{2}} + \Delta_{y}u
+ug_{\pm}^{\infty}(y,u),~~~~~(x,y,t)\in\mathbb{R}\times\Omega\times(0,\infty),\\
\frac{\partial u}{\partial \nu} &= 0,~~~~~\mbox{on}~\mathbb{R}\times\partial\Omega\times(0,\infty).
\end{array} \right.
\end {equation}
Define  $P:\mathbb{R}_+\times C_+\to C_+$ by
$P[g;t,\phi](x,y)=u^{\phi}(t,x,y;g)$ for all $(t,x,y,\phi)\in \mathbb{R}_+\times \mathbb{R}\times\overline{\Omega} \times C_+$.
Let $Q[g;t,\phi]$, $Q_\pm[g_\pm^\infty;t,\phi]$ and $\Phi_\pm[g_\pm^\infty;t,\phi]$ be the mild solutions of~\eqref{6.5.1}, \eqref{6.5.2} and \eqref{6.5.3} with the initial value $u(0,\cdot) = \phi\in C_+$,
 respectively. For simplicity, we denote $P[g;t,\phi]$, $Q[g;t,\phi]$, $Q_\pm[g_\pm^\infty;t,\phi]$,
and $\Phi_\pm[g_\pm^\infty;t,\phi]$ by $P[t,\phi]$, $Q[t,\phi]$, $Q_\pm[t,\phi]$, and $\Phi_\pm[t,\phi]$

\begin{prop} \label{prop6.14} Let $t\in \mathbb{R}$ and $\phi\in C_+$. Then the following statements are valid:
\begin{itemize}
\item [{\rm (i)}]  $Q[t,\phi](x,y)=P[t,\phi](x+ct,y)$ for all $(x,y)\in \mathbb{R}\times\overline{\Omega}$.

\item [{\rm (ii)}]  $ Q_\pm[t,\phi]=\lim\limits_{y\to \pm \infty} Q[t,T_{y}[\phi]](\cdot+y,\cdot),$ that is,  $Q[t,T_{y}[\phi]](\cdot+y,\cdot)\to  Q_\pm[t,\phi]$ in $C$ as $y\to \pm \infty$.

\item [{\rm (iii)}]   $\Phi_\pm[t,\phi]=Q_\pm[t,\phi](\cdot-ct,\cdot).$
\end{itemize}
\end{prop}

\noindent
{\bf Proof.}  By straightforward  computations, we can directly verify  (i) and (iii).

(ii)  It is easy to see that for any  $z\in \mathbb{R}$, $u(t,x,y):=Q[t,T_{z}[\phi]](x+z,y)$ satisfies
\begin {equation*}
\begin{array}{rcl}
    &&\frac {\partial u}{\partial t}=\frac {\partial^{2} u}{\partial x^{2}}+\Delta_{y}u+c\frac{\partial u}{\partial x}+ug(x+z,y,u),~~~~~(t,x,y)\in(0,\infty)\times\mathbb {R}\times\overline{\Omega},\\
   && u(x,y) = \phi(x,y), \qquad (x,y)\in    \mathbb{R}\times\overline{\Omega}.
\end{array}
\end {equation*}
 This, together with the fact that $g(x+z,y,u) \to g_{\pm}^{\infty}(y,u)$ locally uniform for $(x,y,u)$  as $z\to \infty$, yields (ii).  \qed

\

 According to \cite{lz2007}, system \eqref{6.5.3} with $g_+^{\infty}(y,u)$ admits a
unique positive $x$-independent steady state $\beta(y)$ and $c^*:=2\sqrt{\lambda_0}$ is the spreading speed
for its solutions with initial data having compact supports.

\begin{prop} \label{prop6.15}
Let $c^*=2\sqrt{\lambda_0}$ and set $c^*_+=\bar{c}_+=c^*-c$, $c_-^*=\bar{c}_-=c+c^*$,  $r^*(y)=\beta(y)$, $Q_t:=Q[t,\cdot]$, and $Q_\pm:=Q_\pm[1,\cdot]$. Then the following statements are valid:
\begin{itemize}
\item [{\rm (i)}]  $\{Q_t\}^{\infty}_{t=0}$  is a subhomogeneous, continuous-time semiflow on $C_+$  and satisfies (A1-A4), (SP), and (SC).

\item [{\rm (ii)}]  $\{(T_{-z}\circ Q_1)^n[r^*]:n\in \mathbb{N}\}$ is precompact in $C$  for all $z\in \mathbb{R}$.

\item [{\rm (iii)}]    $Q_+$ satisfies  (UC) and (AA).

\item [{\rm (iv)}]    $Q_-$ satisfies (UAA).

\end{itemize}

\end{prop}

\noindent
{\bf Proof.}  Since (i), (ii), and (iv) are obvious, we only verify (iii).

(iii) Applying  \cite[Theorem 5.5]{lz2007} to \eqref{6.5.3},   we know that
any $\varepsilon>0$ and $\phi\in C_+$,
\[
\lim\limits_{t\to \infty} \Big[\sup\{\Phi[t,\phi](x,y):|x|\geq t(c^*+\varepsilon)\mbox{ and } y\in \overline{\Omega}\}\Big]=0 \mbox{ whenever } \phi \mbox{ has a compact support},
\]
\mbox{ and }
\[
\lim\limits_{t\to \infty} \max\{|\Phi[t,\phi](x,y)-r^*|:  |x|\leq t (c^*-\varepsilon)\mbox{ and } y\in \overline{\Omega} \}=0 \mbox{ whenever } \phi\neq 0.
\]
These, together with  Proposition~\ref{prop6.14}-(iii), imply (iii).
 \qed

\

As a straightforward consequence of  Proposition~\ref{prop6.15} and Theorems~\ref{thm5.1} and~\ref{thm5.2}, we have the following result
for system  ~\eqref{6.13}.
\begin{thm}\label{thm6.5}
Assume that  $g$ satisfies (C1)--(C6). Let $c^*=2\sqrt{\lambda_0}$ and $r^*(y)=\beta(y)$.
Then the following statements are valid:
\begin{itemize}
\item [{\rm (i)}] If $c<c^*$, then for any
$\varepsilon\in (0,\min\{c^*-c,c^*\})$
 and $\varphi> 0$, we have $$\lim\limits_{n\rightarrow \infty}
\max\{||P[t,\varphi](x,\cdot)-r^*||_{L^{\infty}(\overline{\Omega})}:t\max\{c+\varepsilon,-c^*+\varepsilon\}\leq x \leq t(c^*-\varepsilon)\}= 0.$$
\item [{\rm (ii)}] If $\varphi \in C_+$ has a compact support, then $\lim\limits_{t\rightarrow \infty}
\Big[\sup\{
||P[t,\varphi](x,\cdot)||_{L^{\infty}(\overline{\Omega})}:x\in   t(\mathbb{R}\setminus [-c^*-\varepsilon,c^*+\varepsilon])\}\Big]=0$ for all $\varepsilon>0$.

\item [{\rm (iii)}]  For any $\varepsilon>0$ and $\varphi \in C_+$, we have
$\lim\limits_{t\rightarrow \infty}
\Big[\sup\{
||P[t,\varphi](x,\cdot)||_{L^{\infty}(\overline{\Omega})}:x\in  (-\infty,t(c-\varepsilon)]\}\Big]=0$.

 \item [{\rm (iv)}] For any $c>-c^{*} $, $\{P[t,\cdot]\}_{t\in\mathbb{R}_+}$ has  a  travelling wave $W(x-ct,\cdot)$ connecting
$0$ to $r^*:=\beta(\cdot)$.
\end{itemize}
\end{thm}

We remark that the forced traveling waves and spreading properties of  reaction-diffusion equations in a cylinder and with a shifting habitat
were studied in \cite{BR2009} and \cite{BG2019}, respectively, under different
assumptions on the reaction term $f(x,y,u)$.  In particular, it was
assumed that $\lim\limits_{r\to \infty}\sup\limits_{|x|>r, y\in \Omega}
\partial_uf(x,y,0)<0$ in \cite{BR2009} and that $\lim_{x\to \pm \infty}
f(x,y,u)=f^{\infty}_{\pm}(u)$ uniformly for $y\in \Omega$ in \cite{BG2019}.

\subsection{The Dirichlet problem for a  time-delayed equation}
In this subsection, we focus on  the following time-delayed reaction-diffusion equation subject to the Dirichlet boundary condition:
\begin {equation} \left\{
\begin{array}{rcl}
	\frac{\partial u}{\partial t}(t,x)&=& d\frac{\partial^2 u}{\partial x^2}(t,x)-\mu u(t,x)+ \mu f(u(t-\tau,x)),
	\quad t,\, x>0,\\
	u(t,x)& =& \varphi(t,x), \quad (t,x)\in [-\tau,0]\times \mathbb{R}_+,\\
	u(t,0)& =& 0, \quad t\in   [-\tau,\infty),
\end{array} \right.
\label{6.10}
\end {equation}
where $d>0$, $\tau>0$ and $\varphi\in \mathcal{C}_+:=\{\psi\in C([-\tau,0]\times\mathbb{R}_+,\mathbb{R}_+)\cap L^\infty([-\tau,0]\times\mathbb{R}_+,\mathbb{R}):\psi(\theta, 0)=0\}$. Without loss of generality, we may assume $\tau\in (0,1)$.
The   reaction term $f\in C^1(\mathbb{R}_+,\mathbb{R}_+)$ satisfies that $f'(0)>1$, $f'(u)\geq 0$ and $f(u)<f'(0)u$ for all $u>0$, and  $f$ has one unique fixed point $u^{*}>0$.

It is well-known that equation~\eqref{6.10} has a
unique mild solution on $[0, \infty)$, denoted by $u^{\varphi}(f,\mu;t,x)$ or $(u^\varphi)_t$, which is
also the classical solution of~\eqref{6.10} on  $(\tau,  \infty)$.
To obtain the propagation dynamics of~\eqref{6.10}, we consider the following
integral equation with the given initial function:
\begin {equation}\label{6.10-1}
\left\{
\begin{array}{rcl}
v(t,\cdot) & =& S_{\mu}(t)[\varphi(0,\cdot)]+\int^t_0
S_{\mu}({t-s})[\mu f(v(s-\tau,\cdot))] \mathrm {d} s, \qquad  t\in
\mathbb{R}_+,
\\
v_0 & = & \varphi\in \mathcal{X}_+.
\end{array}
\right.
\end {equation}
Here $\mathcal{X}_+=\{\varphi\in L^\infty(M\times \mathbb{R},\mathbb{R}):\varphi(\theta,\cdot)\in X_+, \, \forall  \theta \in M \mbox{ and } \sup|\varphi|(M\times \mathbb{R})<\infty$\} with the norm $||\varphi||_{\mathcal{X}}\triangleq
\sup\limits_{\theta\in M}\{||\varphi(\theta,\cdot)||_{X}$\} and
$S_{\mu}(t):X\to X$ is defined by
\begin{equation}\label{eqn6.10-2}
\left\{
\begin{array}{rcl}
S_{\mu}(0)[\phi](x) & = & \phi(x), x\in \mathbb{R},\\
S_{\mu}(t)[\phi](x) & = & \frac{\exp(-\mu t)}{\sqrt{4d\pi t}}\int_0^\infty
\phi(y)\left[\exp\left(-\frac{(x-y)^2}{4 d
t}\right)-\exp\left(-\frac{(x+y)^2}{4 dt}\right)\right]{{\rm{d}}}y, (t,x)\in (0,\infty)\times \mathbb{R}_+,\\
S_{\mu}(t)[\phi](x) & = &0, (t,x)\in (0,\infty)\times (-\infty,0),
\end{array}
\right.
\end{equation}
for any  $\phi\in X$. Note that for any $\phi\in X$ with $\phi((-\infty,0])=0$,
 $S_{\mu}(t)[\phi](x)$ solves the following linear system:
\[
\left\{
\begin{array}{rcll}
   \frac{\partial u}{\partial t}&=&d \frac{\partial^2 u}{\partial x^2}-\mu u, \qquad & t>0 \mbox{ and } x>0,
    \\
    u(t,0) & = & 0,    & t \in \mathbb{R}_+,
    \\
    u(0,x)& =& \phi(x), & x\in \mathbb{R}_+.
\end{array}
\right.
\]
By the method of steps, it follows that for any given $\varphi\in \mathcal{X}_+$, equation~\eqref{6.10-1} has a
unique solution on a maximal interval $[0, \infty)$, denoted by $v^{\varphi}(t,x)$ or $(v^\varphi)_t$.

According to the definitions of  $(u^{\varphi})_t$ and $(v^{\varphi})_t$,  we easily obtain the following relation between ~\eqref{6.10} and ~\eqref{6.10-1}.

\begin{prop} \label{prop6.7-1}
 If  $\varphi\in C_+$ and $\varphi|_{[-\tau,0]\times (-\infty,0]}=0$, then
$v^{\varphi}(t,x)=u^{\varphi|_{[-\tau,0]\times \mathbb{R}_+}}(t,x)$ for all $(t,x)\in \mathbb{R}_+\times \mathbb{R}_+$.

\end{prop}

With  the definitions of $S_\mu(t)$ and $(v^\phi)_t$,  we can verity the following result about the boundedness, positive invariance, semigroup property, continuity, and monotonicity.

\begin {prop} \label{prop6.7} Let  $Q:\mathbb{R}_+\times \mathcal{X}_+\ni (t,\varphi)\mapsto (v^\phi)_t\in \mathcal{X}_+$  and $C_{+,0}:=\{\varphi\in C_+:\varphi([-\tau,0]\times (-\infty,0])\subseteq \{0\}\}$.
Then the following statements are valid:

\begin{itemize}
\item[{\rm (i)}] Let $\varphi \in \mathcal{X}_+$ and $t\in \mathbb{R}_+$. Then $\limsup\limits_{t\to\infty}||Q_t[\varphi]||_{L^\infty([-\tau,0]\times \mathbb {R}, \mathbb {R})}\leq u^*$, and $Q_t[\phi]:=Q[t,\phi]\in \mathcal{X}_+$. Moreover, $Q_t[\varphi]\in C_+ $ and   $Q_t[\varphi](\theta,x)\in Int(\mathbb{R}^N_+) $ for all $(\varphi,t,\theta,x)\in \mathcal{X}_+\times (\tau,\infty)\times M\times (0,\infty)$.
\item[{\rm (ii)}]  $Q_t\circ Q_s=Q_{t+s}$ for all $s,t\in \mathbb{R}_+$.
\item[{\rm (iii)}] Let $\psi,\varphi \in \mathcal{X}_+$ with $\psi-\varphi \in \mathcal{X}_+$. Then
$0\leq v^\varphi(t,x)\leq v^\psi(t,x)$ for all $(t,x)\in\mathbb{R}_+\times \mathbb{R}$.

\item[{\rm (iv)}]  For any $t^*,r>0$,  $Q[t,\cdot]$ is continuous in $ C_r$ uniformly for $t\in [0,t^*]$.  In particular, $Q_t|_{C_r}:C_r\to C_+$ is continuous  for any $r>0$ and $t>\tau$.

\item[{\rm (v)}] $Q[t,C_r]$ is precompact in $C$ for all $(r,t)\in (0,\infty)\times(\tau,\infty)$.

\end{itemize}
\end{prop}

We should point out that $Q$ is not continuous at $(0,\varphi)\in \mathbb{R}_+\times C_+$ whence $\varphi(0,0)>0$. This is because  $Q[t,\psi](0,0)=0$ for all $t>0$ and $\psi\in C_+$.

To continue our study, we introduce two auxiliary  semigroups (for any given $z\in \mathbb{R}$, the former $S_{\mu,z}(t)$  is not continuous at $(0,\phi)\in \mathbb{R}_+\times X$ with $\phi(-z)\neq 0
$).  For any given $z\in \mathbb{R}$, define
$S_{\mu,z}(t), S_{\mu,\infty}(t):X\to X$ by
\begin{equation}\label{eqn6.10-2}
\left\{
\begin{array}{rcl}
S_{\mu,z}(0)[\phi](x) & = & \phi(x), x\in \mathbb{R},\\
S_{\mu,z}(t)[\phi](x) & = & \frac{\exp(-\mu t)}{\sqrt{4d \pi t}}\int_{-z}^\infty
\phi(y)\left[\exp\left(-\frac{(x-y)^2}{4 d
t}\right)-\exp\left(-\frac{(x+y+2z)^2}{4 d t}\right)\right]{{\rm{d}}}y, (t,x)\in (0,\infty)\times [-z,\infty),\\
S_{\mu,z}(t)[\phi](x) & = &0, (t,x)\in (0,\infty)\times (-\infty,-z),
\end{array}
\right.
\end{equation}
and
\begin{equation}\label{eqn6.10-2-1}
\left\{
\begin{array}{rcl}
S_{\mu,\infty}(0)[\phi](x) & = & \phi(x), x\in \mathbb{R},\\
S_{\mu,\infty}(t)[\phi](x) & = & \frac{\exp(-\mu t)}{\sqrt{4 d \pi t}}\int_{\mathbb{R}}
\phi(y)\exp\left(-\frac{(x-y)^2}{4 d
t}\right){{\rm{d}}}y, (t,x)\in (0,\infty)\times \mathbb{R},\\
\end{array}
\right.
\end{equation}
for any  $\phi\in X$.

Now we consider  the following two integral equations with the given initial function:
\begin {equation}\label{6.3.1}
\left\{
\begin{array}{rcl}
v(t,\cdot) & =& S_{\mu,z}(t)[\varphi(0,\cdot)]+\int^t_0
S_{\mu,z}({t-s})[\mu f(v(s-\tau,\cdot))] \mathrm {d} s, \qquad  t\in
\mathbb{R}_+,
\\
v_0 & = & \varphi\in   \mathcal{X}_+,
\end{array}
\right.
\end {equation}
and
\begin {equation}\label{6.3.2}
\left\{
\begin{array}{rcl}
v(t,\cdot) & =& S_{\mu,\infty}(t)[\varphi(0,\cdot)]+\int^t_0
S_{\mu,\infty}({t-s})[\mu f(v(s-\tau,\cdot))] \mathrm {d} s, \qquad  t\in
\mathbb{R}_+,
\\
v_0 & = & \varphi\in  \mathcal{X}_+
\end{array}
\right.
\end {equation}
It is easy to see that for all $t>\tau$, the solution of \eqref{6.3.2}  is also a classical solution of  the following equation
\begin {equation} \left\{
\begin{array}{rcl}
	\frac{\partial u}{\partial t}(t,x)&=& d\frac{\partial^2 u}{\partial x^2}(t,x)-\mu u(t,x)+ \mu f(u(t-\tau,x)),
	\quad (t,x)\in (0,\infty)\times \mathbb{R},\\
	u(t,x)& =& \varphi(t,x), \quad (t,x)\in [-\tau,0]\times \mathbb{R}.
\end{array} \right.
\label{6.3.2001}
\end {equation}
Let $Q[t,\phi;z]$  and $\Phi[t,\phi]$ be the solutions of~\eqref{6.3.1} and~\eqref{6.3.2} with the initial value $v_0= \phi\in  \mathcal{X}_+$, respectively.
By straightforward computations, we then have  the following result.
\begin{prop} \label{prop6.8} Let $t\in \mathbb{R}_+$ and $\phi\in  \mathcal{X}_+$. Then the following statements are valid:
	\begin{itemize}
		\item [{\rm (i)}]  $Q[t,\phi;z](\theta,x)=Q[t,T_z[\phi]](\theta,x+z)$ and $Q[t,\phi;z](\theta,x)\leq Q[t,\phi;\tilde{z}](\theta,x)$ for all $(\theta,x,z,\tilde{z})\in [-\tau,0]\times \mathbb{R}^3$ with $z\leq \tilde{z}$.
		
		\item [{\rm (ii)}]  $ \Phi[t,\phi]=\lim\limits_{z\to  \infty} Q[t,\phi;z],$ that is,  $Q[t,\phi;z]\to  \Phi[t,\phi]$ with respect to the compact open topology  as $z\to  \infty$.

		\item [{\rm (iii)}]  $ \lim\limits_{z\to - \infty} Q[t,\phi;z]=0$, that is,  $Q[t,\phi;z]\to  0$ with respect to the compact open topology  as $z\to - \infty$.
	\end{itemize}

\end{prop}

\noindent

\begin{prop} \label{prop6.9}  Let $c^*>0$ be the spreading speed  of system \eqref{6.3.2001} (see \cite{lz2007}).
 For a given number $t_0>\tau$,
 let  $c^*_+=\bar{c}_+=c_-^*=\bar{c}_-=t_0 c^*$,  $r^*=u^*$, $Q_t:=Q[t,\cdot]$, $Q_-:=0$, and $Q_+:=\Phi[t_0,\cdot]$ for all $t\in \mathbb{R}_+$. Then the following statements are valid:
\begin{itemize}
\item [{\rm (i)}]  $Q_\pm$ satisfies (A4).

\item [{\rm (ii)}]    $Q_+$ satisfies  (UC) and (AA).
\item [{\rm (iii)}]    $Q_-$ satisfies (UAA).
\item [{\rm (iv)}]   $Q_{t_0}|_{C_+}$ satisfies (ASH-UC-SP).
\end{itemize}
\end{prop}

\noindent
{\bf Proof.}  Since (i) and (iii) are obvious, we only prove (ii) and (iv).

(ii)  follows from the proof of Proposition~\ref{prop6.3}-(iii) with $k(\cdot)=\delta(\cdot)$.

(iv)  Take  $\gamma_l\in (0,  f'(0)-1)$ with $\lim\limits_{l\to \infty} \gamma_l=0$ and $\gamma_l>\gamma_k$ for all  positive integers $k>l$. Then there exists $\delta_l>0$ such that $f(u)\geq (f'(0)-\gamma_l) u$ for all $u\in [0,\delta_l]$. Let us define
\[
f_l(u)=\left\{
\begin{array}{ll}
(f'(0)-\gamma_l)(u-\frac{u^2}{2\delta_l}), & (s,u)\in \mathbb{R}\times  [0,\delta_l],
 \\
\frac{f'(0)-\gamma_l}{2}\delta_l, \qquad & (s,u)\in \mathbb{R}\times (\delta_l,\infty).
\end {array}
\right.
\]
and \[
r^*_l=\left\{
\begin{array}{ll}
\frac{f'(0)-\gamma_l}{2}\delta_l, & f'(0)>2+\gamma_l,
 \\
2\frac{f'(0)-\gamma_l-1}{f'(0)-\gamma_l}\delta_l, \qquad & f'(0)\leq 2+\gamma_l.
\end {array}
\right.
\]
Then $f(u)\geq f_l(u)$ for all $u\in\mathbb{R}_+$ and $f(r^*_l)=r_l^*$.
Let $Q_l[t,\phi]=Q[f_l;t,\phi]$, and let $c_{\pm l}^*$ be the spreading speed  of  \eqref{6.1.3} with $f_l$.  Then (iv) follows from   Propositions~\ref{prop6.7}, ~\ref{prop6.8}-(ii), and  ~\ref{prop6.9}-(i,ii). \qed

\

As a consequence of  Remark~\ref{rem4.1}, Theorem~\ref{thm4.1}, Theorem~\ref{thm4.2}, and Corollary~\ref{cor4.1}, we have  the following result for  system  ~\eqref{6.10-1}.
\begin{thm}\label{thm6.3}
 Let $r^*=u^*$ and $c^*>0$ be the spreading speed  of
  system \eqref{6.3.2001}. Then the following statements are valid:
	\begin{itemize}
		\item [{\rm (i)}]  For any $\varepsilon\in (0,c^*)$ and $\varphi> 0$, we have $$\lim\limits_{n\rightarrow \infty}
		\max\{|Q[t,\varphi](\theta, x)-r^*|:\theta\in [-\tau,0]\mbox{ and }t\varepsilon\leq x \leq t(c^*-\varepsilon)\}= 0.$$
		\item [{\rm (ii)}] If $\varphi
		$ has a  compact support, then $\lim\limits_{t\rightarrow \infty}
		\Big[\sup\{
		Q[t,\varphi](\theta,x):(\theta,x)\in [-\tau,0]\times t(\mathbb{R}\setminus [-\varepsilon,c^*+\varepsilon])\}\Big]=0$ for all $\varepsilon>0$.
		
		\item [{\rm (iii)}]  $\{Q[t,\cdot]\}_{t\in\mathbb{R}_+}$ has  a nontravial fixed point $W(x)$ connecting
		$0$ to $r^*$.
		
		\item [{\rm (iv)}]   For any $\varepsilon\in (0,c^*)$ and $\varphi> 0$, we have $$\lim\limits_{t\rightarrow \infty}
		\max\{|Q[t,\varphi](\theta, x)-W(x)|:\theta\in [-\tau,0]\mbox{ and } x \leq t(c^*-\varepsilon)\}= 0.$$
	\end{itemize}
	\end{thm}

In view of Theorem \ref{thm6.3} and  Proposition \ref{prop6.7-1}, we  have the following result.

\begin{thm}\label{thm6.3-1}
 Let  $c^*>0$ be the spreading speed  of  system \eqref{6.3.2001}.
  Then system  \eqref{6.10} has
a unique nontrivial steady state $W(x)$ connecting
		$0$ to $u^*$ such that  for any  $\varepsilon\in (0,c^*)$ and $\varphi\in \mathcal{C}_+\setminus \{0\}$, there holds $\lim\limits_{n\rightarrow \infty}
		\max\{|u^\varphi(t,x)-W(x)|:\theta\in [-\tau,0]\mbox{ and } x \leq t(c^*-\varepsilon)\}= 0.$	
\end{thm}

In the case where  $\tau=0$, we  can remove the assumptions  that  $f'(u)\geq 0$ for all $u>0$ to obtain the following result.
\begin{cor}\label{cor6.5}
System  \eqref{6.10} has a unique nontrivial steady state $W(x)$ connecting
		$0$ to $u^*$ such that  for any  $\varepsilon\in (0,2\sqrt{\mu [f'(0)-1]})$ and $\varphi\in BC(\mathbb{R_+,\mathbb{R}_+})\setminus \{0\}$ with $\varphi(0)=0$, there holds $\lim\limits_{n\rightarrow \infty}
		\max\{|u^\varphi(t,x)-W(x)|:\theta\in [-\tau,0]\mbox{ and } x \leq t(2\sqrt{\mu [f'(0)-1]}-\varepsilon)\}= 0.$		
		\end{cor}
\noindent
{\bf Proof.}   Fix  $\varepsilon\in (0,2\sqrt{\mu [f'(0)-1]})$ and $\varphi\in BC(\mathbb{R_+,\mathbb{R}_+})\setminus \{0\}$ with $\varphi(0)=0$. Let $M=\max\{u^*,||\varphi||_{L^\infty(\mathbb{R}_+,\mathbb{R})}\}$,  $k:=k_M=\sup\{|f'(u)|:u\in [0,M]\}>0$, $\tilde{M}=\inf\{ u\geq M:f'(u)+k= 0\}$, $\mu_k=\mu(1+k)\geq M$,
and
 \[
f_k(u)=\left\{
\begin{array}{ll}
\frac{f(u)+ku)}{1+k}, & u\in  [0,\tilde{M}),
 \\
\frac{f(\tilde{M})+k\tilde{M})}{1+k}, & u\in  [\tilde{M},\infty).
\end {array}
\right.
\]
By applying Theorem~\ref{thm6.3-1} to $f=f_k,\mu=\mu_k$ and $\tau=0$, we know that \eqref{6.3} has a unique nontrivial steady state $W_k(x)$ connecting
		$0$ to $u^*$ such that
		$$\lim\limits_{n\rightarrow \infty}
		\max\{|u^\varphi(f_k,\mu_k;t,x)-W_k(x)|:\theta\in [-\tau,0]\mbox{ and } x \leq t(2\sqrt{\mu_k [f_k'(0)-1]}-\varepsilon)\}= 0.$$
This, together with the fact that $u^\psi(f,\mu;t,x)=u^\psi(f_k,\mu_k;t,x)$ for all $(t,x,\psi)\in \mathbb{R}_+\times \mathbb{R}_+\times C(\mathbb{R}_+,[0,M])$ with $\psi(0)=0$, yields the desired  result. \qed

\

We note that Corollary~\ref{cor6.5}  was  obtained in \cite{yc2017}  by using  the iteration method of travelling wave maps.

\subsection{ A KPP-type equation  in  spatially inhomogeneous media}

Consider the following asymptotically homogeneous KPP-type equation:
\begin {equation}\label{6.26}
\left\{
\begin{array}{ll}
&\frac{\partial u}{\partial t}  =  d\frac{\partial^{2}u}{\partial x^{2}} + h(x,u),~~~~~(x,t)\in\mathbb{R}\times\mathbb{R}_{+},\\
&u(0,\cdot) = \phi\in C(\mathbb{R},\mathbb{R}_{+})\cap L^{\infty}(\mathbb{R},\mathbb{R}),
\end{array} \right.
\end {equation}
where $d>0$ and $h\in C(\mathbb{R}\times\mathbb{R}_{+},\mathbb{R})$. Assume that
\begin {enumerate}
\item[(D1)] $h(s,0)=0\mbox{~for~all~}s\in\mathbb{R}$;
\item[(D2)] $\mbox{There~exists~}M^{*}>0\mbox{~such~that~}h(s,u)\leq0\mbox{~for~all~}(s,u)\in\mathbb{R}\times[M^{*},\infty)$;
\item[(D3)] $h(s,\cdot)\rightarrow h_{\pm}^{\infty}(\cdot)\mbox{~in~}
C_{loc}^{1}(\mathbb{R}_{+},\mathbb{R})\mbox{~as~}s\rightarrow\pm \infty$,
where  $h_{-}^{\infty}(\cdot)$ and  $h_{+}^{\infty}(\cdot)$  are two KPP-type functions with $u_{-}^{*}$  and $u_{+}^{*}\in(0,\infty)$, respectively.
\end{enumerate}

We introduce the following auxiliary KPP-type equations:
\begin {equation}\label{6.28}
\begin{array}{ll}
\frac{\partial u}{\partial t}  &=  d\frac{\partial^{2}u}{\partial x^{2}} + h_{\pm}^{\infty}(u).
\end{array}
\end {equation}
It is well-known that for any given $\phi\in C_+:=C(\mathbb{R},\mathbb{R}_+)\cap L^\infty(\mathbb{R},\mathbb{R})$, equation~\eqref{6.26} has a
unique solution on its maximal interval $[0, \eta_{\phi})$, denoted by $u^{\phi}(t,x)$, with $\limsup\limits_{t\rightarrow \eta_\phi^-
}||u^\phi(t,x)||=\infty$ whenever $\eta_\phi<\infty$.
In order to emphasize the dependence on nonlinear reaction terms, we also use  $u^{\phi}(t,x;h)$ and $u^{\phi}(t,x;h_\pm^\infty)$ to
represent the solutions of the initial value problem of \eqref{6.26} and \eqref{6.28}, respectively.
By the standard arguments, we have the following result.

\begin {prop} \label{prop6.20} Let  $Q:[0,\eta_\phi)\times C_+\ni (t,\varphi)\mapsto u^\phi(t,\cdot)\in C:=C(\mathbb{R},\mathbb{R})\cap L^\infty(\mathbb{R},\mathbb{R})$.
Then the following statements are valid:
\begin{itemize}
\item[{\rm (i)}] Let $\phi \in C_+$ and $t\in \mathbb{R}$. Then
$\eta_\phi=\infty$,  $\limsup\limits_{t\to
\infty}||u^\phi(t,\cdot)||_{L^\infty(\mathbb{R},\mathbb{R})}\leq M^*$, and $Q_t[\phi]:=Q[t,\phi]\in C_+$.
Moreover, $Q_t[\phi]\in C_+^\circ:=\{\phi\in C_+:\phi(\mathbb{R})\subseteq (0,\infty)\}$ for all $(t,\phi)\in(0,\infty)\times(C_+\backslash\{0\})$.

\item[{\rm (ii)}] Let $\psi,\phi \in C_+$ with $\phi\leq \psi $. Then
$0\leq u^\phi(t,x)\leq u^\psi(t,x)$ for all $(t,x)\in
\mathbb{R}_+\times \mathbb{R}$. Moreover, if $\phi<\psi$, then $u^\phi(t,x)<u^\psi(t,x)$ for all $(t,x)\in (0,\infty)\times\mathbb{R}$.

\item[{\rm (iii)}]  $Q_t[C_r]\subseteq C_r$ and $Q_t\big|_{C_{r}}: C_r\rightarrow C_r$ is a compact semiflow, where $t>0$, $r\geq M^*$ and $C_r:=\{\phi\in C_+:\phi\leq r\}$ equipped with
the compact open topology.
\end{itemize}
\end{prop}

By the arguments similar to those in \cite{ycw2019}, we can prove the following
two technical results on $h(x,u)$ and its limiting functions $h_{\pm}^{\infty}(u)$.

\begin{lemma}\label{lem6.2}
For any $M\geq M^{*}$ and $\gamma>0$, there exist $l_{\pm}=l_{\pm}(M,\gamma)$, $L_{\pm}=L_{\pm}(M,\gamma)$, $\mathcal{S}=\mathcal{S}(M,\gamma)>0$ such that
for any $(s,u)\in\mathbb{R}\times[0,M]$, there holds
$$(\frac{\partial h_{\pm}^{\infty}}{\partial u}(0)-\gamma)u-l_{\pm}u^{2}\leq h(s,u)\leq(\frac{\partial h_{\pm}^{\infty}}{\partial u}(0)+\gamma)u-L_{\pm}u^{2},$$
for all $\pm s\geq\mathcal{S}$.
\end{lemma}

\begin{lemma}\label{lem6.3}
For any $M\geq M^{*}$ and $\gamma>0$, there exist $k^*\leq 0$, $r^{\pm}=r_{\gamma,M}^{\pm}(\cdot)\in C(\mathbb{R},\mathbb{R})$, and $K_{\pm}=K_{\pm}(\gamma,M)>0$ such that
\begin{itemize}
\item[{\rm (i)}] $r^{\pm}$ is nondecreasing with  $r^{\pm}(s)=\frac{k^*}{K_\pm}\leq0$ for all $s\leq 0$ and $r^{\pm}(+\infty)=\frac{\frac{{\rm d}h_{\pm}^{\infty}(0)}{{\rm d}u}-\gamma}{K_{\pm}}$;

\item[{\rm (ii)}] $h(s,u)\geq \max\{K_{+}u(r^{+}(s)-u),K_{-}u(r^{-}(- s)-u)\}$  for all $(s,u)\in \mathbb{R}\times[0,M]$.
\end{itemize}
\end{lemma}

We  define $R_{\pm}=R_{\pm}^{\gamma,M}:\mathbb{R}\times\mathbb{R}_{+}\rightarrow\mathbb{R}$ by
\begin{equation}\label{Rfunctions}
R_{\pm}(s,u)=K_{\pm}(r^{\pm}(s)u-u^{2}) \mbox{ for all } (s,u)\in\mathbb{R}\times[0,M].
\end{equation}
It then follows that
$$
h(s,u)\geq R_{+}(s,u) \mbox{ and } h(s,u)\geq R_{-}(-s,u)  \mbox{ for all } (s,u)\in\mathbb{R}\times[0,M].
$$
Based on Lemma~\ref{lem6.2}, we are able to prove  the following  result on the extinction behavior of solutions.
\begin{prop}\label{prop6.19}
Assume that  $h(x,u)$  satisfies (D1)--(D3) and let $c_{-}^*=2\sqrt{d \frac{{\rm d}h_{-}^{\infty}(0)}{{\rm d}u}}$ and $c_{+}^*=2\sqrt{d \frac{{\rm d}h_{+}^{\infty}(0)}{{\rm d}u}}$. Then the following statements are
valid:
\begin{itemize}
\item[{\rm (i)}]
If $\phi\in C_+$ with $\phi \equiv0$ for all sufficiently large negative $x$, then for any  $\varepsilon >0$,
\[
\lim_{t\to  \infty}\left[\sup_{x\leq -t(c_{-}^*+\varepsilon)}u^\phi(t,x)\right]=0.
\]
\item[{\rm (ii)}]
If $\phi\in C_+$ with $\phi \equiv0$ for all sufficiently large positive $x$, then for any  $\varepsilon >0$,
\[
\lim_{t\to  \infty}\left[\sup_{x\geq t(c^*_++\varepsilon)}u^\phi(t,x)\right]=0.
\]
\end{itemize}
\end{prop}

\noindent \textbf {Proof.}
We only prove (i) since the proof of (ii) is similar.
Fix $\varepsilon>0$ and $\phi\in C_{+}$ with $\phi\equiv0$ for all sufficiently large negative $x$.
According (D1)--(D3) and Proposition \ref{prop6.20}-(i), we easily see that there exist $M_{0}\geq M^{*}$, $\xi_{0}>0$ and $\delta_{0}>0$ such that $0\leq u^{\phi}(t,x)\leq M_{0}$
for all $(t,x)\in \mathbb{R}_{+}\times\mathbb{R}$ and
$$h(s,u)\leq (\frac{{\rm d}h_{-}^{\infty}(0)}{{\rm d}u}+\delta_{0})u:=\frac{(c_{-}^{*}+\frac{\varepsilon}{3})^{2}}{4d}u$$
for all $(s,u)\in(-\infty,-\xi_{0}]\times [0,M_{0}]$.

By using  \eqref{6.26} and taking $r^{*}=\frac{{\rm d}h_{-}^{\infty}(0)}{{\rm d}u}+\delta_{0}$, we have
$$u^{\phi}(t,x)=e^{r^{*}t}T(t)[\phi](x)+\int_{0}^{t}e^{r^{*}(t-s)}T(t-s)\big[ h(\cdot,u^{\phi}(s,\cdot)) - r^{*}u^{\phi}(s,\cdot)  \big](x){\rm d}s$$
for all $(t,x)\in\mathbb{R}_{+}\times\mathbb{R}$, where
\begin {equation*}
T(t)[\phi](x)= \left\{
\begin{array}{rcl}
&\phi(x),~~~~~&t=0,\\
&\frac{1}{\sqrt{4\pi dt}}\int_{\mathbb{R}}\phi(y)\exp(-\frac{(x-y)^{2}}{4dt}){\rm d}y,~~~~~&t>0.
\end{array} \right.
\end {equation*}
Let $M^{**}=\sup\{ h(s,u):(s,u)\in\mathbb{R}\times[0,M_{0}] \}\in (0,\infty)$.
It then  follows that for any $(t,x)\in\mathbb{R}_{+}\times\mathbb{R}$ with $x\leq-t(c_{-}^{*}+\varepsilon)\leq -\frac{\varepsilon t}{10}\leq-4\xi_{0}$,
\begin{align*}
u^{\phi}(t,x)&\leq e^{r^{*}t}T(t)[\phi](x)+\int_{0}^{t}e^{r^{*}(t-s)}
\int_{-\xi_{0}}^{\infty}  \frac{1}{\sqrt{4\pi d(t-s)}}\big[ h(y,u^{\phi}(s,y)) - r^{*}u^{\phi}(s,y)  \big] \exp (-\frac{(x-y)^{2}}{4d(t-s)}){\rm d}y{\rm d}s\\
&\leq e^{r^{*}t}T(t)[\phi](x)+\int_{0}^{t}   \frac{M^{**}e^{r^{*}s-\frac{x^{2}}{4ds}}}{\sqrt{4\pi ds}}
\int_{-\xi_{0}}^{\infty}  e^{    -\frac{y^{2}-2xy}{4ds}    }{\rm d}y{\rm d}s\\
&\leq e^{r^{*}t}T(t)[\phi](x)+ M^{**}\int_{0}^{t}e^{  r^{*}s-\frac{x^{2}}{4ds}  }{\rm d}s  +
\frac{M^{**}}{|x|} \sqrt{\frac{dt}{\pi}}\int_{0}^{t}e^{  r^{*}s-\frac{(x+\xi_0)^{2}}{4ds}  }{\rm d}s\\
&\leq e^{r^{*}t}T(t)[\phi](x)+ M^{**}\int_{0}^{t}e^{  r^{*}s-\frac{   (c_-^*+\frac{2}{3}\varepsilon)^{2}t^2  }{4ds}  }{\rm d}s  +
\frac{M^{**}}{c_-^*+\varepsilon} \sqrt{\frac{d}{\pi t}}\int_{0}^{t}e^{  r^{*}s-\frac{   (c_-^*+\frac{2}{3}\varepsilon)^{2}t^2   }  {4ds}  }{\rm d}s \\
&\leq e^{r^{*}t}T(t)[\phi](x)+ M^{**} \left[ 1+\frac{1}{c_-^*+\varepsilon} \sqrt{\frac{d}{\pi t}}  \right]
\int_{0}^{t}e^{  r^{*}s-\frac{   r^* t^2+\frac{\varepsilon^{2}t^2}{36d}  }{s}  }{\rm d}s  \\
&\leq e^{r^{*}t}T(t)[\phi](x)+ M^{**} \left[ 1+\frac{1}{c_-^*+\varepsilon} \sqrt{\frac{d}{\pi t}}  \right] t e^{ -\frac{ \varepsilon^{2}  } {36d} t}.
\end{align*}
This, together with the fact that $\lim\limits_{t\rightarrow\infty}\big[ \sup\limits_{x\leq-t(2\sqrt{dr^{*}}+\frac{2\varepsilon}{3})} e^{r^{*}t}T(t)[\phi](x) \big]=0$ (see \cite[Proposition 2.3-(ii)]{ycw2019}), yields the desired statement (i). \qed

\

Applying Lemma \ref{lem6.3}, the definitions of $R_{\pm}$, and the standard comparison technique, we easily obtain the following result.

\begin {prop} \label{prop6.19.1}
If $\gamma>0$, $M\geq M^{*}$ and $\phi\in C_{+}$, then the following statements are valid:
 \begin{itemize}
\item[{\rm (i)}] $u^{\phi}(h;t,x)\geq u^{\phi}(R_{+};t,x)$ for all $(t,x)\in \mathbb{R}_+\times \mathbb{R}$;

\item[{\rm (ii)}] $u^{\phi}(h;t,x)\geq u^{\phi}(R_{-}(-\cdot,\cdot);t,x)=u^{\phi(-\cdot)}(R_{-};t,-x)$ for all $(t,x)\in \mathbb{R}_+\times \mathbb{R}$,
\end{itemize}
where $R_{\pm}:=R_{\pm}^{\gamma,M}$ are defined as in \eqref{Rfunctions}.
\end{prop}

\begin {prop} \label{prop6.19.2}
If $0<\gamma<\min\{ \frac{{\rm d}h_{+}^{\infty}(0)}{{\rm d}u},\frac{{\rm d}h_{-}^{\infty}(0)}{{\rm d}u}\}$ and $M\geq M^{*}$,
then there exist $W_{\pm}:=W_{\pm}^{\gamma,M}\in C_{+}^{\circ} $ such that

 \begin{itemize}
\item[{\rm (i)}] $W_{\pm} (\infty)= \frac{\frac{{\rm d}h_{\pm}^{\infty}(0)}{{\rm d}u} - \gamma}{K_\pm} $ and $W_{\pm}(-\infty)=0$;

\item[{\rm (ii)}] $W_{+}$ and $W_{-}$ are nondecreasing;

\item[{\rm (iii)}] for any $\phi\in C_{+}\backslash \{0\}$ and $\varepsilon>0$,
$$\lim_{t\rightarrow\infty}\Big[\sup\big\{  |u^{\phi}(R_{\pm};t,x)-W_{\pm}(x)|:~
x\leq( 2\sqrt{d\left(\frac{{\rm d}h_{\pm}^{\infty}(0)}{{\rm d}u} - \gamma\right)} -\varepsilon)t  \big\}\Big]=0.$$
In particular, $u^{\phi}(R_{\pm};t,\cdot)\rightarrow W_{\pm}(\cdot)$  with respect to the compact open topology as $t\rightarrow\infty$,
where $R_{\pm}:=R_{\pm}^{\gamma,M}$ are defined as in \eqref{Rfunctions}.
\end{itemize}
\end{prop}

\noindent {\bf Proof.} We only consider the case of $"+"$ since the case of $"-"$ is similar.
Applying Corollary \ref{cor6.3}-(iv) with $f\equiv R_{+}$ and $c=0$, we can obtain the existence of $W_{+}$ with (i) and (ii).
Furthermore, we have $W_{+}=\lim\limits_{t\rightarrow\infty}u^{\phi}(R_{+};t,\cdot)$ with $\phi\equiv\frac{\frac{{\rm d}h_{+}^{\infty}(0)}{{\rm d}u}-\gamma}{K_+}$.
By Theorem \ref{thm4.1}-(v) and Proposition \ref{prop6.3} with $f\equiv R_{+}$ and $c=0$, we have
$$\lim_{\alpha\rightarrow\infty}\Big[  \sup\big\{   | u^{\phi}(R_{+};t,x)-\frac{\frac{{\rm d}h_{+}^{\infty}(0)}{{\rm d}u}-\gamma}{K_+}|:
t\geq\alpha\mbox{~and~}\alpha\leq x\leq t(2\sqrt{d(\frac{{\rm d}h_{+}^{\infty}(0)}{{\rm d}u}-\gamma)}-\varepsilon)  \big\}   \Big]=0$$
for any $\varepsilon>0$ and $\phi\in C_{+}\backslash\{0\}$.

Again, by Corollary \ref{cor6.3}-(iii) with $f\equiv R_{+}$ and $c=0$, we have
$$\lim_{t\rightarrow\infty}  \sup\big\{   | u^{\phi}(R_{+};t,x)|:x\in(-\infty,-t\varepsilon) \big\} =0  $$
for any $\phi\in C_{+}$ and $\varepsilon>0$.
It then follows that $W(+\infty)=\frac{\frac{{\rm d}h_{+}^{\infty}(0)}{{\rm d}u}-\gamma}{K_+}$ and $W(-\infty)=0$ provided that $W$ is nontrivial steady state of $u^{\phi}(R_{+};t,x)$.

Now we claim that $u^{\phi}(R_{+};t,x)$ has the unique steady state solution in $C_{+}\backslash\{0\}$.
Otherwise, there exists $W\in C_{+}^{\circ}:=C^+\cap C(\mathbb{R},(0,\infty))$ such that
$$W(-\infty)=0,\,  W(+\infty)=\frac{\frac{{\rm d}h_{+}^{\infty}(0)}{{\rm d}u}-\gamma}{K_+}, \, W\leq W_{+}, \,  W\neq W_+.
$$
Let $I(\alpha)=\big\{ \beta\geq1:~W_{+}(x)\leq \beta W(x)\mbox{ for all } x\in(\alpha,\infty) \big\}$ and $\beta(\alpha)=\inf I(\alpha)
\mbox{ for all } \alpha\in\mathbb{R}\cup\{-\infty\}$.
Then $I(\alpha_{1})\subseteq I(\alpha_{2})$, $\beta(\alpha_{1})\geq \beta(\alpha_{2})$, and $I(\alpha_1)=[\beta(\alpha_1),\infty)$ for all
$\infty>\alpha_{2}\geq\alpha_{1}\geq-\infty$.
We proceed with two cases.

{\it Case 1.} There exists $\widetilde{\alpha}\in(-\infty,0)$ such that $\beta(0)\notin I_{\widetilde{\alpha}}$.
Let $$\alpha^{*}=\inf\{ \alpha\in\mathbb{R}:~\beta(0)\in I(\alpha) \}.$$
Then $0\geq \alpha^{*}>\widetilde{\alpha}>-\infty$. According to the definitions $I(\cdot)$ and $\beta(\cdot)$, we know that
$$W_{+}(x)\leq\beta(0)W(x)\mbox{~~~for~all~}x\in[\alpha^{*},\infty)$$
and
$$W_{+}(\alpha^{*})=\beta(0)W(\alpha^*),~~W'_+(\alpha^+)\leq\beta(0)W'(\alpha^*).$$
It follows from the definition of $R_+$ that
\begin {equation*}
\left\{
\begin{array}{ll}
&dW''_+(x)+k^*W_+(x)-K_+W_{+}^{2}(x)=0,~~~~~x\leq0,\\
&dW''(x)+k^*W(x)-K_{+}W^{2}(x)=0,~~~~~x\leq0,\\
\end{array} \right.
\end {equation*}
and hence,
\begin{align*}
&\int_{-\infty}^{\alpha^*}[dW''_++k^*W_+(x)-K_+W_{+}^{2}(x)]W(x){\rm d}x\\
&=\int_{-\infty}^{\alpha^*}[dW''(x)+k^*W(x)-K_+W^{2}(x)]W_+(x){\rm d}x=0.
\end{align*}
Thus, we have
\begin{align*}
&\int_{-\infty}^{\alpha^*}K_+W_{+}(x)W(x)[W_+(x)-W(x)]{\rm d}x\\
&=d\int_{-\infty}^{\alpha^*}[W_{+}''(x)W(x)-W''(x)W_+(x)]{\rm d}x\\
&=d[W'_{+}(\alpha^*)W(\alpha^*)-W'(\alpha^*)W_+(\alpha^*)]\leq 0,
\end{align*}
which, together with $W_+(\cdot)\geq W(\cdot)$ and $W,W_+ \in C_+^{\circ}$, implies $W_+(\cdot)\equiv W(x)$, a contradiction.

{\it Case 2.} $\beta(0)\in I(\alpha)$ for all $\alpha\in[-\infty,0]$.
In this case, $I(\alpha)=[\beta(0),\infty)\mbox{ for all }\alpha\in [-\infty,0]$, and hence, $\beta(0)W(x)\geq W_+(x)\mbox{ for all } x\in \mathbb{R}$.
Since  $W_+\gneq W$, we have $\beta(0)>1$.
We further claim that $W_+\big|_{(-\infty,0]}\neq\beta(0)W\big|_{(-\infty,0]}$. Otherwise, $W_+(x)=\beta(0)W(x)$ for all $x\in(-\infty,0]$.
This, together with the definition of ${R}_+$, gives rise to
\begin {equation*}
\left\{
\begin{array}{ll}
&dW''(x)+k^*W(x)-K_+W^{2}(x)=0,~~~~~x\in(-\infty,0],\\
&d\beta(0)W''(x)+k^*\beta(0)W(x)-K_{+}\beta^2(0)W^{2}(x)=0,~~~~~x\in(-\infty,0].\\
\end{array} \right.
\end {equation*}
Then we obtain $K_+W^2\big|_{(-\infty,0]}=0$, a contradiction. It follows that there exists $x^*\in(-\infty,0]$ such that $W_+(x^*)<\beta(0)W(x^*)$,
which, together with the definition of $\mathbb{R}_+$, $W_+(\infty)<\beta(0)W(\infty)$ and strong maximum principle, implies that $W_+(x)<\beta(0)W(x)\mbox{ for all }x\geq x^*$.
Since  $W_+(\infty)=W(\infty)>0$, we have $W_+\big|_{\mathbb{R}_+}\leq (\beta(0)-\delta)W\big|_{\mathbb{R}_+}$ for some $\delta\in (0,\beta(0)-1)$,
which  contradicts the definition of $\beta(0)$.

By  the uniqueness of steady state solution and Corollary \ref{cor4.1}, it then follows that statement (iii)
holds true.  \qed

\begin {prop} \label{prop6.19.3}
If $h(x,u)\not\equiv h(0,u)$, then there exists $W\in C_{+}^{\circ} $ such that $u^{W}(t,x)=W(x)$ for all $(t,x)\in \mathbb{R}_{+}\times\mathbb{R}$.
\end{prop}

\noindent {\bf Proof.}
Let $M\geq M^*$, $\gamma=\frac{1}{3}\min\big\{  \frac{{\rm d}h_{+}^{\infty}(0)}{{\rm d}u},\frac{{\rm d}h_{-}^{\infty}(0)}{{\rm d}u} \big\}$, and let $\widetilde{W}=W_{+}^{\gamma,M}$
be defined as in Proposition \ref{prop6.19.2}.
It follows from Proposition \ref{prop6.19.1}-(i) that
$$
M^*\geq u^{M^*}(h;,t,x)\geq u^{\widetilde{W}}(R_+;t,x)\equiv\widetilde{W}(x) \mbox{ for all } (t,x)\in\mathbb{R}_+\times\mathbb{R}.
$$
Thus,  $u^{M^*}(h;,t,\cdot)\rightarrow W$ in $C_{loc}(\mathbb{R},\mathbb{R})$, and hence in $C_{loc}^{2}(\mathbb{R},\mathbb{R})$, which implies
that $M^*\geq W(x)\geq\widetilde{W}(x)>0 \mbox{ for all } x\in\mathbb{R}$ and $u^{W}(t,x)=W(x)$ for all $(t,x)\in \mathbb{R}_{+}\times\mathbb{R}$.
\qed

\begin {prop} \label{prop6.19.4}
Assume that $h(x,u)\not\equiv h(0,u)$, and let $\mathcal{E}$ be the set of  all steady states in $C_{+}\backslash \{0\}$.  Then the following statements are
valid:
\begin{itemize}
\item[{\rm (i)}] $\emptyset\neq\mathcal{E}\subseteq C_{+}^{\circ}\cap C(\mathbb{R},[0,M^{*}])$;

\item[{\rm (ii)}] $\phi\geq W_{\pm}^{\gamma,M^*}(\pm\cdot)$ for all $\phi\in\mathcal{E}$ and
$\gamma\in(0,\min\{ \frac{{\rm d}h_{+}^{\infty}(0)}{{\rm d}u},\frac{{\rm d}h_{-}^{\infty}(0)}{{\rm d}u} \})$;

\item[{\rm (iii)}] $\phi(\pm\infty)=u_{\pm}^*$ for all $\phi\in\mathcal{E}$;

\item[{\rm (iv)}] $\inf\mathcal{E}$, $\sup\mathcal{E}\in\mathcal{E}$, where $\inf\mathcal{E}(x):=\inf\{\phi(x):\phi\in\mathcal{E} \}$ and
$\sup\mathcal{E}(x):=\sup\{\phi(x):\phi\in\mathcal{E} \}$ for all $x\in\mathbb{R}$;

\item[{\rm (v)}] for any $\phi\in C_{+}\backslash \{0\}$ and $\varepsilon>0$,
there holds
$$\lim_{t\rightarrow\infty}\Big[\max\big\{   \inf\{ |u^{\phi}(h;t,x)-\widetilde{u}|:~\widetilde{u}\in [\inf\mathcal{E}(x),\sup\mathcal{E}(x)]\}:
(-c_{-}^{*}+\varepsilon)t\leq x\leq(c_{+}^{*}-\varepsilon)t \big\}\Big]=0,$$ where $c_-^*$ and $c_+^*$ are defined as in Proposition~\ref{prop6.19}.
\end{itemize}
\end{prop}

\noindent {\bf Proof.} (i) follows from Propositions \ref{prop6.20}-(i) and \ref{prop6.19.3}.

(ii) By Proposition \ref{prop6.19.2}-(iv), $u^{\phi(\pm\cdot)}(R_{\pm}^{\gamma,M^*}(\pm\cdot);t,\cdot)\rightarrow
W_{\pm}^{\gamma,M^*}(\pm\cdot)$ in $L_{loc}^{\infty}(\mathbb{R},\mathbb{R})$ as $t\rightarrow\infty$.
Proposition \ref{prop6.19.1} shows that $u^{\phi}(t,\cdot)\geq u^{\phi(\pm\cdot)}(R_{\pm}^{\gamma}(\pm\cdot);t,\cdot)$ for all $t\geq0$.
Thus $\phi\geq W_{\pm}^{\gamma}(\pm\cdot)$.

(iii) Fix $\phi\in\mathcal{E}$. By  (ii) and the fact that $$W_{\pm}^{\gamma,M^*}(+\infty)=\frac{\frac{{\rm d}h_{\pm}^{\infty}(0)}{{\rm d}u}-\gamma}{K_+} \mbox{ for all } \gamma\in(0,\min\{ \frac{{\rm d}h_{+}^{\infty}(0)}{{\rm d}u},\frac{{\rm d}h_{-}^{\infty}(0)}{{\rm d}u} \}),
$$
it follows that $\varliminf_{ x\rightarrow\pm\infty }\phi(x)\geq\frac{\frac{{\rm d}h_{\pm}^{\infty}(0)}{{\rm d}u}-\gamma}{K_\pm}\mbox{~~~for~all~}\phi\in\mathcal{E}$.
Letting $\gamma\rightarrow0^+$, we have $\varliminf\limits_{ x\rightarrow\pm\infty }\phi(x)\geq\frac{\frac{{\rm d}h_{\pm}^{\infty}(0)}{{\rm d}u}}{K_\pm}$.

We now claim that $\lim\limits_{x\rightarrow\pm\infty}\phi(x)=u_{\pm}^*$. Otherwise, without loss of generality, we may assume that,
there exists a sequence $x_n\rightarrow+\infty$ such that $0<\lim\limits_{n\rightarrow\infty}\phi(x_{n})\neq u_{+}^{*}$.
In view of \eqref{6.26}, we have $d\phi''(x)+h(x,\phi(x))=0$, and hence, $d\phi''(x+x_n)+h(x+x_n,\phi(x+x_n))=0$.
According the standard elliptic estimates, we have
$\| \phi(\cdot+x_n) \|_{C^{2,\alpha}(I)}\leq C_I$,
where $I$ is any given bounded open interval with $0\in I$ and $C_I$ depend on $I$ and $h$. By the standard diagonal argument, we may assume $\phi(\cdot+x_n)\rightarrow\psi(\cdot)$ in $C_{loc}^{2}(\mathbb{R},\mathbb{R})$
as $n\rightarrow\infty$.
Letting $n\rightarrow\infty$, we get $\psi(0)=\lim\limits_{n\rightarrow\infty}\phi(x_n)\in (0,\infty)\backslash\{u_{+}^{*}\}$, while
$d\psi''(x)+h_{+}^{\infty}(\psi(x))=0\mbox{ for all } x\in\mathbb{R}$, which implies $\psi(x)\equiv u_{+}^{*}$, a contradiction.

(iv) By the proof of Proposition \ref{prop6.19.3} and (i), it follows that
$$u^{M^*}(t,\cdot)\rightarrow\sup\mathcal{E}\in\mathcal{E},~~~\mbox{as~}t\rightarrow\infty.$$
Let $$\mathcal{D}=\{ \phi\in C_{+}:W_{+}^{\gamma^*,M^*}(\cdot)\leq\phi\leq\mathcal{E} \}$$
with $\gamma^*=\frac{1}{2}\min\{ \frac{{\rm d}h_{+}^{\infty}(0)}{{\rm d}u},\frac{{\rm d}h_{-}^{\infty}(0)}{{\rm d}u} \}$.
By virtue of  (ii), we see that $\mathcal{D}\neq\phi$.

By the choices of $W_{+}^{\gamma^*,M^*}$ and $\mathcal{E}$, and Propositions \ref{prop6.19.1}, \ref{prop6.19.2}, it follows that
$u^{\phi}(h;t,\cdot)\in \mathcal{D}$ for all $(t,\phi)\in\mathbb{R}_+\times\mathcal{D}$.
This, together with the compactness of $u^{\phi}(h;t,\cdot)$ and the Schauder fixed point theorem, implies that there exists $\phi^*\in\mathcal{D}$ such that
$u^{\phi^*}(h;t,\cdot)\equiv\phi^*$ for all $t\in\mathbb{R}_+$. Thus, by the definition of $\mathcal{D}$, we have $\phi^*=\inf\mathcal{E}\in\mathcal{E}$.

(v) Fix $\phi\in C_+\backslash\{0\}$ and $\varepsilon>0$.
Take $M=\max\{ \|\phi\|_{L^\infty(\mathbb{R},\mathbb{R})},M^* \}$ and
$\gamma\in\big( 0,\frac{1}{3}\min\{ \frac{{\rm d}h_{+}^{\infty}(0)}{{\rm d}u},\frac{{\rm d}h_{-}^{\infty}(0)}{{\rm d}u} \}  \big)$ with
$$2\sqrt{d(\frac{{\rm d}h_{\pm}^{\infty}(0)}{{\rm d}u}-\gamma)}-\frac{\varepsilon}{3}\geq 2\sqrt{d\frac{{\rm d}h_{\pm}^{\infty}(0)}{{\rm d}u}}-\frac{2\varepsilon}{3}.$$
By Propositions \ref{prop6.20}, \ref{prop6.19.1} and \ref{prop6.19.2},
it follows that
\begin{align}\label{6.26-1}
\lim_{t\rightarrow\infty} &\max \bigg\{   \inf\big\{  |u^{\phi}(t,x)-v|:v\in[\max\{ W_+(x),W_-(-x) \},M^*]  \big\}:\nonumber\\
&\Big[-2\sqrt{d(\frac{{\rm d}h_{-}^{\infty}(0)}{{\rm d}u}-\gamma)}+\frac{\varepsilon}{3}\Big]t\leq x
\leq\Big[2\sqrt{d(\frac{{\rm d}h_{+}^{\infty}(0)}{{\rm d}u}-\gamma)}-\frac{\varepsilon}{3}\Big]t   \bigg\}=0,
\end{align}
where $W_\pm(\cdot):=W_{\pm}^{\gamma,M}$ are defined as in Propositions \ref{prop6.19.2}. In particular, $M^*\geq\omega(\phi)\geq W\pm(-\cdot)$,
where $\omega(\phi)$ is the omega limit set of the orbit $u^{\phi}(t,\cdot)$  with respect to the compact open topology.

 Let
$$\mathcal{D}=\{ \psi\in C_+: \omega(\phi)\geq\psi\geq W_\pm(-\cdot) \}.$$
Since $\omega(\phi)\leq \sup\mathcal{E}$, it follows that
$\mathcal{D}$ is a positively invariant subset of $u^{\phi}(t,\cdot)$, which implies $\mathcal{D}\cap\mathcal{E}\neq\emptyset$ and hence
$\sup\mathcal{E}\geq\omega(\phi)\geq\inf\mathcal{E}$.
Thus, by the definition of $\omega(\phi)$ and the statement (iii), we only need to prove that
$$V_\pm(\delta):=\lim_{\alpha\rightarrow\infty}\Big[  \sup\big\{ |u^{\phi}(t,x)-u_{\pm}^*|:t\geq\alpha  \mbox{~and~}
\alpha\leq\pm x\leq [2\sqrt{d\frac{{\rm d}h_{\pm}^{\infty}(0)}{{\rm d}u}}-\delta]t     \big\}  \Big]=0$$
for all $\delta\in\Big(0,\min\big\{ \sqrt{d\frac{{\rm d}h_{+}^{\infty}(0)}{{\rm d}u}},\sqrt{d\frac{{\rm d}h_{-}^{\infty}(0)}{{\rm d}u}}  \big\}\Big)$.
It suffices to prove $V_+(\delta)=0$ since the same method leads to $V_-(\delta)=0$.
Let $$\overline{P}(\delta)=\limsup_{\alpha\rightarrow\infty}\big[ \sup\{ u^{\phi}(t,x):t\geq\alpha\mbox{~and~}\alpha\leq x\leq
[2\sqrt{d\frac{{\rm d}h_{+}^{\infty}(0)}{{\rm d}u}}-\delta]t \}  \big]$$
and
$$\underline{P}(\delta)=\liminf_{\alpha\rightarrow\infty}\big[ \inf\{ u^{\phi}(t,x):t\geq\alpha\mbox{~and~}\alpha\leq x\leq
[2\sqrt{d\frac{{\rm d}h_{+}^{\infty}(0)}{{\rm d}u}}-\delta]t \}  \big]$$
for all $\delta\in\Big(0,\min\big\{ \sqrt{d\frac{{\rm d}h_{+}^{\infty}(0)}{{\rm d}u}},\sqrt{d\frac{{\rm d}h_{-}^{\infty}(0)}{{\rm d}u}}  \big\}\Big)$.
From Proposition \ref{prop6.20}-(i), we easily see that
$$\overline{P}(\delta)\leq M^*\mbox{ for all }
\delta\in\Big(0, \min\big\{ \sqrt{d\frac{{\rm d}h_{+}^{\infty}(0)}{{\rm d}u}},\sqrt{d\frac{{\rm d}h_{-}^{\infty}(0)}{{\rm d}u}}  \big\}\Big).
$$
In view of  \eqref{6.26-1}, we obtain
$$\underline{P}(\delta)\geq \frac{\frac{{\rm d}h_{+}^{\infty}(0)}{{\rm d}u}-\gamma}{K_+}>0\mbox{ for all }
\delta\in\Big(0,\min\big\{ \sqrt{d\frac{{\rm d}h_{+}^{\infty}(0)}{{\rm d}u}},\sqrt{d\frac{{\rm d}h_{-}^{\infty}(0)}{{\rm d}u}}  \big\}\Big),
$$
where $K_+=K_+^{\gamma,M}$ is defined as in Lemma \ref{lem6.3}.
In view of  the definitions of $\overline{P}(\cdot)$ and $\underline{P}(\cdot)$,
it then follows  that
\begin{enumerate}
\item[(i)]  $\overline{P}(\cdot)$ and $-\underline{P}(\cdot)$ are nonincreasing in
$\Big(0,\min\big\{ \sqrt{d\frac{{\rm d}h_{+}^{\infty}(0)}{{\rm d}u}},\sqrt{d\frac{{\rm d}h_{-}^{\infty}(0)}{{\rm d}u}}  \big\}\Big)$;
\item[(ii)] $M^*\geq \overline{P}(\cdot)\geq \underline{P}(\cdot)\geq \frac{\frac{{\rm d}h_{+}^{\infty}(0)}{{\rm d}u}-\gamma}{K_+}>0   \mbox{ for all }
\delta\in\Big(0,\min\big\{ \sqrt{d\frac{{\rm d}h_{+}^{\infty}(0)}{{\rm d}u}},\sqrt{d\frac{{\rm d}h_{-}^{\infty}(0)}{{\rm d}u}}  \big\}\Big)$.
\end{enumerate}
Now it suffices to prove $\overline{P}(\delta)=\underline{P}(\delta)=u_+^*$ for all
$\delta\in\Big(0,\min\big\{ \sqrt{d\frac{{\rm d}h_{+}^{\infty}(0)}{{\rm d}u}},\sqrt{d\frac{{\rm d}h_{-}^{\infty}(0)}{{\rm d}u}}  \big\}\Big)$.
Otherwise, there exists $\delta_1\in\Big(0,\min\big\{ \sqrt{d\frac{{\rm d}h_{+}^{\infty}(0)}{{\rm d}u}},\sqrt{d\frac{{\rm d}h_{-}^{\infty}(0)}{{\rm d}u}}  \big\}\Big)$
such that either $\overline{P}(\delta)>\underline{P}(\delta)\mbox{~for~all~}\delta\in(0,\delta_1),
$
or $\overline{P}\equiv \underline{P}(\delta)\triangleq P^*\neq u_+^*\mbox{~for~all~}\delta\in(0,\delta_1)$.
By the monotonicity of $\overline{P}(\cdot)$ and $\underline{P}(\cdot)$, there exists $\delta_2\in(0,\delta_1)$ such that
$\overline{P}(\cdot)$ and $\underline{P}(\cdot)$ are continue at $\delta_2$.
Note that
$$\{ u_+^* \}\neq\{ \overline{P}(\delta_2),\underline{P}(\delta_2) \}\subseteq[\frac{\frac{{\rm d}h_{+}^{\infty}(0)}{{\rm d}u}-\gamma}{K_+},M^*].$$
Take $M^{**}=1+M$ and
$\mu=1+\max\left\{  \left|  \frac{{\rm d}h_+^\infty(u)}{{\rm d}u} \right|:u\in[0,M^{**}] \right\}$.
Then $\mu>1$, $\frac{{\rm d}h_+^\infty(u)}{{\rm d}u}+\mu\geq1$ for all $u\in[0,M^{**}]$, and $0<u(t,x)<M^{**}$ for all $(t,x)\in (0,\infty)\times\mathbb{R}$.

It follows from the choices of $\delta_2,\mu$ and the KPP-type property of $h_{+}^{\infty}(\cdot)$  that either
$$\inf h_{+,\mu}^{\infty}\big([\underline{P}(\delta_2),\overline{P}(\delta_2)]\big)>\underline{P}(\delta_2)$$
or
$$\sup h_{+,\mu}^{\infty}\big([\underline{P}(\delta_2),\overline{P}(\delta_2)]\big)<\overline{P}(\delta_2),$$
where $h_{+,\mu}^{\infty}(u)=\frac{1}{\mu}h_+^\infty(u)+u$ for all $u\in \mathbb{R}_+$.
Now we  finish the proof by distinguishing two cases.

{\it Case 1. } $\inf h_{+,\mu}^\infty\big([\underline{P}(\delta_2),\overline{P}(\delta_2)]\big)>\underline{P}(\delta_2)$.

In this case, by the choice of $h_{+,\mu}^{\infty}$ and the fact that $h(s,\cdot)\rightarrow h_{+}^{\infty}$ in $L^{\infty}([0,M^{**}],\mathbb{R})$ as $s\rightarrow\infty$,
there exist $\mathcal{S}_{1}>0$ and $\gamma_{1}\in (0,\frac{\underline{P}(\delta_2)}{2})$ such that
\[
 h^*
\triangleq  \inf \left \{ \frac{h(s,\xi)}{\mu}+\xi:(s,\xi)\in I_{\mathcal{S}_1,\gamma_1}\right\} >    \underline{P}(\delta_2),
\]
where $I_{\mathcal{S}_1,\gamma_1}=(\mathcal{S}_{1},\infty]\times[\underline{P}(\delta_2)-\gamma_{1},\overline{P}(\delta_2)+\gamma_{1}]$.
Fix $\tau\in(0,\delta_{1}-\delta_{2})$. Take $\alpha_{0}>\mathcal{S}_1+1$ and $t_0>\max\{\alpha_0,\frac{\alpha_0}{2\sqrt{d\frac{{\rm d}h_{+}^{\infty}(0)}{{\rm d}u}}-\delta_2}\}$
such that $u^{\phi}(t,x)\in(\underline{P}(\delta_2)-\gamma_{1},\overline{P}(\delta_2)+\gamma_{1})$
when $t\geq t_{0}$ and $\alpha_0\leq x\leq(2\sqrt{d\frac{{\rm d}h_{+}^{\infty}(0)}{{\rm d}u}}-\delta_{2})t$.
It follows from~\eqref{6.26} that for any $t\geq t_{0}$,
\begin{eqnarray*}
	u(t,\cdot)
&\geq &  \int_{t_{0}}^{t}e^{-\mu(t-s)}\int_{\mathbb{R}}\frac{1}{\sqrt{4\pi d(t-s)}}e^{-\frac{(x-y)^{2}}{4\pi d(t-s)}}[\mu u^\phi(s,y)+h(y,u^\phi(s,y))]{\rm d}y{\rm d}s
\\
&\geq & \mbox{$\mu\int_{t_{0}}^{t}e^{-\mu(t-s)}\int_{\alpha_0}^{(2\sqrt{d\frac{{\rm d}h_{+}^{\infty}(0)}{{\rm d}u}}-\delta_{2})s}\frac{1}{\sqrt{4\pi d(t-s)}}
e^{-\frac{(x-y)^{2}}{4\pi d(t-s)}}h^{*}{\rm d}y{\rm d}s$}
\\
&=&\mu h^{*}\int_{0}^{t-t_{0}}e^{-\mu s}\int_{\alpha_0-x}^{(2\sqrt{d\frac{{\rm d}h_{+}^{\infty}(0)}{{\rm d}u}}-\delta_{2})(t-s)-x}\frac{1}{\sqrt{4\pi ds}}e^{-\frac{y^{2}}{4\pi ds}}{\rm d}y{\rm d}s.
\end{eqnarray*}
Moreover, when $\alpha_0<\alpha\leq x\leq(2\sqrt{d\frac{{\rm d}h_{+}^{\infty}(0)}{{\rm d}u}}-\delta_{2}-\tau)t$ and
$t\geq t_{0}:=\frac{2\sqrt{d\frac{{\rm d}h_{+}^{\infty}(0)}{{\rm d}u}}-\delta_{2}-\tau}{2\sqrt{d\frac{{\rm d}h_{+}^{\infty}(0)}{{\rm d}u}}-\delta_{2}-\frac{\tau}{2}}t>\max\{\alpha_0,\frac{\alpha_0}{2\sqrt{d\frac{{\rm d}h_{+}^{\infty}(0)}{{\rm d}u}}-\delta_2}\}$, there holds
\begin{align*}
u(t,x)&\geq \mu h^{*}\int_{0}^{t-t_{0}}e^{-\mu s}\int_{\alpha_0-\alpha}^{\frac{\tau(2\sqrt{d\frac{{\rm d}h_{+}^{\infty}(0)}{{\rm d}u}}-\delta_{2}-\tau)}{2(2\sqrt{d\frac{{\rm d}h_{+}^{\infty}(0)}{{\rm d}u}}-\delta_{2}-\frac{\tau}{2})}t}
\frac{1}{\sqrt{4\pi ds}}e^{-\frac{y^{2}}{4\pi ds}}{\rm d}y{\rm d}s,\\
&\geq \mu h^{*}\int_{0}^{\frac{\tau}{2(2\sqrt{d\frac{{\rm d}h_{+}^{\infty}(0)}{{\rm d}u}}-\delta_{2}-\frac{\tau}{2})}t}
e^{-\mu s}\int_{\alpha_0-\alpha}^{\frac{\tau(2\sqrt{d\frac{{\rm d}h_{+}^{\infty}(0)}{{\rm d}u}}-\delta_{2}-\tau)}{2(2\sqrt{d\frac{{\rm d}h_{+}^{\infty}(0)}{{\rm d}u}}-\delta_{2}-\frac{\tau}{2})}t}
\frac{1}{\sqrt{4\pi ds}}e^{-\frac{y^{2}}{4\pi ds}}{\rm d}y{\rm d}s.
\end{align*}
Note that $\int_{0}^{\frac{\tau}{2(2\sqrt{d\frac{{\rm d}h_{+}^{\infty}(0)}{{\rm d}u}}-\delta_{2}-\frac{\tau}{2})}t}e^{-\mu s}\int_{\alpha_0-\alpha}^{\frac{\tau(2\sqrt{d\frac{{\rm d}h_{+}^{\infty}(0)}{{\rm d}u}}-\delta_{2}-\tau)}{2(2\sqrt{d\frac{{\rm d}h_{+}^{\infty}(0)}{{\rm d}u}}-\delta_{2}-\frac{\tau}{2})}t}
\frac{1}{\sqrt{4\pi ds}}e^{-\frac{y^{2}}{4\pi ds}}{\rm d}y{\rm d}s\rightarrow\frac{1}{\mu}$ as
$t\geq\alpha\rightarrow\infty$, according to the choices of $\delta_{2}$ and $\tau$.
By the definition of $\underline{P}(\delta_2+\tau)$, we have
\[
\underline{P}(\delta_2+\tau)\geq h^{*}>\underline{P}(\tau).
\]
Letting $\tau\rightarrow0^{+}$, we obtain $\underline{P}(\delta_2)>\underline{P}(\delta_2)$,
a contradiction.

{\it Case 2.} $\sup h_{+,\mu}^{\infty}\big( \underline{P}(\delta_2),\overline{P}(\delta_2)] \big)<\overline{P}(\delta_2)$.

In this case, by the definition of $h_{+,\mu}^{\infty}$ and the fact that
$h(s,\cdot)\rightarrow h_{+}^{\infty}$ in $L^{\infty}([0,u^{**}],\mathbb{R})$ as $s\rightarrow\infty$,
there exists $\mathcal{S}_{2}>0$ and $\gamma_{2}\in (0,\frac{\underline{P}(\delta_2)}{3})$ such that
\[ h^{**}
\ \triangleq
\sup\left\{ \frac{h(s,\xi)}{\mu}+\xi:(s,\xi)\in I_{\mathcal{S}_2,\gamma_2}\right\}
 <    \overline{P}(\delta_2),
\]
where $I_{\mathcal{S}_2,\gamma_2}=[\mathcal{S}_{2},\infty)\times[\underline{P}(\delta_2)-\gamma_{2},\overline{P}(\delta_2)+\gamma_{2}]$.
By the definition of $ \overline{P}(\cdot)$ and $\underline{P}(\cdot)$ and their continuity at $\delta_{2}$, it follows  that
there exist $\tau_{0}\in(0,\delta_{2})$, $\alpha_{0}>\mathcal{S}_{2}+1$ and $t_0>\max\{\alpha_0,\frac{\alpha_0}{c-\tau_0}\}$ such that
\[
[\underline{P}(\tau_0),  \overline{P}(\tau_0)]\subseteq(\underline{P}(\delta_2)-\frac{\gamma_{2}}{3},  \overline{P}(\delta_2)+\frac{\gamma_{2}}{3}\big ]
\]
and
$$u^\phi(t,x)\in(\underline{P}(\delta_2)-\gamma_{2},  \overline{P}(\delta_2)+\gamma_{2})$$
whenever $t\geq t_{0}$ and $\alpha_{0}\leq x\leq(2\sqrt{d\frac{{\rm d}h_{+}^{\infty}(0)}{{\rm d}u}}-\tau_{0})t$.
Note that (D3) implies that there exists $H^*>0$ such that $H^*=\sup\{\frac{h(s,u)}{\mu}+ u: (s,u)\in \mathbb{R}\times [0, M^{**}]\}$.

In view of   \eqref{6.26}, we easily see that for any $(t,x)\in\mathbb{R}_{+}\times\mathbb{R}$ with $t\geq t_{0}$,
\begin{eqnarray*}
u^\phi(t,x)&=&
e^{-\mu (t-t_0)}S(t-t_0)[u^\phi(t_0,\cdot)](x)
\mbox{$+\int^{t}_{t_0}
 \frac{e^{-\mu(t-s)}}{\sqrt{4\pi d (t-s)}}\int_{(2\sqrt{d\frac{{\rm d}h_{+}^{\infty}(0)}{{\rm d}u}}-\tau_{0})s}^\infty
[h(y,u^\phi(s,y))+\mu u^\phi(s,y)]$}
\\
&&\times e^{-\frac{(x-y)^2}{4d (t-s)}}{\rm{d}}y \mathrm {d} s
\mbox{$+\int^{t}_{t_0}
 \frac{e^{-\mu(t-s)}}{\sqrt{4\pi d (t-s)}}\int^{(2\sqrt{d\frac{{\rm d}h_{+}^{\infty}(0)}{{\rm d}u}}-\tau_{0})s}_{\alpha_0}
[h(y,u^\phi(s,y))+\mu u^\phi(s,y)]e^{-\frac{(x-y)^2}{4d (t-s)}}{\rm{d}}y \mathrm {d} s$}
\\
&&\mbox{$+\int^{t}_{t_0}
 \frac{e^{-\mu(t-s)}}{\sqrt{4\pi d (t-s)}}\int^{\alpha_0}_{-\infty}
[h(y,u^\phi(s,y))+\mu u^\phi(s,y)]e^{-\frac{(x-y)^2}{4d (t-s)}}{\rm{d}}y \mathrm {d} s$}
\\
&\leq&
M^{**}e^{-\mu (t-t_0)}
+\mu H^{**}\int^{t}_{t_0}
 \frac{e^{-\mu(t-s)}}{\sqrt{4\pi d (t-s)}}[\int_{(2\sqrt{d\frac{{\rm d}h_{+}^{\infty}(0)}{{\rm d}u}}-\tau_{0})s}^\infty
e^{-\frac{(x-y)^2}{4d (t-s)}}{\rm{d}}y+\int_{-\infty}^{\alpha_0}
e^{-\frac{(x-y)^2}{4d (t-s)}}{\rm{d}}y ]\mathrm {d} s
\\
&&\mbox{$+\int^{t}_{t_{0}}
 \frac{e^{-\mu(t-s)}}{\sqrt{4\pi d (t-s)}}\int^{(2\sqrt{d\frac{{\rm d}h_{+}^{\infty}(0)}{{\rm d}u}}-\tau_{0})s}_{\alpha_0}
[h(y,u^\phi(s,y))+\mu u^\phi(s,y)]e^{-\frac{(x-y)^2}{4d (t-s)}}{\rm{d}}y \mathrm {d} s$}.
\end{eqnarray*}
For any $(t,x)\in \mathbb{R}_+\times   \mathbb{R}$   with  $t>t_0:=\frac{2\sqrt{d\frac{{\rm d}h_{+}^{\infty}(0)}{{\rm d}u}}-\frac{\tau_0+\delta_2}{2}}{2\sqrt{d\frac{{\rm d}h_{+}^{\infty}(0)}{{\rm d}u}}-\tau_0}t>\max\{\alpha_0,\frac{\alpha_{0}}{\delta_2-\tau_{0}}\} $ and $\mathcal{S}_2+1<\alpha_0<\alpha\leq x\leq t(2\sqrt{d\frac{{\rm d}h_{+}^{\infty}(0)}{{\rm d}u}}-\delta_2)$, we have
\begin{eqnarray*}
u^\phi(t,x)&\leq &
M^{**}e^{-\mu (t-t_0)}
+\mu H^{**}\int^{t-t_0}_{0}
 \frac{e^{-\mu s}}{\sqrt{4\pi d s}}[\int_{(2\sqrt{d\frac{{\rm d}h_{+}^{\infty}(0)}{{\rm d}u}}-\tau_{0})(t-s)}^\infty
e^{-\frac{(x-y)^2}{4d s}}{\rm{d}}y+\int_{-\infty}^{\alpha_0}
e^{-\frac{(x-y)^2}{4d s}}{\rm{d}}y ]\mathrm {d} s
\\
&&\mbox{$+\int^{t-t_0}_{0}
 \frac{e^{-\mu s}}{\sqrt{4\pi d s}}\int^{(2\sqrt{d\frac{{\rm d}h_{+}^{\infty}(0)}{{\rm d}u}}-\tau_{0})(t-s)}_{\alpha_0}
[h(y,u^\phi(t-s,y))+\mu u^\phi(t-s,y)]e^{-\frac{(x-y)^2}{4d s}}{\rm{d}}y \mathrm {d} s$}
\\
&\leq &
M^{**}e^{-\mu (t-t_0)}
+\mu H^{**}\int^{t-t_0}_{0}
 \frac{e^{-\mu s}}{\sqrt{4\pi d s}}\int_{(2\sqrt{d\frac{{\rm d}h_{+}^{\infty}(0)}{{\rm d}u}}-\tau_{0})(t-s)}^\infty
e^{-\frac{[\frac{(\delta_2-\tau_0)t}{2}]^2}{8dt}}e^{-\frac{(x-y)^2}{8ds}}{\rm{d}}y \mathrm {d} s
\\
&&
+\mu H^{**}\int^{t-t_0}_{0}
 \frac{e^{-\mu s}}{\sqrt{4\pi d s}}\int_{-\infty}^{\alpha_0-\alpha}
e^{-\frac{y^2}{4d s}}{\rm{d}}y \mathrm {d} s
\mbox{$+\mu h^{**}\int^{t-t_0}_{0}
 \frac{e^{-\mu s}}{\sqrt{4\pi d s}}\int^{(2\sqrt{d\frac{{\rm d}h_{+}^{\infty}(0)}{{\rm d}u}}-\tau_{0})(t-s)}_{\alpha_0}
e^{-\frac{(x-y)^2}{4d s}}{\rm{d}}y \mathrm {d} s$}
\\
&\leq &
M^{**}e^{-\frac{\delta_2-\tau_0}{2(2\sqrt{d\frac{{\rm d}h_{+}^{\infty}(0)}{{\rm d}u}}-\tau_0)}t}
+2 H^{**} e^{-\frac{(\delta_2-\tau_0)^2t}{32
d}}
+\frac{ H^{**} e^{\sqrt{\frac{\mu}{d}}(\alpha_0-\alpha)}}{2}
+ h^{**}
\end{eqnarray*}
Letting $t\geq \alpha\to \infty$, we get $\overline{P}(\delta_2)\leq \emph {h}^{~**}<\overline{P}(\delta_2)$, a contradiction. \qed

\begin{thm}\label{6.6.1}
Assume that  $h(\cdot,\cdot)\neq h(0,\cdot)$ and $h(s,\alpha u)\geq\alpha h(s,u)$ for all $(s,u,\alpha)\in\mathbb{R}\times (0,M^{*})\times(0,1)$.
Then the following statements are valid:
\begin{itemize}
\item [{\rm (i)}]  \eqref{6.26} has a unique steady state $W\in C_{+}^{\circ}$ with $W(\pm\infty)=u_{\pm}^*$;

\item [{\rm (ii)}] If $\phi\in C_{+}\backslash \{0\}$ and $\varepsilon>0$, then
$$\lim_{t\rightarrow\infty}\Big[\max\big\{  |u^{\phi}(h;t,x)-W(x)|:~(-c_{-}^{*}+\varepsilon)t\leq x\leq(c_{+}^{*}+\varepsilon)t \big\}\Big]=0;$$

\item [{\rm (iii)}] If $\phi$ has a compact support, then
$$\lim_{t\rightarrow\infty}\Big[\sup\big\{ u^{\phi}(h;t,x):~\pm x\geq(c_{\pm}^{*}+\varepsilon)t \big\}\Big]=0.$$
\end{itemize}
\end{thm}

\noindent {\bf Proof.} (i)  In view of Proposition \ref{prop6.19.4},
we only need to prove the uniqueness. Otherwise,
$$\widetilde{W}:=\inf\mathcal{E}\lneq W:=\sup\mathcal{E}.$$
Let $$\alpha^*=\inf\{ \alpha\geq1:~W\leq\alpha \widetilde{W} \}.$$
Then $\alpha^*>1$ and $W\leq\alpha^*\widetilde{W}$.
By Proposition \ref{prop6.19.4}-(iii), there exist $s_{0}>0$ and $\delta_0\in(0,\alpha^*)$ such that
$$W<\alpha^*\widetilde{W}\mbox{~~and~~}W(x)\leq(\alpha^*-\delta_0)\widetilde{W}(x)$$
for all $|x|\geq s_0$. By virtue of Proposition \ref{prop6.20}-(ii), we have
$$\widetilde{W}(x)=u^{\widetilde{W}}(h;t,x)>u^{\frac{W}{\alpha^*}}(h;t,x)\geq\frac{1}{\alpha^*}u^{W}(h;t,x)=\frac{W(x)}{\alpha^*}$$
for all $x\in\mathbb{R}$ and $t\in(0,\infty)$.
Thus, there exists $\delta_1>0$ such that $(\alpha^*-\delta_1)\widetilde{W}(x)\geq W(x)$ for all  $|x|\leq s_0$.
Let $\delta=\min\{ \delta_0,\delta_1 \}$. It then follows that $(\alpha^*-\delta)\widetilde{W}(x)\geq W(x) \mbox{ for all } x\in\mathbb{R}$,
which contradicts the choice of $\alpha^*$.

(ii) follows from (i) and Proposition \ref{prop6.19.4}-(v), and (iii) directly follows from Proposition \ref{prop6.19}.  \qed

\

We should point out that the positive stationary solutions and spatial spreading
speeds were studied in \cite{ks2011} for the KPP-type evolution equations in a locally spatially inhomogeneous media. Our Theorem \ref{6.6.1} extends such results to the one-dimensional case of asymptotically inhomogeneous media. We also note that the method in this subsection may be used to
remove  the monotonicity condition on functions $f(s,u)$ and $g(s,y,u)$ with respect to  $s$, as assumed in (B3) and (C3).

\

\noindent
{\bf Acknowledgements.} T.  Yi's research is supported by the National Natural Science Foundation of  China (Grant No. 11971494),  and X.-Q. Zhao's research is
supported in part by the NSERC of Canada.
This work was initiated  during Dr.\ Yi's visit to Memorial University of Newfoundland, and he would like to thank the Department of Mathematics
and Statistics there for its kind hospitality.

\end{document}